\documentclass[[preprint,12pt]{elsarticle}
\usepackage{setspace}
\usepackage{float}
\usepackage{cite}
\usepackage{xspace}
\usepackage[toc,page]{appendix}
\usepackage{nomencl}
\usepackage{natbib}
\usepackage{fullpage}

\usepackage{amsmath}
\usepackage{amssymb}
\usepackage{latexsym}
\usepackage{esint}

\usepackage{boxedminipage}

\usepackage{listings}

\usepackage{minitoc}
\usepackage[mathlines]{lineno} % Control line numbers
\usepackage{ifpdf}
\usepackage{xspace}
\usepackage{bm}
\usepackage{natbib}

\usepackage{boxedminipage}

\usepackage{graphicx}

\usepackage{latexsym,amssymb,amsmath,color}

\newcommand{\xxi}{\vec{\xi}}

\newcommand{\norm}[1]{\left\| {#1} \right\|_{\scriptscriptstyle L^{2}(\Omega,D)}}

\renewcommand{\vec}[1]{{\mathchoice
                     {\mbox{\boldmath$\displaystyle{#1}$}}
                     {\mbox{\boldmath$\textstyle{#1}$}}
                     {\mbox{\boldmath$\scriptstyle{#1}$}}
                     {\mbox{\boldmath$\scriptscriptstyle{#1}$}}}}
                    
\newcommand{\Loeve}{Lo\`{e}̀ve\xspace}
\newcommand{\Nobs}{{N_o}}

\journal{Computer Methods in Applied Mechanics and Engineering}
\begin{document}

\ifpdf
\DeclareGraphicsExtensions{.pdf, .png, .jpg, .tif}
\else
\DeclareGraphicsExtensions{.png, .jpg, .tif, .eps}
\fi

\begin{frontmatter}

\title{Coordinate Transformation and Polynomial Chaos for the Bayesian Inference of a Gaussian Process with Parametrized Prior Covariance Function}
%\date{\today}

\author[kaust1]{Ihab Sraj}
\corref{cor1}
\ead{ihab.sraj@kaust.edu.sa}
\author[CNRS]{Olivier P. Le {Ma\^{i}tre}}

\author[duke,kaust2]{Omar M. Knio}
\author[kaust1,kaust2]{Ibrahim Hoteit}
\address[kaust1]{Division of Physical Sciences and Engineering, King Abdullah University of Science and Technology, Thuwal, Saudi Arabia}
\address[CNRS]{LIMSI-CNRS, BP 133, Bt 508, 91403 Orsay Cedex, France}
\address[duke]{Department of Mechanical Engineering and Materials Science, Duke University, 144
Hudson Hall, Durham, North Carolina 27708, USA}
\address[kaust2]{Division of Computer, Electrical and Mathematical Sciences and Engineering, King Abdullah University of Science and Technology, Thuwal, Saudi Arabia}

\cortext[cor1]{Corresponding author}

%\date{\today}

\begin{abstract}
This paper addresses model dimensionality reduction for Bayesian inference based on prior Gaussian fields
with uncertainty in the covariance function hyper-parameters. 
The dimensionality reduction is traditionally achieved using the Karhunen-\Loeve expansion of a prior Gaussian process 
assuming covariance function with fixed hyper-parameters, despite the fact that these are uncertain in nature.
The posterior distribution of the Karhunen-\Loeve coordinates is then inferred using available observations. 
The resulting inferred field is therefore dependent on the assumed hyper-parameters. 
Here, we seek to efficiently estimate both the field and covariance hyper-parameters using
Bayesian inference. To this end, a generalized Karhunen-\Loeve expansion is derived using a coordinate transformation to account for the dependence with respect to the covariance hyper-parameters.
Polynomial Chaos expansions are employed for the acceleration of the Bayesian inference using similar coordinate transformations, 
enabling us to avoid expanding explicitly the solution dependence on the uncertain hyper-parameters. 
We demonstrate the feasibility of the proposed method on a transient diffusion equation 
by inferring spatially-varying log-diffusivity fields from noisy data. The inferred profiles 
were found closer to the true profiles when including the hyper-parameters' uncertainty
in the inference formulation.
\end{abstract}

\begin{keyword}
 Karhunen-\Loeve expansion \sep dimensionality reduction \sep Markov Chain Monte Carlo \sep polynomial chaos \sep Bayesian inference
\end{keyword}

\end{frontmatter}

%\linenumbers
%\maketitle

%\tableofcontents
%\printnomenclature

%Department of Physical Sciences and Engineering, King Abdullah University for Science and Technology, Thuwal, Saudi Arabia
%\author[kaust]{}
%
%Department of Physical Sciences and Engineering, King Abdullah University for Science and Technology, Thuwal, Saudi Arabia
%\doublespacing

\section{Introduction}
\label{sec:intro}
Inverse problems arise in many applications whenever we seek to find some information 
about a physical system based on some observations. From a computational point of view, 
a major challenge of inverse problems is their ill-posedness where there 
is no guarantee that a solution exists, multiple solutions may exist, or even the 
solution does not depend continuously on the observations. 
This can be significantly affected by measurement errors,
and inferring a suitable solution from noisy observations is an important and challenging topic.       

In this paper, we are only concerned with Bayesian approaches to inverse problems.
This is motivated by their ability of providing complete posterior statistics
and not just a single value for the quantity of interest. The multi-dimensional posterior 
can be directly explored via Markov Chain Monte Carlo (MCMC). This, however, requires repeated 
simulations (sometimes hundreds of thousands) of the forward model, once for every proposed set 
of parameters of the Markov chain~\citep{Malinverno2002}. This practice
renders Bayesian methods computationally prohibitive for large-scale applications.
Acceleration techniques have been proposed in the literature
in which a surrogate model is constructed that requires a much smaller ensemble of 
forward model runs which is then used in the sampling MCMC step instead at a significantly reduced 
computational cost. Marzouk  \textit{et al.}~\citep{MarzoukNajmRahn:2007} for instance proposed
a spectral projection method that uses spectral expansion of the prior model in Polynomial 
Chaos (PC) basis. The PC method has been extensively investigated in the literature, and its suitability  
for large-scale models has been demonstrated in various settings, including ocean~\citep{winokur:2012,sraj:2013a,sraj:2013b,MatternFennelDowd2012}, 
tsunami~\citep{sraj:2014,Cheung2011}, climate modeling~\citep{sargsyan2011} and subsurface flow modeling~\citep{Elsheikh2014}.

The PC method has been shown to be efficient for inverse problems involving a limited 
number of stochastic parameters; yet in some cases the unknown quantity is a spatial or temporal field
in which the number of stochastic parameters is quite large. Computational challenges in this case arise in the surrogate model construction as PC suffers from the curse of dimensionality~\citep{LeMaitreKnio2010}. In addition, convergence is hard to achieve using the Bayesian inference
due to the high dimensionality of the posterior. To overcome this numerical issue, Marzouk  \textit{et al.}~\citep{MarzoukNajm2009} 
introduced truncated Karhunen-\Loeve (KL) expansions to parametrize the stochastic field, endowed with a hierarchical Gaussian process prior.
The idea is to transform the high-dimensional stochastic forward problem into a smaller problem
whose solution captures that of the deterministic forward model over the support 
of the prior. Galerkin projection on a PC basis was used to seek the solution of the problem,
and a reduced-dimensionality surrogate posterior density was constructed that is inexpensive to evaluate.

The Gaussian process prior assumed in Marzouk  \textit{et al.}~\citep{MarzoukNajm2009} is associated with
hyper-parameters that are rarely known in practice. Assuming otherwise renders the quantification of prior uncertainty unrealistic
and incomplete. Hierarchical Bayesian inference is proposed in the literature for calibration in presence of uncertain hyper-parameters
but is done \textit{a priori}~\citep{Rasmussen:2005}. The method proposed by Marzouk \textit{et al.}~\citep{MarzoukNajm2009}
does not explicitly consider the effect of length-scales and only includes one hyper-parameter accounting for prior variance. 
An attempt to extend the method proposed by Marzouk  \textit{et al.}~\citep{MarzoukNajm2009} 
for priors with uncertain hyper-parameters has been recently proposed by Tagade and Choi~\citep{Tagade:2014}. 
In their work, a methodology is introduced to obtain a KL expansion of a stochastic process in terms of functions of the
hyper-parameters. The prior uncertainty in these hyper-parameters was expanded in a PC basis, and Galerkin projection was used to evaluate PC coefficients of the surrogate model. 
The hyper-parameters hence become part of the inference problem and are estimated from the observations. 

This paper proposes an extension of the method of Marzouk  \textit{et al.}~\citep{MarzoukNajm2009} that is also
an alternative to Tagade and Choi method~\citep{Tagade:2014}. Our proposed method explores the origin of the KL expansion
where it is based on the eigen-functions and eigen-values of a given covariance function. These 
eigen-functions form a basis in a space dictated by the covariance hyper-parameters. 
Our method utilizes a change of basis methodology and therefore
transforming the KL expansion based on one certain set of hyper-parameters into another.
A fundamental distinction of the present work is that we avoid constructing a PC 
expansion for the uncertain hyper-parameters, and instead use the PC expansion constructed for a reference 
set of hyper-parameters and apply transformations to obtain PC expansion for any another set of hyper-parameters.
The advantage of the proposed method is that the dimensionality of the PC expansion is not augmented by the number of hyper-parameters
of the covariance function. Also, our method avoids cases when the hyper-parameters have complex distributions
and PC bases may not even exist. 

To outline the proposed developments, we start in Section~\ref{sec:inference} by providing a statistical 
formulation of the inverse problem based on Bayesian inference.  
Section~\ref{sec:KL} then presents the KL expansion and its generalization to account for uncertain hyper-parameters by means of change of basis.
Section~\ref{sec:pc} describes the role of PC in Bayesian inference acceleration. 
In Section~\ref{sec:results}, numerical results for the calibration of a one-dimensional 
toy problem are presented and Section~\ref{sec:conc} concludes the paper with a summary of the results, discussion and conclusion.

%%%%%%%%%%%%%%%%%%%%%%%%%%%%%%%%%%%%%%%%%%%%%%%%%%%%%%%%%%%%%%%%%%%%%%%%%%%%%%%%%%%%%%%%%%%%%%%%%%%%%%%%%%%
%%%%%%%%%%%%%%%%%%%%%%%%%%%%%%%%%%%%%%%%%%%%%%%%%%%%%%%%%%%%%%%%%%%%%%%%%%%%%%%%%%%%%%%%%%%%%%%%%%%%%%%%%%%%%%%%
%%%%%%%%%%%%%%%%%%%%%%%%%%%%%%%%%%%%%%%%%%%%%%%%%%%%%%%%%%%%%%%%%%%%%%%%%%%%%%%%%

\section{Bayesian Inference}
\label{sec:inference}
Bayesian inference is a statistical approach to inverse problems that has gained much interest in different applications including
ocean~\citep{Alexanderian2011a,Zedler2012,sraj:2013a}, climate~\citep{OlsonEtAl2012}  and geophysical~\citep{Malinverno2002} modeling. 
We review the Bayesian approach briefly  below and discuss its implementation to our problem.

Our objective is to infer a deterministic field $m(\vec{x})$, for some $\vec x\in D$, from a finite set of $\Nobs \ge 1$ observations $\vec d \in \mathbb R^\Nobs$.
We consider situations where the observations $\vec d$ are not direct measurements of $m(\vec x)$, but are derived quantities that can be predicted using a 
model-problem (typically a set of partial differential equations), often called the forward model, relating the $m(\vec x)$ to the model predictions: 
$m(\vec x) \mapsto \vec u(m) \in \mathbb R^\Nobs$.
The Bayesian formula updates our prior knowledge of the $m$ introducing an error model for the discrepancy between the model predictions $\vec u(m)$ and the observations 
$\vec d$; the Bayes'~rule is expressed as~\citep{sivia}:
\begin{equation}
 	p(m,\sigma_o^2|\vec d) \propto p(\vec d | m,\sigma_o^2) p_m(m) p_o(\sigma_o^2),
\end{equation}
where $p(\vec d|m,\sigma_o^2)$ is the likelihood of the observations, given $m$ and $\sigma_o^2$ the error model 
hyper-parameter with prior $p_o(\sigma_o^2)$, and $p_m(m)$ is the field's prior. 
For simplicity, an unbiased additive Gaussian error model will be considered, 
\begin{equation}
\vec \epsilon \doteq \vec d - \vec u(m),\quad \vec \epsilon \sim {\cal N}(0,\sigma_o^2 I_\Nobs),
\end{equation}
where $N(0,\sigma_o^2I_\Nobs)$ denotes the centered multivariate Gaussian distribution with diagonal covariance $\sigma_o^2 I_\Nobs$. 
In other words, the errors in the observations are assumed independent. For this choice, the likelihood becomes
\begin{equation}
	p(\vec d|m,\sigma_o^2) = \prod_{i=1}^\Nobs p_\epsilon ( d_i - u_i (m),\sigma_o^2), \quad p_\epsilon ( x,\sigma_o^2) \doteq
	\frac{1}{\sqrt{2\pi\sigma_o^2}} \exp \left[-\frac{x^2}{2\sigma_o^2}\right]. \label{eq:likelihood1}
\end{equation}

The main difficulties with the posterior above are the infinite dimensional character of $m(\vec x)$ and its prior definition. 
A discretization of $m(\vec x)$ is needed to perform the inference and setting a finite dimensional prior distribution. 
If $m(\vec x)$ is endowed with a Gaussian prior, it is fully characterized by its second-order properties, namely its mean $\mu(\vec x)$ and covariance function $\cal C(\vec x,\vec x')$.
From $\mu$ and $\cal C$, one can rely on truncated Karhunen-\Loeve (KL) decomposition to represent $m(\vec x)$ as a convergent series involving a finite set of KL coordinates (or expansion coefficients) $\eta_k$, $k=1,\ldots,K$ as discussed in Section~\ref{sec:KL}.
The inference problem can then be reformulated for the vector $\vec \eta$ of coordinates $\eta_k$, leading to
\begin{equation}
 	p(\vec \eta,\sigma_o^2|\vec d) \propto p(\vec d|\vec \eta,\sigma_o^2) p_\eta(\vec \eta) p_o(\sigma_o^2),
\end{equation}
where $p_\eta(\vec \eta) = \exp( - \vec \eta^T \vec \eta / 2) / (2\pi)^{K/2}$ is the Gaussian prior of the KL coordinates.

As discussed below, the covariance function $\cal C(\vec x,\vec x')$ is generally selected on the basis of limited knowledge 
and the inference of $m(\vec x)$ can be improved by introducing additional hyper-parameters $\vec q$ in the definition of $\cal C$
i.e. $\cal C(\vec x,\vec x', \vec q)$.
This yields the generalized Bayes' formula,
\begin{equation}
 	p(\vec \eta, \vec q, \sigma_o^2|\vec d) \propto p(\vec d|\vec \eta,\vec q ,\sigma_o^2) p_\eta(\vec \eta)p_q(\vec q) p_o(\sigma_o^2),
\end{equation}
where $p_q$ is the prior distribution of the covariance parameters.
For the case of covariance with hyper-parameters $\vec q$, the likelihood takes the following general form,
\begin{equation}
	p(\vec d|\vec \eta, \vec q,\sigma_o^2) = \prod_{i=1}^\Nobs p_\epsilon ( d_i - u_i (\vec \eta,\vec q),\sigma_o^2), 
\end{equation}
with $p_\epsilon$ defined in Eq.~\eqref{eq:likelihood1}, and $u_i(\vec \eta,\vec q)$ being a short-hand notation for the model prediction 
$u_i(m)$ with $m(\vec \eta, \vec q)$ the reconstructed field.

\begin{figure}[htb]
\centering
\includegraphics[width=0.75\textwidth]{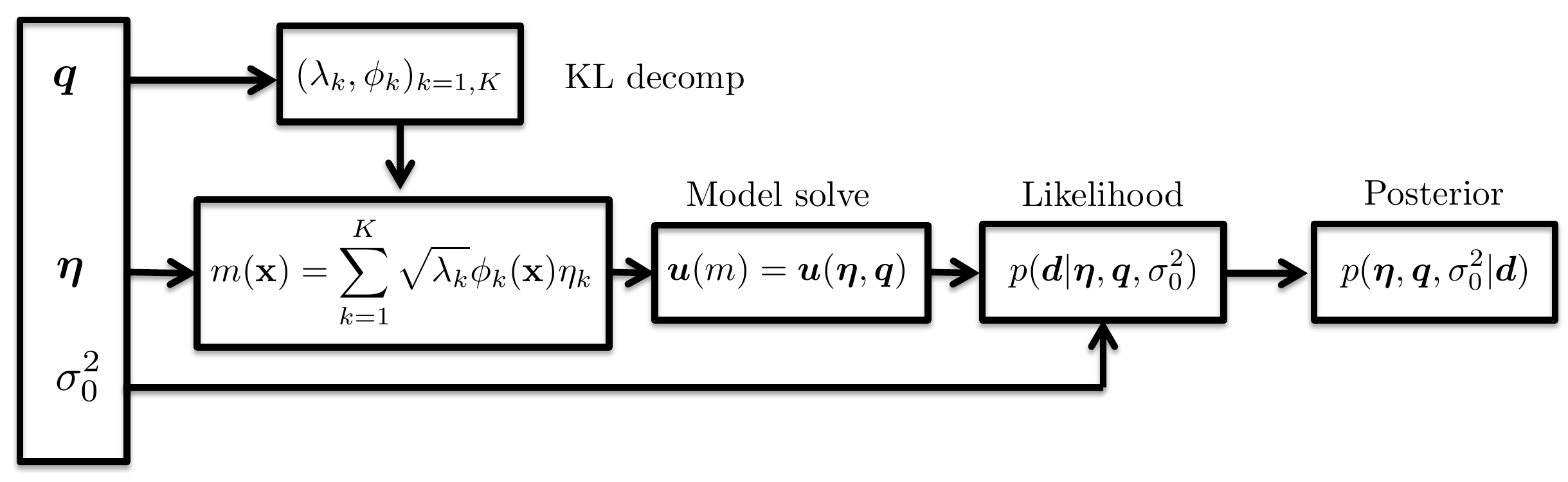}
\caption{Flow-chart for the evaluation of the posterior distribution in the inference problem.}
\label{fig:flowchart1}
\end{figure}

Inferring the field then amounts to sampling the posterior of KL coordinates $\vec \eta$ and hyper-parameters $\vec q$. 
In general, the sample space is high-dimensional and suitable computational strategy is the Markov chain Monte Carlo (MCMC) method. 
In this work, we rely on an adaptive Metropolis-Hastings MCMC algorithm~\citep{Haario2001,Gareth2009} to accurately and efficiently sample the posterior 
distribution $p(\vec \eta, \vec q, \sigma_o^2|\vec d)$.
This requires the evaluation of the posterior (up to its normalization constant) for multiple sample values of $(\vec \eta, \vec q, \sigma_o^2)$.
The computational flow-chart for an evaluation of the posterior is schematically illustrated in Figure~\ref{fig:flowchart1}. Briefly, given a sample value of $\vec q$, the dominant KL modes of ${\cal C}(\vec x,\vec x',\vec q)$ are computed and the corresponding field $m(\vec x)$ is constructed using the sampled value of $\vec \eta$. This field is fed into the solver to compute the model predictions $\vec u (\vec \eta,\vec q)$ which are used, together with the sample value of the model error parameter $\sigma_o^2$, to successively compute the likelihood and finally the posterior.

In general, the most computationally demanding part for sampling the posterior is the computation of the model predictions, given 
$(\vec \eta,  \vec q)$ (that is the realization of the field). 
This is particularly the case when the predictions involve the solution of partial differential equations.
This computational cost motivates the substitution of $\vec u(\vec \eta,\vec q)$ with a polynomial surrogate model $\tilde{\vec u}(\vec \eta,\vec q)$, whose evaluation is inexpensive compared to the solution of the complete model. The surrogate is constructed offline and subsequently used on-line when running the MCMC algorithm. Specifically, the likelihood of the observations is approximated using
\begin{equation}
	p(\vec d | \vec \eta,\vec q,\sigma_o^2) = \prod_{i=1}^\Nobs p_\epsilon ( d_i - u_i (\vec \eta,\vec q),\sigma_o^2) \approx
	\prod_{i=1}^\Nobs p_\epsilon ( d_i - \tilde u_i (\vec \xi(\vec \eta,\vec q)),\sigma_o^2),
\end{equation}
where, as mentioned previously, $\tilde {\vec u}(\vec \xi)$ is a polynomial and 
$\xxi : (\vec \eta, \vec q) \mapsto \xxi(\vec \eta,\vec q)$ is an explicit change of coordinates. 
Construction of the surrogate model for the predictions is detailed in the next two sections; Section~\ref{sec:KL} introduces the 
$\vec q$-dependent coordinate transformation, while the polynomial approximation $\tilde{\vec u}(\xxi)$ is discussed in 
Section~\ref{sec:pc}, together with the resulting surrogate-based sampling scheme.
%%%%%%%%%%%%%%%%%%%%%%%%%%%%%%%%%%%%%%%%%%%%%%%%%%%%%%%%%%%%%%%%%%%%%%%%%%%%%%%%%%%%%
%%%%%%%%%%%%%%%%%%%%%%%%%%%%%%%%%%%%%%%%%%%%%%%%%%%%%%%%%%%%%%%%%%%%%%%%%%%%%%%%%%
%%%%%%%%%%%%%%%%%%%%%%%%%%%%%%%%%%%%%%%%%%%%%%%%%%%%%%%%%%%%%%%%%%

\section{Coordinate transformation for Uncertain Covariance Functions}
\label{sec:KL}

\subsection{Karhunen-\Loeve expansion}

Let $D\subset \mathbb R^d$, $d\ge 1$, be a bounded domain, and denote $X\doteq L^2(D)$ equipped with inner product $(\cdot,\cdot)_X$ and norm $\|\cdot\|_X$:
\begin{equation}
 u\in X \Leftrightarrow \| u\|_X < \infty, \quad \|u\|_X^2 = (u,u)_X = \int_D |u(\vec x)|^2 d\vec x.
 \end{equation}

Consider a real-valued stochastic process $M(\vec x,\omega)$ with mean $\mu(\vec x)$ and continuous covariance function $\cal C(\vec{x},\vec{x}')$ on $D \times D$;
$\omega$ is a random event belonging to a sample space $\Omega$ of a probability space $(\Omega,\Sigma,P)$. 
The covariance function is defined as
\begin{equation}
	{\cal C}(\vec{x},\vec{x}') = \mathbb E\left[(M(\vec x,\cdot) - \mu(\vec x)) (M(\vec x',\cdot) - \mu(\vec x')) \right],
\end{equation}
where $\mathbb E$ denotes the expectation operator. The covariance function $\cal C$ is symmetric positive semi-definite and thus by Mercer's theorem~\citep{Grigoriu2002} it has the following spectral decomposition:
\begin{equation}
	{\cal C}(\vec{x},\vec{x}') =\sum_{k=1}^{\infty}\lambda_k \phi_k(\vec{x}) \phi_k(\vec{x}'),\label{eq:covdec}
\end{equation}
where the $\lambda_k$ and $\phi_k(\vec{x})$ are the eigen-values and associated (normalized) eigen-functions of the linear operator corresponding to the covariance function $\cal C$; they satisfy the Fredholm equation of the second kind:
\begin{equation}
\int_D {\cal C} (\vec{x},\vec{x}') \phi_k(\vec{x}') dx = \lambda_k \phi_k(\vec{x}), \quad \|\phi_k \|_X = 1.
\label{eq:fred}
\end{equation}
The eigen-values $\lambda_k$ are real and countable and the eigen-functions $\phi_k(\vec x)$ are continuous and constitute an orthonormal basis in $L^2(D)$.
Ordering the eigen-values in a decreasing sequence $\lambda_1 \ge \lambda_2 \ge \cdots \ge 0$, the truncated Karhunen-\Loeve (KL) expansion $M_K$ of $M$ is given by~\citep{GhanemSpanos1991}
\begin{equation}
M(\vec{x},\omega) \approx M_{K}(\vec{x},\omega) \doteq \mu(\vec x) + \sum_{k = 1}^{K}\sqrt{\lambda_k}\phi_k(\vec{x})\eta_k(\omega),
\label{eq:kldec}
\end{equation}
where $K$ is the number of expansion terms retained in the spectral approximation. The stochastic coefficients
\begin{equation}
	\eta_k(\omega) = (M(\vec x,\omega) - \mu(\vec x), \phi_k(\vec x))_X, \label{def:projcoeff}
\end{equation}
are mutually uncorrelated random variables with zero mean and unit variance, such that 
$\mathbb E\left[ \eta_k \eta_{k'}\right] = \delta_{kk'}$.
Under the assumption that $M$ is a Gaussian Process ($\cal{GP}$) denoted by $M \sim \cal{GP}$ $(\mu,{\cal C})$, the $\eta_k$'s are Gaussian and also independent. 
The truncated KL expansion is optimal in the mean square sense, meaning that of all possible $K$-term expansions, the $M_K$ in 
Eq.~\eqref{eq:kldec} with $\lambda_k$ and $\phi_k(\vec x)$ satisfying Eq.~\eqref{eq:fred} minimizes the mean-squared error in  the approximation of $M$~\citep{GhanemSpanos1991}. 
While it is known that the KL decomposition of $M$ converges uniformly as $K\to \infty$~\citep{adler2007}, the truncation error has implicit dependence on the covariance function $\cal C$.

The KL expansion is often employed to reduce the dimensionality in inverse problems, considering the expansion coefficients $\eta_{k=1,\dots,K}$ in Eq.~\eqref{eq:kldec} as reduced coordinates for the field $m(\vec x)$ to be inferred from the collected observations~\citep{Elsheikh2014a}.
%Specifically, one updates the posterior distribution of the expansion coefficients $\eta_{k=1,\dots,K}$ from the collected observations~\citep{Elsheikh2014a}, instead of updating the whole infinite dimensional process.
%Then, Eq.~\eqref{eq:kldec}, together with the posterior distribution for the $\eta_k$, approximates the posterior of $M$ and so $Q$ (defined as the Maximum A Posteriori Probability of $M$, its expectation, median, \dots). 
In the Bayesian framework, this amounts to the determination of the posterior distribution of the expansion coefficients vector $\vec \eta$, which can be sampled or analyzed to estimate the characteristics of the field $m(x)$ (in particular, retrieving the median, MAP value, confidence intervals,\dots). 
%  where $\eta_{k=1,\dots,K}$ are inferred from collected observations, 
% instead of inferring directly the infinite dimensional $M$~\citep{Elsheikh2014a}. 
% In this case, the inferred $M_K$ has 
However, the posterior and so the inferred field $m$ have implicit dependencies on the assumed prior covariance structure. 
This point has motivated the introduction of parametrized covariance families, as discussed in the following section, where the covariance parameters are treated as hyper-parameters in the inference procedure.

%%%%%%%%%%%%%%%%%%%%%%%%%%%%%%%%%%%%%%%%%%%%%%%%%%%%%%%%%
%%%%%%%%%%%%%%%%%%%%%%%%%%%%%%%%%%%%%%%%%%%%%%%%%%%%%%%%%%
%%%%%%%%%%%%%%%%%%%%%%%%%%%%%%%%%%%%%%%%%%%%%%%%%%%%%%%%%%%

\subsection{Covariance function with uncertain hyper-parameters}\label{sec:hyper}

From now on, we assume the prior of $m(\vec x)$ to be Gaussian and so completely characterized by its mean $\mu(\vec{x})$ and covariance function $\cal C$.
However, in many applications, not all the aspects of the covariance function are well-known \textit{a priori}. 
The stationarity of the covariance function can be easily determined and confirmed, yet, we have a large uncertainty in the other 
characteristics such as the values of the hyper-parameters. An example of a parametrized covariance function is
\begin{equation}
{\cal{C}} (\vec{x},\vec{x}')= \sigma^2_f \exp \left(-\frac{1}{2}(\vec{x}-\vec{x}')^T \bold{M}  (\vec{x}-\vec{x}') \right) 
+ \sigma^2_d\vec{x}^T\vec{x}' +\sigma^2_b + \sigma^2_n \delta_{pq}
	\label{eq:c_example}
\end{equation}
where $\bold{M}$ is a symmetric positive definite matrix.
The covariance hyper-parameters $\bold{M}$, $\sigma^2_f$, $\sigma^2_b$, $\sigma^2_d$, $\sigma^2_n$ are usually not exactly known \textit{a priori} and should be treated as uncertain quantities.
For many covariance functions it is easy to interpret the meaning of the hyper-parameters, which is of great importance when trying to understand the data. 
Traditionally, the hyper-parameters are estimated using Gaussian Process Regression (GPR) before inferring the model parameters~\citep{Rasmussen:2005}.
To this end, a set of possibly noisy observations of the field $m(\vec x)$ are used to perform stochastic interpolation of static data collected at few locations and maximize the marginal likelihood function using Bayesian inference or optimization techniques.
Optimal values of the inferred hyper-parameters are then used in the covariance function and KL expansion is applied as described in Eq.~\eqref{eq:kldec}.
The uncertainty bound can be estimated using GPR but is usually not considered in the expansion due to the complexity of the resulting model. 
This paper addresses the uncertainty in the hyper-parameters of covariance models. Specifically we develop a formulation that enables inferring the covariance function 
hyper-parameters along with the KL stochastic coordinates $\eta_k$. The formulation is based on basis transformations as described below.

\subsection{Stochastic coordinate transformation}
\label{sec:coordchg}
Without loss of generality, we assume that the stochastic prior process $M$ is centered ($\mu(\vec x)=0$) and has a parametrized covariance function 
${\cal C}(\vec x, \vec{x'}, \vec q)$ defined  by a random vector $\vec q\subset \mathbb R^h$ of hyper-parameters ($h$ is the number of hyper-parameters, \textit{e.g.}\ $\vec q=\{\bold{M},\sigma^2_f,\sigma^2_b,\sigma^2_d,\sigma^2_n\}$ for the example in Eq.~\eqref{eq:c_example}), with joint density $p_q$. 
Because of the dependence of the covariance function on $\vec q$, the KL expansion of $M$ in Eq.~\eqref{eq:kldec} becomes: 
\begin{equation}
	M_K(\vec x, \omega, \vec q) =  \sum_{k=1}^{K} \sqrt{\lambda_k (\vec q)} \phi_k(\vec x, \vec q) \eta_k (\omega), \quad 
	\int_D {\cal C}(\vec x, \vec{x'}, \vec q) \phi_k(\vec{x'},\vec q) dx = \lambda_k(\vec q) \phi_k(\vec x,\vec q).
\end{equation}
To simplify the notation, we drop the $\vec x$ and $\vec x'$  dependence and introduce the scaled eigen-functions $\Phi_k(\vec q)$:
\begin{equation}
  \Phi_k(\vec q) \doteq \sqrt{\lambda_k(\vec q)} \phi_k (\vec q), \quad \mbox{so} \quad M_K(\omega,\vec q) = \sum_{k=1}^K \Phi_k(\vec q) \eta_k(\omega).
	\label{eq:par:mode}
 \end{equation}

We further assume the continuity of the scaled eigen-functions $\Phi_k$ with respect to $\vec q$, in the sense (see~~\citep{Lax:1996})
$\exists D_k(\vec q)>0$, $ \|\Phi_k(\vec q) - \Phi_k(\vec q+\delta \vec q) \|_X^2 \le 
D_k(\vec q) \|\delta \vec q\|^2_{\ell^2_h}$, and  $\sum^{K}_{k=1} D_k(\vec q) \doteq D^{(K)}(\vec q) < \infty$ uniformly,
so that 
\begin{equation}
 \mathbb E \left[ \| M_K(\vec q) - M_K(\vec q+\delta \vec q) \|_X^2 \right]\le D_K(\vec q)  \|\delta \vec q\|^2_{\ell^2_h}.
\end{equation}
In practice, when decomposing a covariance function ${\cal C}(\vec q)$, the normalized eigen-functions are defined up to a factor of $\pm 1$.
To ensure the $\vec q$-continuity of the $\phi_k(\vec q)$'s, we have to select a consistent orientation of eigen-functions. A possibility, followed in this work, is to 
define the orientation of the eigen-functions with respect to a reference set of eigen-functions $\{\phi^r_k, k=1,\dots,K\}$, \textit{e.g.}\ the reference 
set defined below, such that $(\phi_k(\vec q),\phi^r_k)_X$ has a constant sign for all $\vec q$~\citep{Salloum:2012}. 
The dependence of the eigen-functions on the hyper-parameter $\vec q$ is further illustrated in Section~\ref{sec:exKL} below. 

Let ${\cal C}^r$ be a covariance function representative of the $\vec q$-dependent covariance function ${\cal C}(\vec q)$ of $M$. 
As investigated below, a possible choice for ${\cal C}^r$ can be
\begin{equation}
 {\cal C}^r = \overline{\cal C}\doteq \int {\cal C}(\vec q) p_q(\vec q) d\vec q,
 \end{equation}
that is the $\vec q$-averaged of ${\cal C}(\vec q)$,
or a particular realization of $\cal C$ corresponding to a deterministic value
${\vec{q}^r}$ of the random parameters (\textit{e.g.}\ nominal values obtained using GPR~\citep{Rasmussen:2005}). We denote $\phi^r_k$ the ordered and normalized eigen-vectors of ${\cal C}^r$. Note that $\{\phi^r_k, k=1,2,\dots\ ,\infty \}$ is an orthonormal basis of $X$; 
as a result, any scaled eigen-function $\Phi_k(\vec q)$ can be expressed in this basis:
\begin{equation}
 \Phi_k (\vec q) = \sum_{k'=1}^{\infty} b_{kk'} (\vec q) \phi^r_{k'}, \quad b_{kk'}(\vec q) 
 = \left( \phi^r_{k}, \Phi_{k'}(\vec q)\right)_X. \label{eq:chbasis}
  \end{equation}
The continuity of the scaled eigen-functions implies the continuity of the projection coefficients $b_{kk'}(\vec q)$.
For computational purposes, the expansion in Eq.~\eqref{eq:chbasis} needs to be truncated to the first $K^r$ terms. 
Without loss of generality we shall use in the following $K^r=K$, allowing for convergence analysis with respect to a single parameter $K$.

Further, the change of basis gives:
\begin{eqnarray}
	M_K(\omega,\vec q) = \sum_{k=1}^{K} \Phi_k (\vec q) \eta_k(\omega) \approx \sum_{k=1}^{K}\left(\sum_{k'=1}^{K} b_{kk'}(\vec q) \phi^r_{k'} \right) \eta_k(\omega)
	= \sum_{k=1}^{K} \phi^r_k \hat{\eta}_k(\omega, \vec q),
\label{eq:cbasis}
\end{eqnarray}
where we have denoted
\begin{equation}
\hat{\eta}_k(\omega, \vec q) = \sum_{k'=1}^{K} b_{k'k}(\vec q)\eta_{k'}(\omega). \label{chgcoord}
\end{equation}
The transformation shows that the $\vec q$-dependence of $\cal C$ can be translated into an expansion $M_K$ with $\vec q$ dependent scaled 
eigen-functions, see Eq.~\eqref{eq:par:mode}, or approximated by a $\vec q$-dependent linear transformation of the random variables in Eq.~\eqref{eq:cbasis}. Specifically, denoting the latter approximation $\hat{M}_K$ we have the approximations
\begin{equation}
	M(\omega,\vec q) \approx M_K(\omega,\vec q) = \sum_{k=1}^K \Phi_k(\vec q) \eta_k(\omega) \approx \hat M_K(\omega,\vec q) = \sum_{k=1}^K
	\phi^r_k \hat\eta_k(\omega,\vec q),
\end{equation}
with $\hat \eta_k$ related to the $\eta_k$'s by Eq.~\eqref{chgcoord}. 

We observe that $\hat{\eta}_k(\omega, \vec q)$ is a linear combination of standard Gaussian random variables, 
so it is also Gaussian (with zero mean). However, the $\hat{\eta}_k(\omega, \vec q)$ are generally correlated. 
The change of random coordinates in Eq.\eqref{chgcoord} can be cast in matrix form:
\begin{equation}
\hat{\vec {\eta}} (\omega,\vec q) = {\cal{B}} (\vec q) \vec \eta(\omega).
\label{eq:trans}
\end{equation}
The covariance matrix for the random coefficients $\hat{\vec \eta}$, denoted $\Sigma^2(\vec q)$, can be expressed as
\begin{equation}
 \Sigma^2 (\vec q) = \mathbb E \left[ \hat{\vec\eta}(\vec q) \hat{\vec\eta}^t(\vec q) \right] = {\cal B }(\vec q){\cal B}^t(\vec q).
\end{equation}
We shall assume that $\Sigma^2(\vec{q})$ is invertible (for almost every $\vec q$); a sufficient condition is that  $\Phi_{1\le k\le K}(\vec q)$ is not orthogonal to $\mbox{span} \{ \phi^r_1, \dots, \phi^r_{K} \}$. In addition, the conditional distribution of $\hat{\vec \eta}$, given $\vec q$, is
\begin{equation} 
	p_{\hat{\eta}} (\hat{\vec \eta}| \vec{q}) = \frac{1}{\sqrt{2\pi^K |\Sigma^2 (\vec q)|}} 
\exp\left[- \frac{\hat{\vec\eta}^t (\Sigma^{2})^{-1}(\vec{q}) \hat{\vec\eta}}{2} \right], \label{cond:dist:eta}
\end{equation} 
where $|\Sigma^2 (\vec q)|$ is the determinant of $\Sigma^2 (\vec q)$.

\subsection{Example}
\label{sec:exKL}
We now provide a brief illustration of the convergence of the error in the approximation of $M(\vec q)$. 
To this end, we consider $D=[0,1]$ and a centered Gaussian process $M$ with covariance function 
\begin{equation}
{\cal{C}} (x,x',\vec q)= \sigma^2_f \exp \left(-\frac{(x-x')^2}{2l^2} \right), 
\label{eq:kernelSE}
\end{equation}
with hyper-parameter vector $\vec q = \{\sigma_f^2,l \}$. In this case, only the correlation length $l$ affects the shape of the eigen-functions, 
while the process variance $\sigma_f^2$ simply scales the eigen-values. Therefore, we fix $\sigma^2_f=0.5$  through-out the section and assume uncertainty in $l$ only, that is $\vec q = \{ l\}$. Specifically, we assume the hyper-parameter $l$ to have a uniform distribution in the range $[0.1,1]$.
It is important to note that the number of KL modes needed for convergence highly depends on the hyper-parameter $l$. 
In particular, if $M$ has small-scale features (small $l$) a large number of KL modes will be needed.

\begin{figure}[htb]
\centering
\begin{tabular}{clcl}
\includegraphics[width=0.45\textwidth]{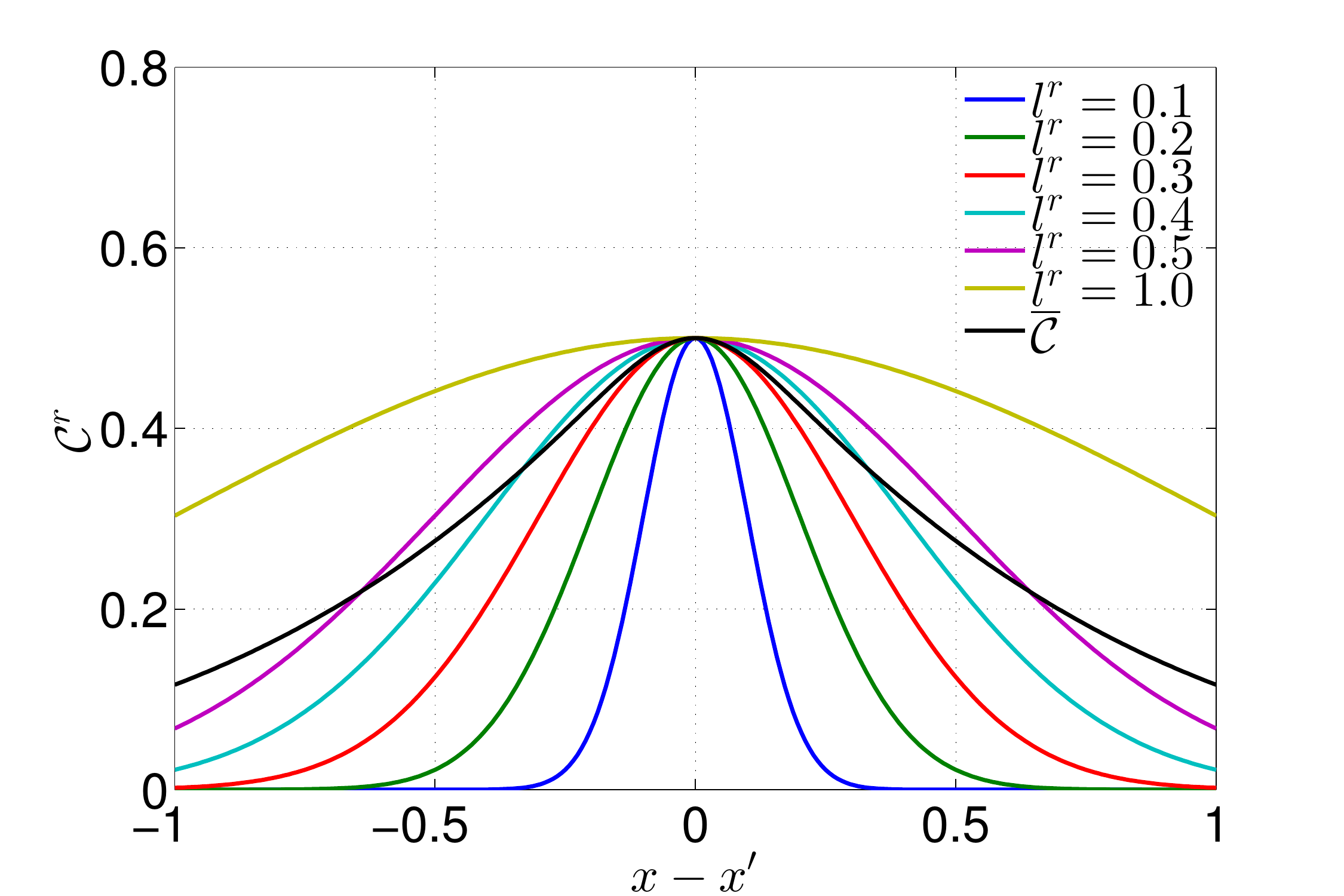} &
\includegraphics[width=0.45\textwidth]{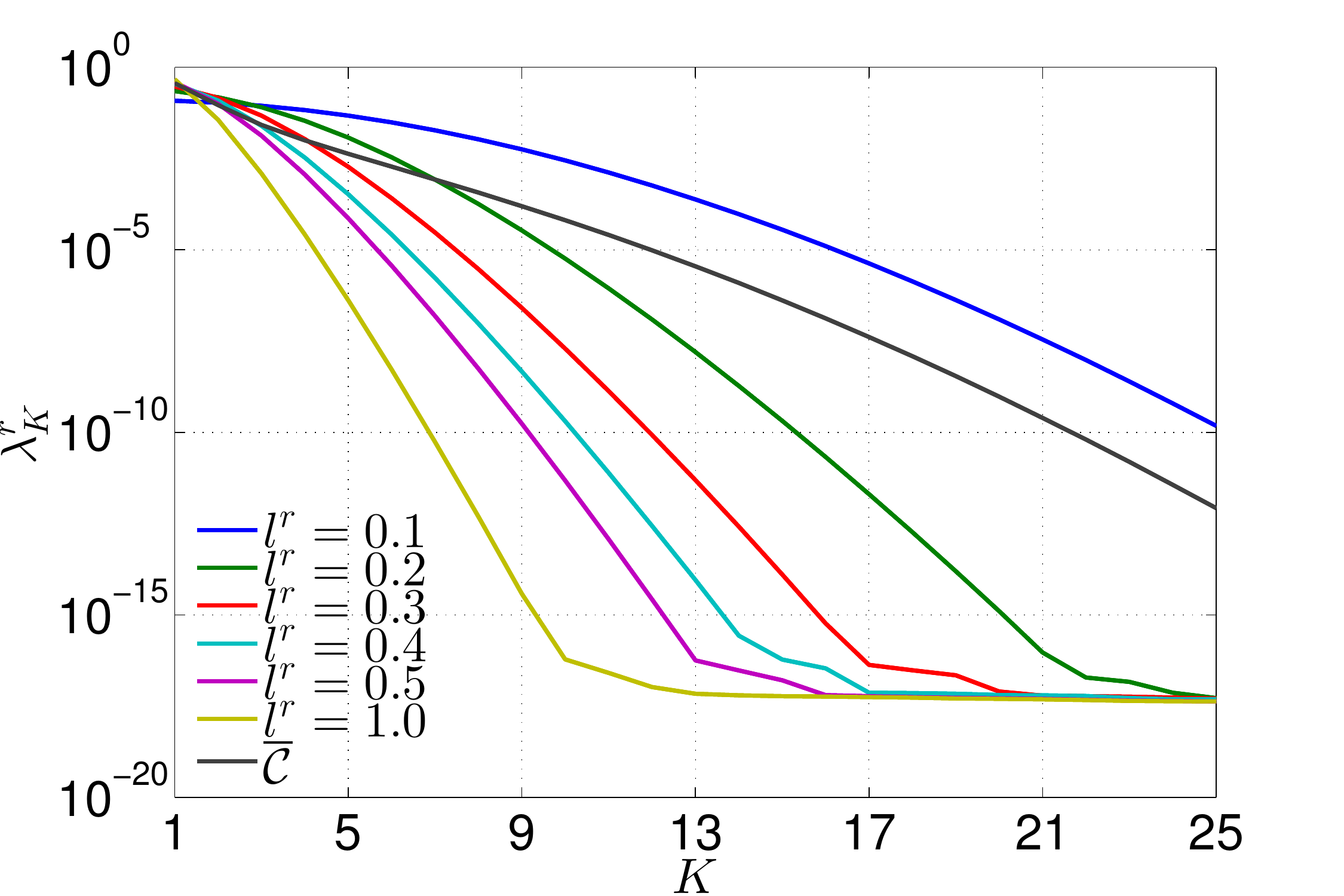} 
\end{tabular}
\caption{(Left) Reference covariance functions ${\cal C}^r={\cal C}(l^r)$ for different values of $l^r$, as indicated. Also plotted is the $\vec q$-averaged covariance $\overline{\cal C}$ . (Right) Spectra of the corresponding eigen-values decay with $K$.
}
\label{fig:errorinMa}
\end{figure}

For the selection of the reference covariance function, we contrast the choice ${\cal C}^r = {\cal C}(l^r)$, for several values 
$l^r \in [0.1,1]$, with the case ${\cal C}^r = \overline{\cal C}$.
The KL decompositions are numerically approximated with Galerkin piecewise constant modes over a uniform grid having $N=128$ elements in space. 
Figure~\ref{fig:errorinMa} compares in the left plot the considered reference covariance functions ${\cal C}^r$  
and in the right plot the respective decay rates with $k$ of their eigen-values ${\lambda}^r_k$. 
When using ${\cal C}(l^r)$, it is seen that the smaller $l^r$ the slowest the decay rate, as expected, whereas for the $\vec q$-averaged covariance $\overline{\cal C}$ the decay rate is asymptotically similar (but with a lower magnitude) to the lowest 
$l^r$ in the uncertainty range. 
Also, note that $\overline{\cal C}$ is evidently not Gaussian.

To quantify the error in the approximation of $M(\omega,\vec q)$ by the proposed transformation method, we introduce the following relative error measure
\begin{equation}
	\epsilon_M(K,{\vec q}) = \frac{\norm{M(\vec q)-\hat M_{K} (\vec q)}}{\norm{M(\vec q)}}, 
\end{equation}
where 
\begin{equation}
	\norm{U}^2 \doteq \mathbb E\left[ (U,U)_X \right].
\end{equation}
The error $\epsilon_M(K,\vec q)$ integrates both the truncation error in approximating $M$ with $M_K$, and the subsequent projection error of $M_K$ into the space of reference modes ${\phi}^r_k$. The error $\epsilon_M$ is estimated by means of Monte Carlo sampling where realizations of $M$ are generated given $\vec q$; these realizations are projected on the $K$-dimensional dominant space of ${\cal C}(\vec q)$ in order to compute the coordinates $\eta_k$ (see Eq.~\eqref{def:projcoeff}) which are transformed using Eq.~\eqref{chgcoord} to obtain the corresponding realizations of $\hat M_K$. Observe also that $\norm{M(\vec q)} = \sigma_f$. Finally, the local (squared) error $\epsilon_M^2(K,\vec q)$ can be averaged over $\vec q$ to yield the averaged error, which we denote $E_M(K)$.

\begin{figure}[htb]
\centering
\begin{tabular}{cc}
\includegraphics[width=0.45\textwidth]{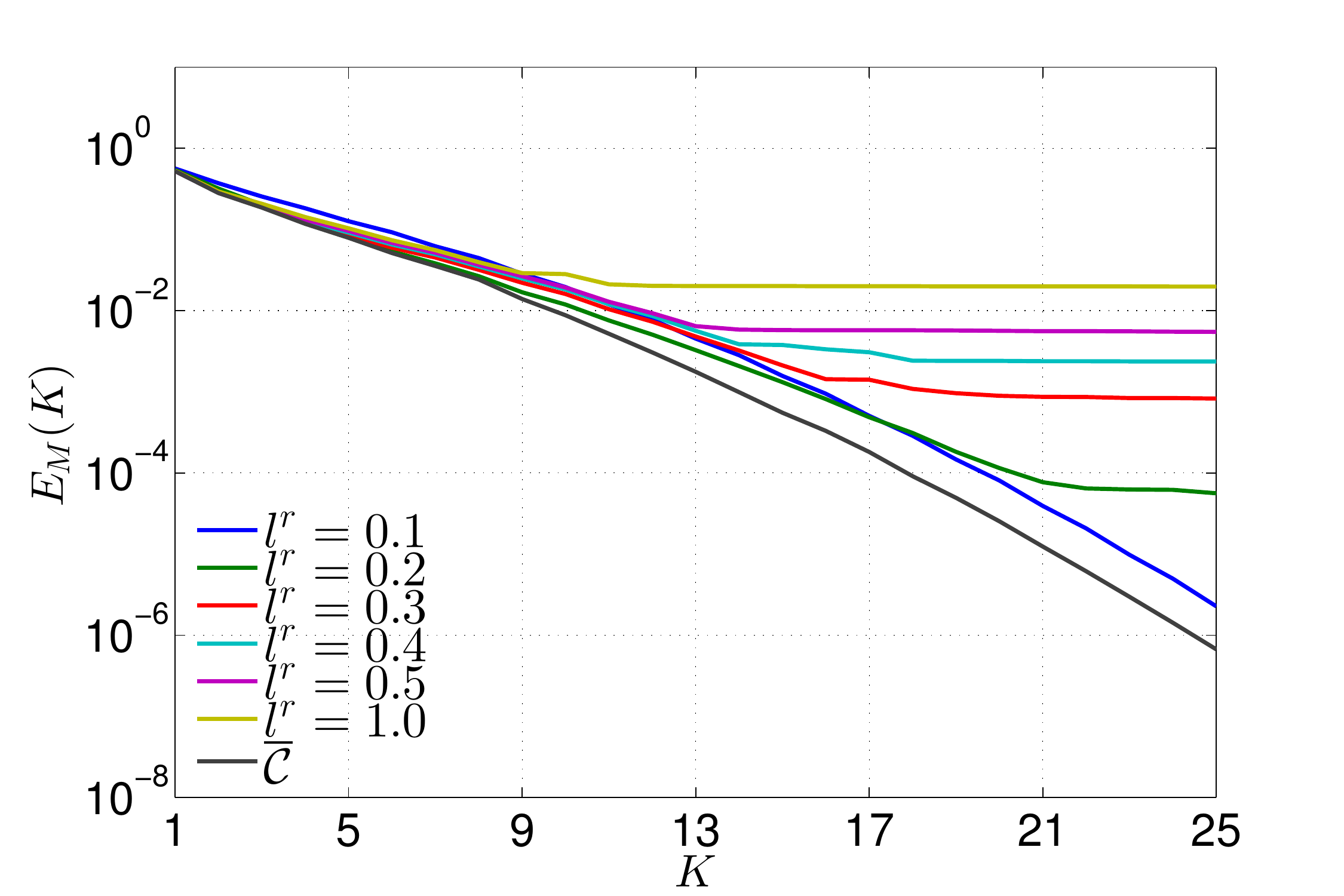} &
\includegraphics[width=0.45\textwidth]{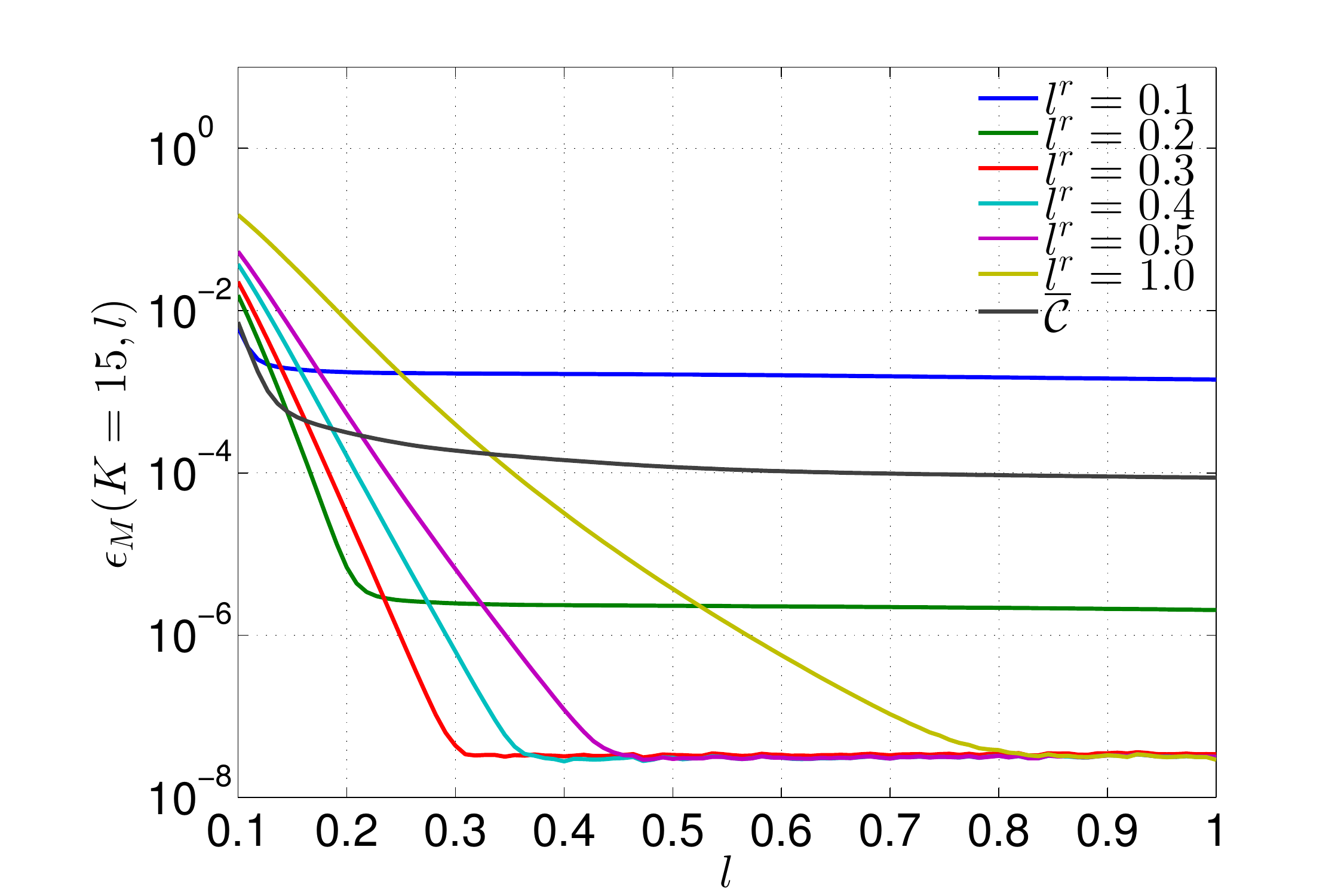}
\end{tabular}
\caption{(Left) Error $E_M(K)$ in approximating the Gaussian Process $M$ by $\hat M_K$ for different 
reference covariance functions based on selected correlation lengths $l^r$ as indicated. Also plotted are results obtained with $\overline{\cal C}$. (Right) Relative error $\epsilon_M(K=15,l)$ for the same cases as in the left plot.}
\label{fig:errorinM}
\end{figure}

\begin{figure}[htb]
\centering
\begin{tabular}{clclcl}
\includegraphics[width=0.30\textwidth]{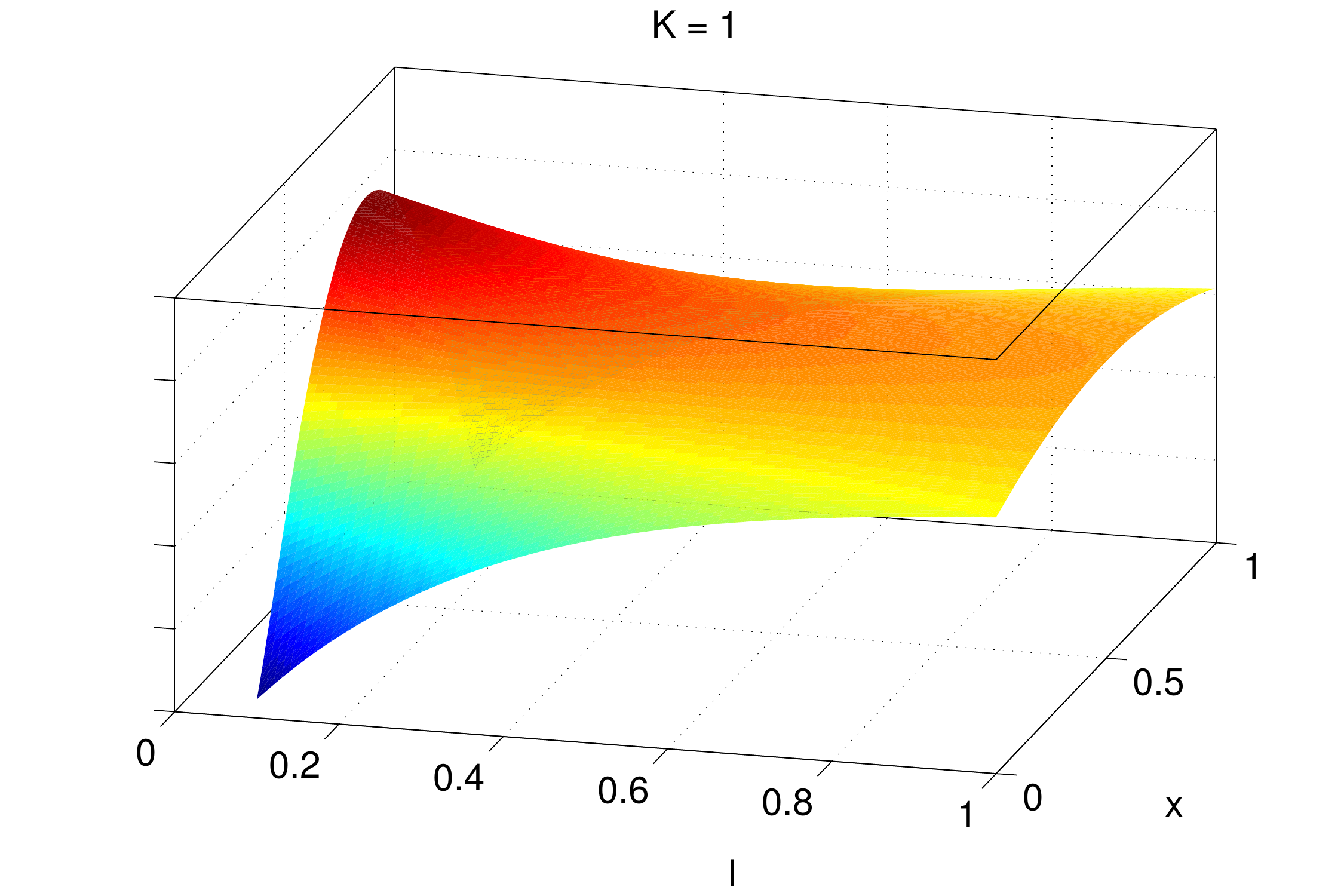} &
\includegraphics[width=0.30\textwidth]{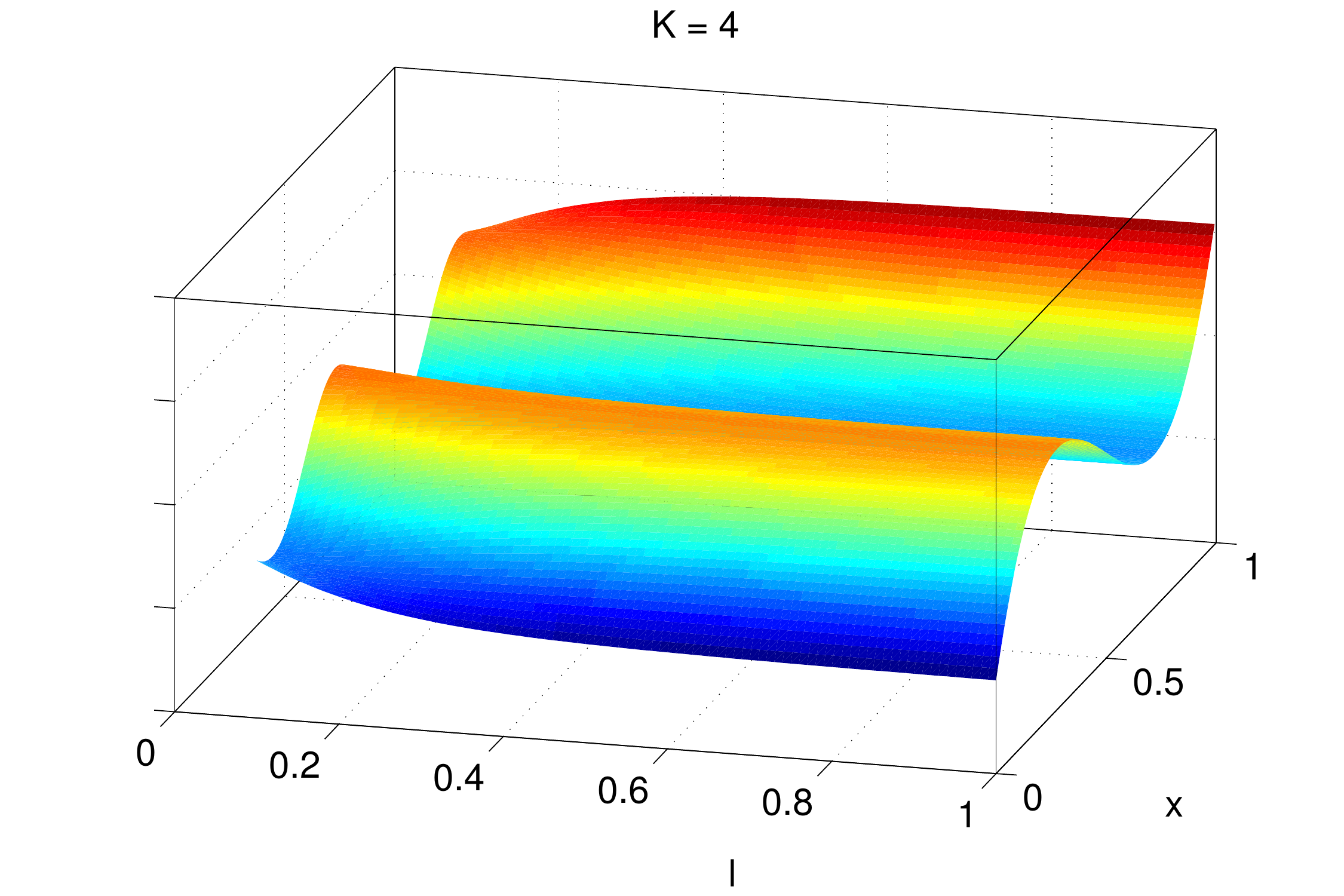} &
\includegraphics[width=0.30\textwidth]{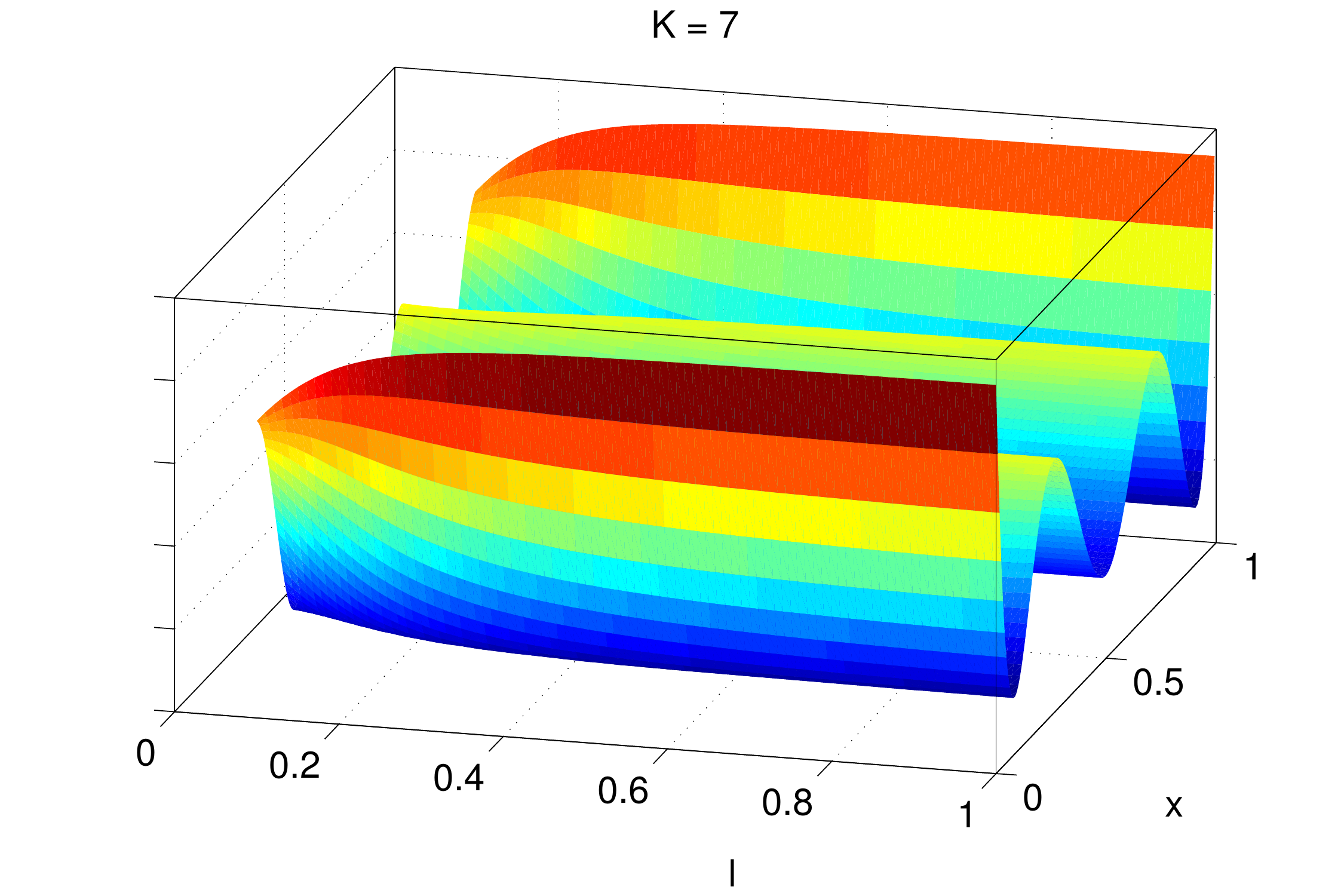} \cr
\includegraphics[width=0.30\textwidth]{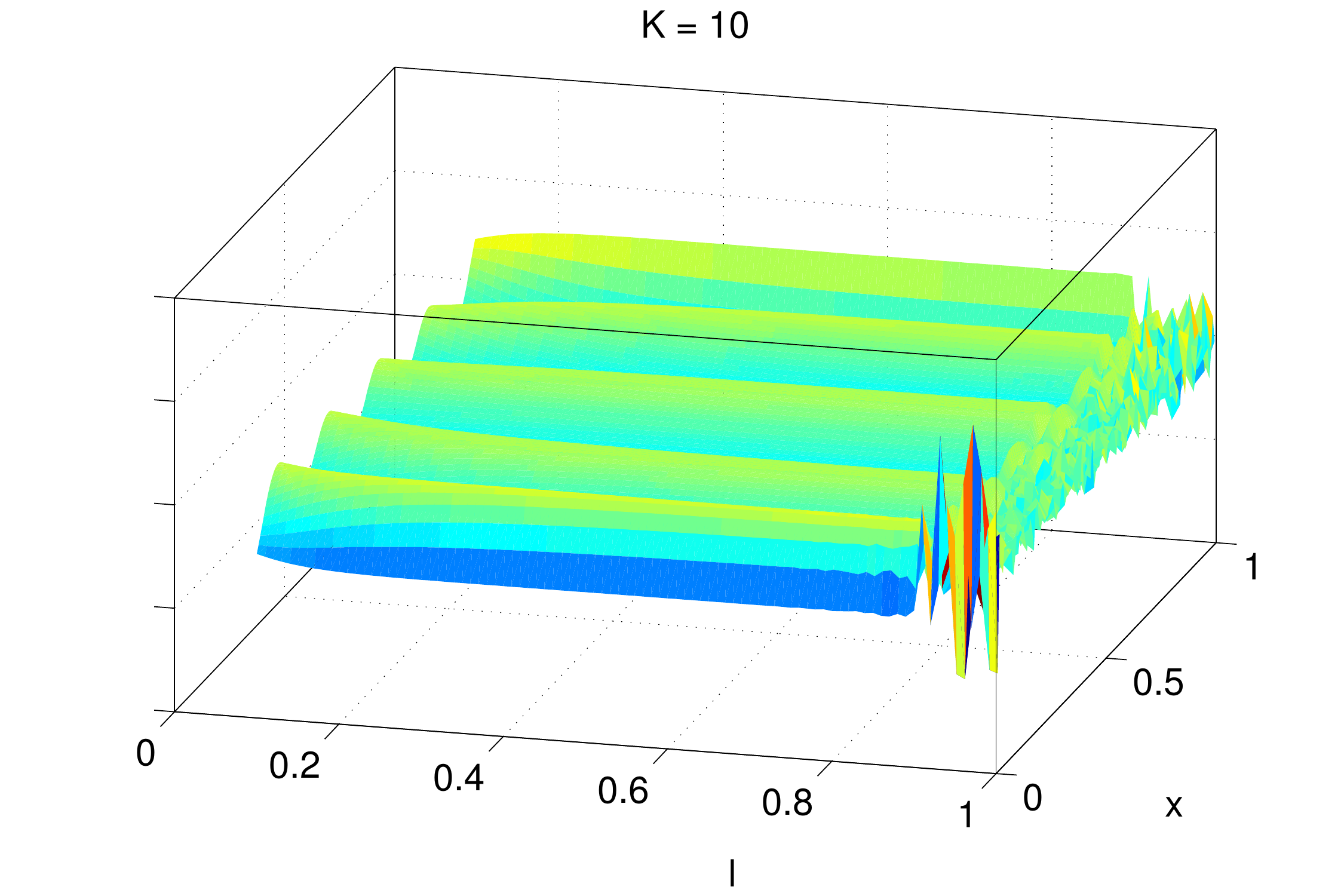} &
\includegraphics[width=0.30\textwidth]{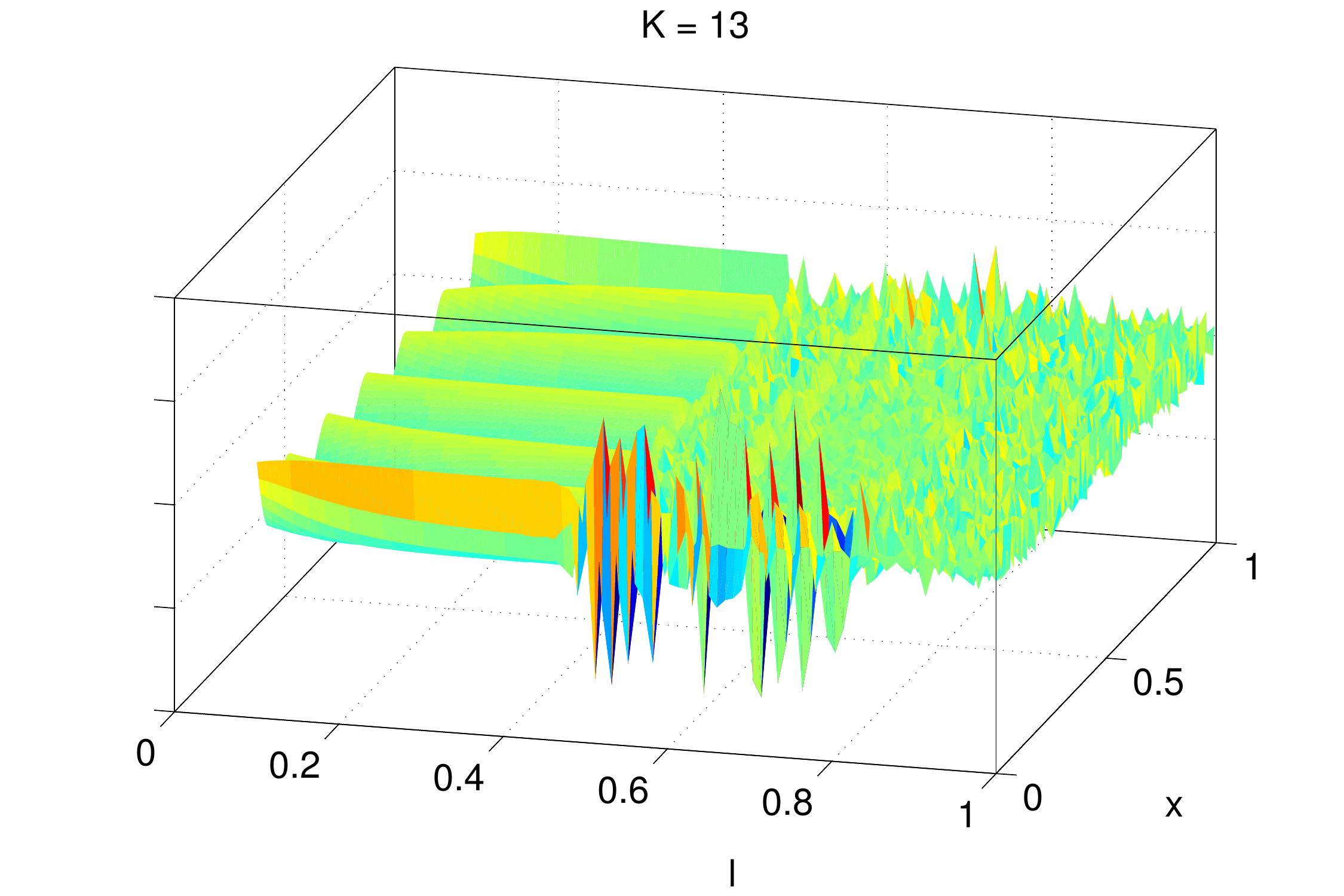} &
\includegraphics[width=0.30\textwidth]{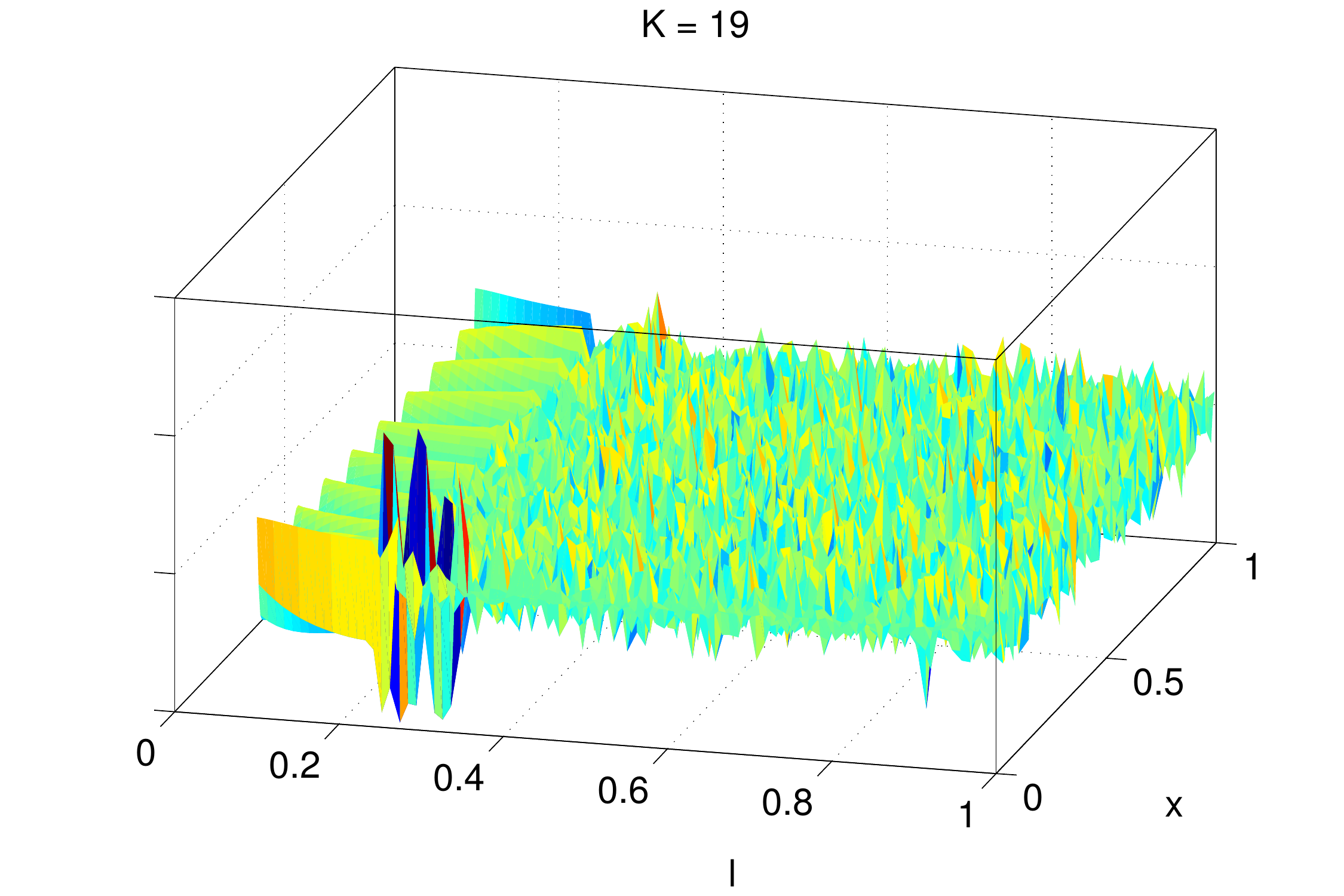}  
\end{tabular}
\caption{Dependence of eigen-functions $\phi_k(\vec q)$ with the length-scale hyper-parameter $l$ and selected $k$ as indicated.
}
\label{fig:eigdepl}
\end{figure}

The mean square error $E_M(K)$ is shown in the left plot of Figure~\ref{fig:errorinM}.
Plotted are curves for different reference bases: using ${\cal C}^r={\cal C}(l^r)$ with selected correlation lengths $l^r$ within $[0.1,1.0]$, and the $\vec q$-averaged covariance function $\overline{\cal C}$.
A first comment from these curves is that the error decreases as $K$ increases as expected.
However, for $l^r>0.1$, the error $E_M(K)$ stagnates as $K$ increases when using ${\cal C}^r = {\cal C}(l^r)$.
The stagnation occurs at lower $K$ when $l^r$ increases. This stagnation can be explained from the spectra reported in Figure~\ref{fig:errorinMa} which shows
that when using $l^r>0.1$ the magnitude of $\lambda^r_k$ quickly decays with $k$ to zero machine precision, such that subsequent modes are not correctly estimated and cannot provide a suitable projection basis. 
To further illustrate the effect of finite numerical accuracy, we provide in Figure~\ref{fig:eigdepl} plots of eigen-functions $\phi_k(x,l)$ for selected $k$ and $(x,l)\in D\times[0.1,1]$. 
It is seen that for $k=1,$ 4 and 7, the dependence on $l$ of the numerical eigen-functions is smooth. In contrast, for $k=10$ (resp. 13 and 19) the computed eigen-functions are seen to be noisy for $l\gtrsim 0.9$ (resp. $l\gtrsim 0.5$ and 0.25) because of 
finite numerical accuracy. Clearly, this indicates that under-resolved modes could be disregarded and that the reference basis should include only modes with indices $k$ 
such that $\lambda^r_k / \lambda^r_1$ remains in achievable accuracy ($\approx10^{-16}$ for double precision). To keep the analysis simple, and because our
approach is in fact robust to under-resolved modes, we continue in the following to compare for the same $K$ the different choices of reference covariance functions.  
Note also that for the reference basis using the shortest correlation length, $l^r=0.1$, and the $\vec q$-averaged covariance, this numerical issue has not yet emerged for the range of considered $K$, and the corresponding errors decay monotonically up to $K=25$. 
In addition, it is seen that the error curve corresponding to the  $\vec q$-averaged reference covariance function $\overline{\cal C}$ has the lowest approximation error $E_M(K)$ for all $K$.
 
This is not a surprise since by construction this choice uses eigen-functions $\phi^r_k$ spanning the optimal subspace to represent $M(\omega,\vec q)$ when $\vec q$ varies with law $p_q(\vec q)$. In fact, finding the $K$-term expansion minimizing the $\vec q$-averaged mean square error approximation of $M\sim {\cal GP}(0,{\cal C}(\vec q))$ amounts precisely to the decomposition of $\overline{\cal C}$. In other words, if using the $K$ dominant eigen-modes of ${\cal C}^r = \overline{\cal C}$ to construct the reference basis is non optimal to represent $M(\vec q)$ for any value $\vec q$ (obviously for each $\vec q$ the optimal choice is the eigen-modes of ${\cal C}(\vec q)$), there is no better choice on average over $\vec q$.

To better appreciate the behavior of the error with the hyper-parameter $l$ in the present example, 
the right plot in Figure~\ref{fig:errorinM} reports the evolution of $\epsilon_M(K,{l})$ for $K=15$, using the same reference covariance functions considered previously.
It is seen that for all reference covariance functions, the error $\epsilon_M(K=15,l)$ increases when $l$ decreases, reflecting the increasing truncation error for $K=15$ when $M$ involves smaller features. 
However, different behaviors are reported depending on the choice of ${\cal C}^r$ when $l$ increases. 
When using  ${\cal C}^r = {\cal C}(l^r)$ with $l^r\ge 0.3$, the error converges to machine precision when $l \gtrsim l^r$, meaning that in this situation the $15$-dimensional reference subspace 
${\rm span} \left\{  \phi^r_k = \phi_k (l^r), k=1,\dots,K\right\}$ essentially encompasses the $15$-dimensional dominant subspace of ${\cal C}(l \gtrsim l^r)$. 
Further, this behavior highlights the robustness of the change of coordinates, even for situations where finite numerical accuracy prevents the correct determination of the whole set of eigen-functions. 
On the contrary, the choice ${\cal C}(l^r)$ with $l^r \le 0.2$, while yielding a lower error at small correlation length $l\lesssim l^r$, exhibits a stagnating error for $l\gtrsim l^r$, 
denoting that the corresponding $K=15$-dimensional reference subspace is not rich enough to encompass the dominant subspaces at larger correlation lengths. 
Roughly speaking, the reference eigen-functions are too oscillating to properly represent processes with long-range correlations.
Finally, the selection of $\overline{\cal C}$ for the reference covariance function provides the best compromise, by construction, maintaining  a maximum error $\epsilon_M(K=15,l)$ less than $10^{-2}$ over the whole range of $l$.

%%%%%%%%%%%%%%%%%%%%%%%%%%%%%%%%%%%%%%%%%%%%%%%%%%%%%%%%%%%%%%%%%%%%%%%%%%%%%%%%%%
%%%%%%%%%%%%%%%%%%%%%%%%%%%%%%%%%%%%%%%%%%%%%%%%%%%%%%%%%%%%%%%%%%%%%%%%%%%%%%%%%
%%%%%%%%%%%%%%%%%%%%%%%%%%%%%%%%%%%%%%%%%%%%%%%%%%%%%%%%%%%%%%%%%
\section{Polynomial Chaos Surrogate}
\label{sec:pc}

A suitable Polynomial Chaos (PC) expansion for the model predictions is constructed to accelerate the Bayesian inference process. 
In Section~\ref{pc:framework} we briefly review the PC methodology and provide some details regarding the numerical methods used in the examples provided in Section~\ref{sec:results}.
Then, in Section~\ref{pc:chftbasis} we focus on exploiting the PC surrogates to efficiently handle uncertain hyper-parameter through the 
change of coordinates introduced previously in Section~\ref{sec:coordchg}. Finally, Section~\ref{pc:example} provides a brief analysis
of the PC surrogate error.
 
\subsection{Polynomial Chaos expansion}
\label{pc:framework}
Polynomial Chaos (PC) is a probabilistic methodology that expresses the dependencies of a model solution on some uncertain model inputs, through a truncated 
spectral polynomial expansion~\citep{GhanemSpanos1991,LeMaitreKnio2010}.
Let $U\in Y$ be solution of a mathematical model $\cal L$ (\textit{e.g.}\ Partial Differential Equations), formally expressed as ${\cal L}U = 0$.
We are interested in situations where the model $\cal L$ is uncertain and parametrized with a finite set of independent second-order random 
variables $\xxi=(\xi_1,...,\xi_N)$ with known probability distribution. 
For simplicity, we shall restrict ourselves to the case of i.i.d.\ standard Gaussian random variables $\xi_i$, 
and will denote $p_\xxi$ the density function of $\xxi$, and $L_2(p_\xxi)$ the space of second order random functionals in $\xxi$, that is
\begin{equation}
		v(\xxi) \in L_2(p_\xxi) \Leftrightarrow \dotsint |v(\xxi)|^2 p_\xxi (\xxi) d\xxi < \infty.
\end{equation}
Since the model depends on $\xxi$, its solution is also generally dependent on $\xxi$ and satisfies
\begin{equation}
	{\cal L}(\xxi) U(\xxi) = 0, \quad {a.s.} \label{strong:form}
\end{equation}
Let $\{\Psi_\alpha, \alpha\in \mathbb N\}$ be a complete orthonormal set of $L_2(p_\xxi)$, such that the model solution has an expansion of the form
\begin{equation}
	U(\xxi) = \sum_{\alpha\in \mathbb N} U_\alpha \Psi_\alpha(\xxi), \quad \left<\Psi_\alpha, \Psi_\beta\right> \doteq \dotsint \Psi_\alpha(\xxi) \Psi_\beta(\xxi) p_\xxi(\xxi)d\xxi = \delta_{\alpha,\beta},
\end{equation}
where the equality stands in the mean square sense and the expansion coefficients $U_\alpha \in Y$ are called the stochastic modes of $U$. 
A classical choice for the random functionals $\Psi_\alpha$ are orthonormal multi-variate polynomials in $\xxi$, leading to the so-called PC expansion of $U(\xxi)$. 
The $\xi_i$ being standard Gaussian random variables, the $\Psi_\alpha$ are in fact normalized multi-variate Hermite polynomials~\citep{Wiener:1938}.
For practical purposes, the PC expansion of $U(\xxi)$ needs to be truncated. When the basis is truncated to total order $o$ the total number of terms in the 
PC expansion is given by $P+1=(N+o)!/(N!o!)$ and therefore increases exponentially fast with both the expansion order $o$ and the number $N$ of random variables $\xi_i$.
The series expansion approximating $U(\xxi)$ is then finite and will be denoted $\tilde U(\xxi)$ in the following:
\begin{equation}
	U(\xxi) \approx \tilde U(\xxi) \doteq \sum_{\alpha=0} ^P U_\alpha \Psi_\alpha(\xxi). \label{eq:stochseries}
\end{equation} 
The existence and convergence of this series is asserted by the Cameron-Martin theorem \citep{Cameron:1947} with the condition of $U$ having a finite variance. 
The rate of convergence, and hence the number of terms in the series, depends on the smoothness of $U$ with respect to $\xxi$. 
The series converges spectrally fast with $P$ when $U$ is infinitely smooth.
\par 
Various methods have been proposed for the determination of the PC coefficients $U_\alpha$. They can be distinguished into the Non-intrusive and Galerkin methods.
Non-intrusive methods rely on an ensemble of deterministic model evaluations of $U(\xxi)$, for particular realizations of $\xxi$ selected either at random or deterministically. 
Non-Intrusive methods include Non-Intrusive Spectral and Pseudo-Spectral Projection~\citep{Reagan:2003,Constantine:2012,Conrad:2013}, 
Least-Square-Fit and regularized variants~\citep{Berveiller:2006,Blatman:2011,Doostan:2014}, 
Collocation (interpolation) methods~\citep{Babuska:2007,Xiu_Hesthaven,Nobile:2008a}, 
that are often combined with Sparse-Grid algorithms to reduce computational complexity. 

In the present paper, we instead rely on the Galerkin projection method~\citep{GhanemSpanos1991,LeMaitreKnio2010} for which the expansion coefficients 
$U_\alpha$ are defined through a reformulation of the model Eq.~\eqref{strong:form}, 
using a weak form at the stochastic level. Specifically, Eq.~\eqref{strong:form} is projected on the PC basis, 
a procedure resulting in a set of $P+1$ coupled problems,
\begin{equation}
	\left< {\cal L}(\xxi) \sum_{\alpha=0}^P U_\alpha \Psi_\alpha(\xxi),  \Psi_\beta(\xxi)\right> = 0, \quad \beta=0,\dots,P.
\end{equation}
Numerical algorithms have been proposed to efficiently solve this set of coupled problems, both in the case of linear operators $\cal L(\xxi)$ 
(see \textit{e.g.}\ \citep{LeMaitre:2002b} for elliptic and parabolic problems) and non-linear operators (see \textit{e.g.}\ \citep{LeMaitre:2001,LeMaitre:2004z} 
and references in~\citep{LeMaitreKnio2010}).

\subsection{PC surrogate for a parametrized covariance}
\label{pc:chftbasis}
Returning to the inference problem, we now want to construct a global PC surrogate for the model predictions, that accounts both for randomness of $M_K$, 
through its random coordinates $\vec \eta$, and the uncertainty in its covariance function, through the random hyper-parameter vector $\vec q$. 
We assume that the model problem amounts to solving for $U$ a model depending on $M_K$. Using the notations above, it is written formally as 
\begin{equation}
	{\cal L}(\vec \eta, \vec q) U(\vec \eta,\vec q) =  0.
\end{equation}
The previous equation has motivated the idea of expanding the dependence of $U$ with respect to the random vectors $\vec \eta$ and $\vec q$ on a PC 
basis~\citep{MarzoukNajm2009,Tagade:2014},
that is using $U(\vec\eta,\vec q)\approx \sum_{\alpha} U_\alpha\Psi_\alpha(\vec \eta,\vec q)$. In the following, we consider an alternative approach, taking advantage of 
the change of coordinates discussed in Section~\ref{sec:KL}. 
The change of coordinates allows us to approximate $M_K(\vec q)$ on the fixed reference basis of KL modes $\{\phi^r_k, k=0,\ldots,K\}$, 
through the linear mapping $\vec \eta \mapsto \hat{\vec \eta}(\vec q) = {\cal B}(\vec q)\vec \eta$. Eq.~\eqref{cond:dist:eta} provides the density of
$\hat{\vec \eta}$ conditioned on $\vec q$. The model problem can therefore be recast as
\begin{equation}
	{\cal L}(\hat{\vec \eta}) U(\hat{\vec \eta}) = 0, \quad \mbox{where }\hat{\vec{\eta}} \sim p_{\hat{\eta}} (\hat{\vec \eta}| \vec{q}).
\end{equation}
The last expression shows that we only need to construct an approximation of the mapping $\hat{\vec \eta} \mapsto U(\hat{\vec \eta})$ which is accurate enough with 
respect to the conditional density $p_{\hat{\eta}} (\hat{\vec \eta}| \vec{q})$ when $\vec q$ varies. To get rid of the $\vec q$-dependence of the conditional density, we can
consider averaging $p_{\hat{\eta}}$ over $\vec q$. In the case of the reference covariance function ${\cal C}^r = \overline{\cal C}$, it can be shown that
\begin{equation}
	\dotsint   p_{\hat{\eta}} (\hat{\vec \eta}| \vec{q}) p_q(\vec q) d\vec q = \frac{1}{\sqrt{2\pi^K |\Lambda^2|}} 
	\exp\left[- \frac{\hat{\vec\eta}^t (\Lambda^{2})^{-1}\hat{\vec\eta}}{2} \right],
\end{equation}
where $\Lambda^2 = {\rm diag\:}(\lambda^r_1, \dots, \lambda^r_K)$. In other words, the $\vec q$-marginal of the conditional density yields independent Gaussian random variables.
This suggests constructing an approximate mapping of $\hat{\vec \eta} \mapsto U$, solving the model problem for a \emph{reference} Gaussian field defined as 
\begin{equation}
	\hat M_K^{\rm PC}(\xxi) = \sum_{k=1}^K \sqrt{\lambda^r_k} \phi^r_k \xi_k, \label{def:RGPPC}
\end{equation}
where the $\xi_k$'s are independent standard Gaussian random variables. 
It corresponds to a reference model problem $\hat{\cal L}(\xxi)$ based on the reference Gaussian process
$\hat M_K^{\rm PC}(\xxi)$. As before, we denote $\tilde U(\xxi)$ the PC approximation of the reference model problem $\hat{\cal L}(\xxi)U(\xxi)=0$.
From this PC approximation, we can approximate the model problem solution for couples $(\vec\eta,\vec q)$ through
\begin{equation}
	U(\vec\eta,\vec q) \approx \tilde{U}(\xxi(\vec \eta,\vec q)) = \sum_{\alpha=0}^P U_\alpha\Psi_\alpha(\xxi(\vec\eta,\vec q)),
	\quad \vec \xxi(\vec\eta,\vec q) = \hat{\cal B}(\vec q) \vec \eta,  \label{chgBasis:PC1}
\end{equation}
where the $\vec q$-dependent matrix $\hat{\cal B}(\vec q)$ expresses the change of coordinates $(\vec \eta,\vec q) \mapsto \xxi(\vec \eta,\vec q)$. Based on 
Eq.~\eqref{def:RGPPC}, we propose to use
\begin{equation}
 \hat{\cal B}_{kl}(\vec q) = \begin{cases}\displaystyle \frac{{\cal B}_{kl}(\vec q)}{\sqrt{\lambda^r_k}}, &
\lambda^r_k / \lambda^r_1 >\kappa, \cr
 0, &\mbox{otherwise,} \end{cases}
 \label{chgBasis:PC}
\end{equation}
where $\kappa>0$ is a small constant related to the numerical accuracy (typically $\kappa \sim 10^{-12}$) introduced to avoid ill-definition of the $\xi_k$'s 
associated to negligibly small $\lambda^r_k$. 
The $\kappa$-thresholding leads to transformed coordinates $\xxi(\vec \eta,\vec q)$ where the first $K^{\rm PC}(\kappa)$ components are  non trivial, 
with $K^{\rm PC}(\kappa)\le K$. 
Note that the number of non-trivial components of $\xxi$ only depends on the reference covariance, ${\cal C}^r$, and not on 
$\vec q$. Also, the PC construction of $\tilde U(\xxi)$ can in fact be reduced, considering $\hat M^{\rm PC}_{K^{\rm PC(\kappa)}}$ instead of
$\hat M^{\rm PC}_{K}$, with computational complexity reduction as a result when $K^{\rm PC} < K$. 
However, we shall continue to report results as a function of $K$ for simplicity.

When the reference covariance ${\cal C}^r$ is not the $\vec q$-averaged of ${\cal C}(\vec q)$, the $\vec q$-marginal conditional density $p_{\hat \eta}$ 
remains Gaussian but introduces correlations between components. These correlations can be dealt with by introducing an additional change of basis in order 
to redefine a reference Gaussian process $\hat M_K^{\rm PC}$ in terms of independent standard Gaussian random variables $\xi_k$'s. 
In that case, Eq.~\eqref{chgBasis:PC} must be accordingly modified to account for the additional change of coordinates. 
Alternatively, when using ${\cal C}^r = {\cal C}({\vec q}^r)$, we can continue to define the reference Gaussian process 
$\hat{M}_K^{\rm PC}$ by Eq.~\eqref{def:RGPPC}, which corresponds to solving the uncertain model problem assuming that $\hat{\vec \eta}$ has 
for density $p_{\hat \eta}$ conditioned on $\vec q = {\vec q}^r$. 
Although simpler, the approach is expected to yield a higher approximation error on average (over $\vec q$), as explained below.

\subsection{Example}
\label{pc:example}

We consider the following model-problem consisting in the $1D$ transient diffusion equation,
\begin{equation}
	\frac{\partial U}{\partial t} = \frac{\partial}{\partial x} \left( \nu \frac{\partial U}{\partial x} \right), \label{eq:diffusion_eq}
\end{equation}
where the diffusivity $\nu$ is a stochastic field. Eq.~\eqref{eq:diffusion_eq} is solved for $t\in [0,T]$, in the unit domain $D=[0,1]$, and with deterministic 
boundary conditions $U(x=0,t)=-1$, $U(x=1,t)=1$, and homogeneous initial condition $U(x,t=0)=0$.
We consider a log-normal stochastic diffusivity field of the form, 
\begin{equation}
\nu = \nu_0 + \exp(M), \label{diff:def}
\end{equation}
$M$ is a (centered) Gaussian process with uncertain covariance function ${\cal C}(\vec q)$. With $\nu_0>0$ the diffusivity is bounded away from 0 which ensures the well-posedness of the problem. In the computations we set $\nu_0=0.1$. In addition, we re-use the settings of Section~\ref{sec:exKL} with Gaussian covariance function 
having an uncertain length-scale $l$ with uniform  distribution in $[0.1,1]$ and fixed variance $\sigma_f^2=0.5$.
For the solution of Eq.~\eqref{eq:diffusion_eq} we use a classical ${\rm P1}$-finite element method (continuous piecewise linear approximation) for the spatial 
discretization, with a second order implicit time-integration scheme. 

To investigate the error introduced by approximating $M \mapsto U$ by the PC map $\hat M_K \mapsto \tilde{U}$, we define the following error measures on the model 
problem solution. We first define the relative local error $\epsilon_U(o,K,\vec q)$ as
\begin{equation}
		\epsilon^2_U(o,K,{\vec q}) \doteq \frac{\| U(M(\vec q)) - \tilde U(\xxi(\cdot,\vec q))\|^2_{L_2(\Omega,Y)}}{\| U(M(\vec q)) \|^2_{L_2(\Omega,Y)}},
\end{equation}
where 
\begin{equation}
	\| V\|^2_{L_2(\Omega,Y)} = \mathbb E\left[ \int_0^T \| V(x,t) \|^2_{L_2(D)} dt \right].
\end{equation}
This error measure incorporates the effects of several approximations: the approximation of $M(\vec q)$ on the $K$-dimensional reference subspace, the truncation of 
the PC expansion to finite order $o$, and the spatial and time discretization errors inherent in the numerical resolution of the model problem. 
Because the PC surrogate will be used in place of solving numerically the model problem (given $\vec \eta$ and $\vec q$), we should not be concerned with the 
spatial and time discretization errors, and rather use for $U(M(\vec q))$ its discrete counterpart, provided that the same spatial and time discretizations are used. 
For the tests presented in this section, we use a uniform mesh with 56 elements and a fixed time-step $\Delta t = 10^{-4}$. 
These discretization parameters were selected to ensure that the error measurements reported below are dominated by the $K$ and $o$-order truncation effects.
Doing so, the local error $\epsilon^2_U$ can be estimated by means of Monte Carlo average proceeding as follows. For a sample of $\vec q$, a) we generate a sample of the 
Gaussian process $M(\vec q)$ on the finite-element mesh and solve the corresponding deterministic diffusion problem for the sample of $U(M(\vec q))$; 
b) we project $M(\vec q)$ on the KL subspace, to obtain the KL coordinates $\vec\eta$ which are further translated to $\xxi(\vec \eta, \vec q)$, and the PC approximation 
$\tilde U(\xxi)$ is evaluated (see Eqs.~\ref{chgBasis:PC1}-\ref{chgBasis:PC}); c) we compute $\|U(M(\vec q))-\tilde U(\xxi)\|^2_{L_2(D\times T)}$ for the sample. 
Further, we set $T=0.05$ in order to focus the error measure in the transient period.

\begin{figure}[hbt]
\centering
\begin{tabular}{cc}
\includegraphics[width=0.45\textwidth]{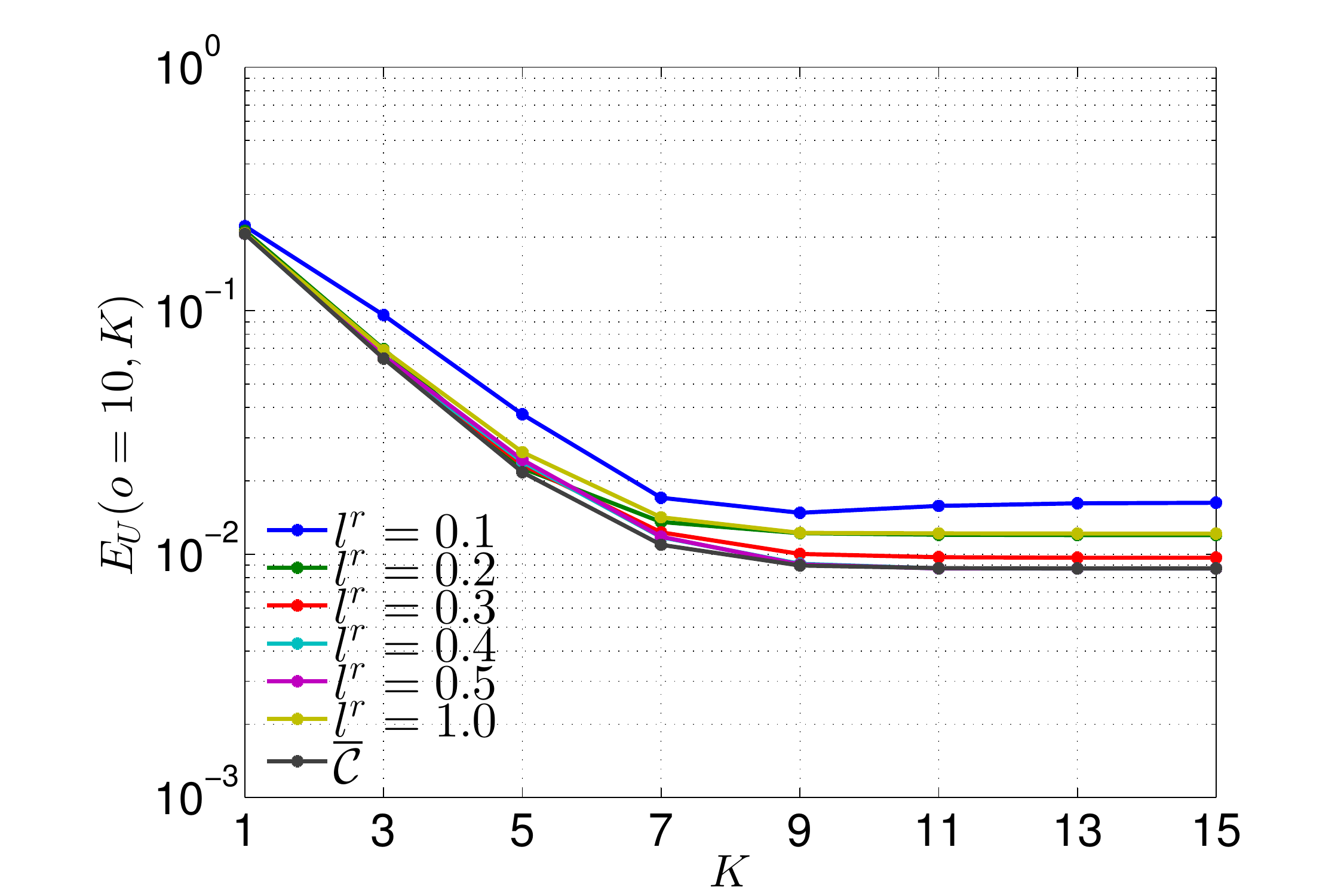}  &
\includegraphics[width=0.45\textwidth]{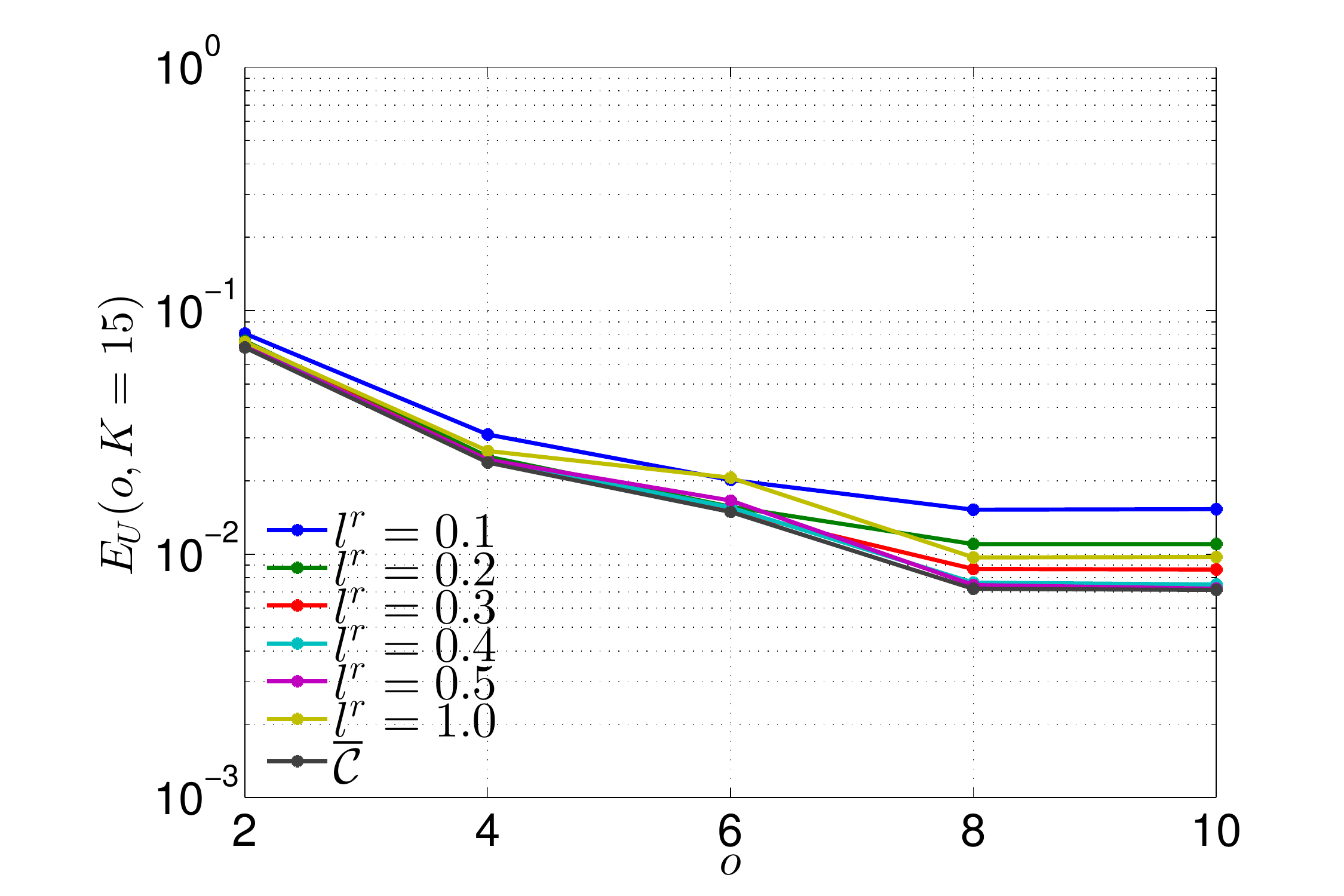}  
\end{tabular}
\caption{
Global error $E_U(o,K)$ of the PC approximation $\tilde U$ of the diffusion model problem solution.
(Left) The plot shows the dependence of the error with $K$ using a PC order $o=10$. (Right) The plot is for different $o$ and $K=15$.
The curves correspond to different definitions of the reference covariance function ${\cal C}^r$ = ${\cal C}(l^r)$ with $l^r$ as indicated or 
the $\vec q$-averaged covariance function $\overline{\cal C}$. }
\label{fig:errorU}
\end{figure}

Similarly, the local error can be $\vec q$-averaged to yield the relative global error counterpart:
\begin{equation}
	E^2_U(o,K) \doteq \frac{\dotsint \| U(M(\vec q)) - \tilde U(\xxi(\cdot,\vec q))\|^2_{L_2(\Omega,Y)} p_q(\vec q) d\vec q}
	{\dotsint \| U(M(\vec q))\|^2_{L_2(\Omega,Y)} p_q(\vec q) d\vec q}.\label{eq:errorU}
\end{equation}
Figure~\ref{fig:errorU} reports the global error for the present test problem. 
The left plot depicts $E_U(o,K)$ as a function of $K$ and for a PC order $o=10$.
Errors are shown for same selection of reference covariance functions ${\cal C}^r$ used in Section~\ref{sec:exKL}.
We observe that for all the selected reference covariance functions, the global error on $U$ stagnates for $K\gtrsim 9$. This indicates that 
the PC truncation becomes the dominant source of error for $K\gtrsim 9$. It is also seen that for all $K$ shown, the error is the lowest when using 
$\overline{\cal C}$ for reference covariance, as expected.

Further, when using ${\cal C}({\vec q}^r)$ as reference covariance function, the dependence of the global error on the reference 
length-scale $l^r$ is non monotonic, but presents a minimum around $l^r=0.4$.
This minimum can be explained by the competition of two effects. 
On the one hand, we have seen that increasing $l^r$ causes an increase in the approximation error of $M$, which translates in a larger 
approximation error on $U$. On the other hand, it can be shown that the lower $l^r$ the more the $\vec q$-marginal of $p_{\hat \eta}$ departs 
from the $K$-variates standard Gaussian distribution assumed for the construction of $\tilde U$, with increasing averaged approximation error on $U$ as a result. 
The right plot of Figure~\ref{fig:errorU} also depicts the global error, but now for a fixed number of KL modes, $K=15$, and increasing PC order $o\in[2,10]$. 
Again, curves are shown for the different reference covariance functions. Similarly to the previous results, the global error is seen to stagnate for 
$o\gtrsim 8$, indicating here that for larger $o$ the KL truncation error is dominant.
In addition, for all shown $o$, using $\overline{\cal C}$ for reference covariance function appears to be superior to the choices
${\cal C}(l^r)$, while the later choice again exhibits a non-monotonic dependence of the error with respect to $l^r$.

\begin{figure}[hbt]
\centering
\begin{tabular}{ccc}
\includegraphics[width=0.5\textwidth]{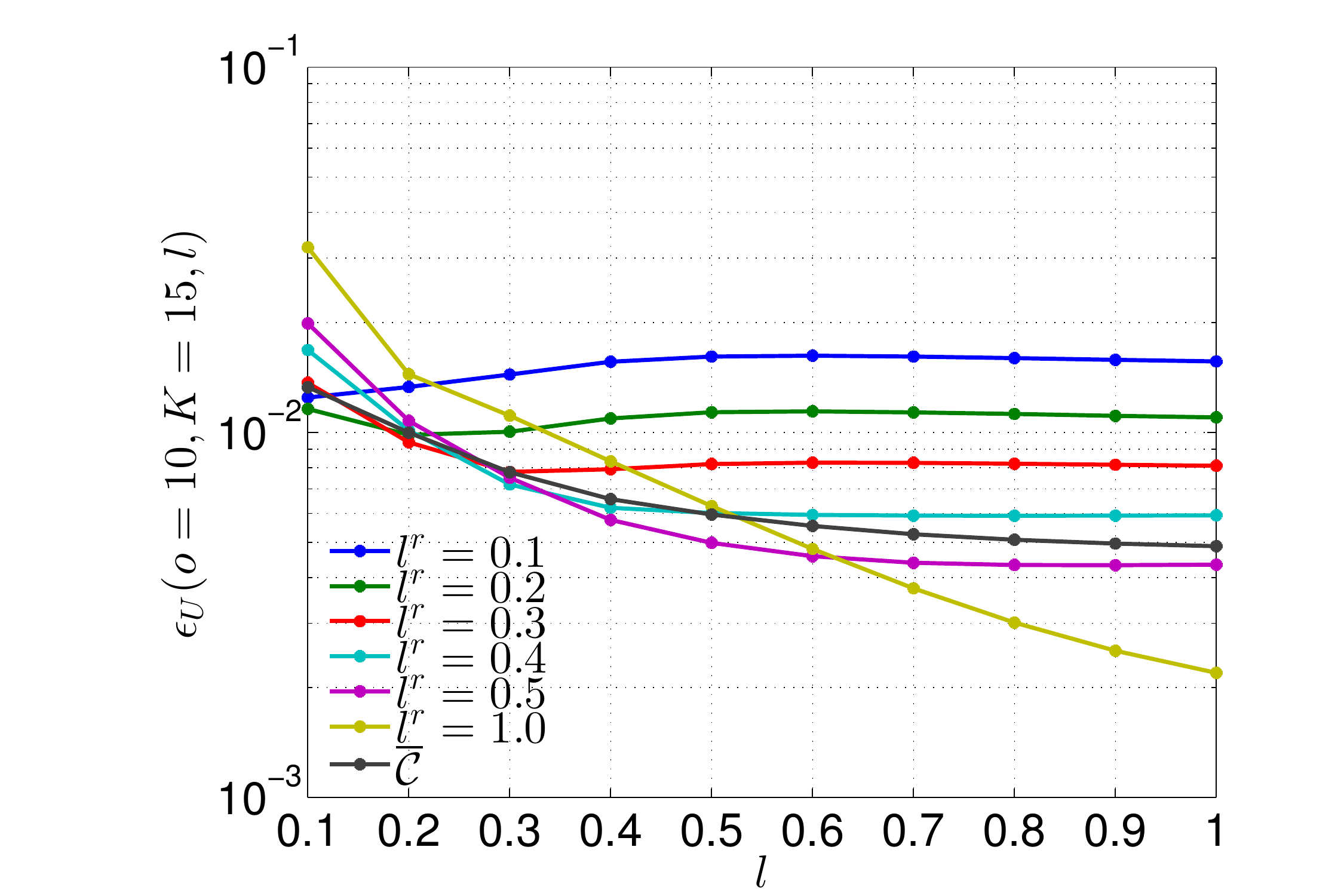}&
\includegraphics[width=0.5\textwidth]{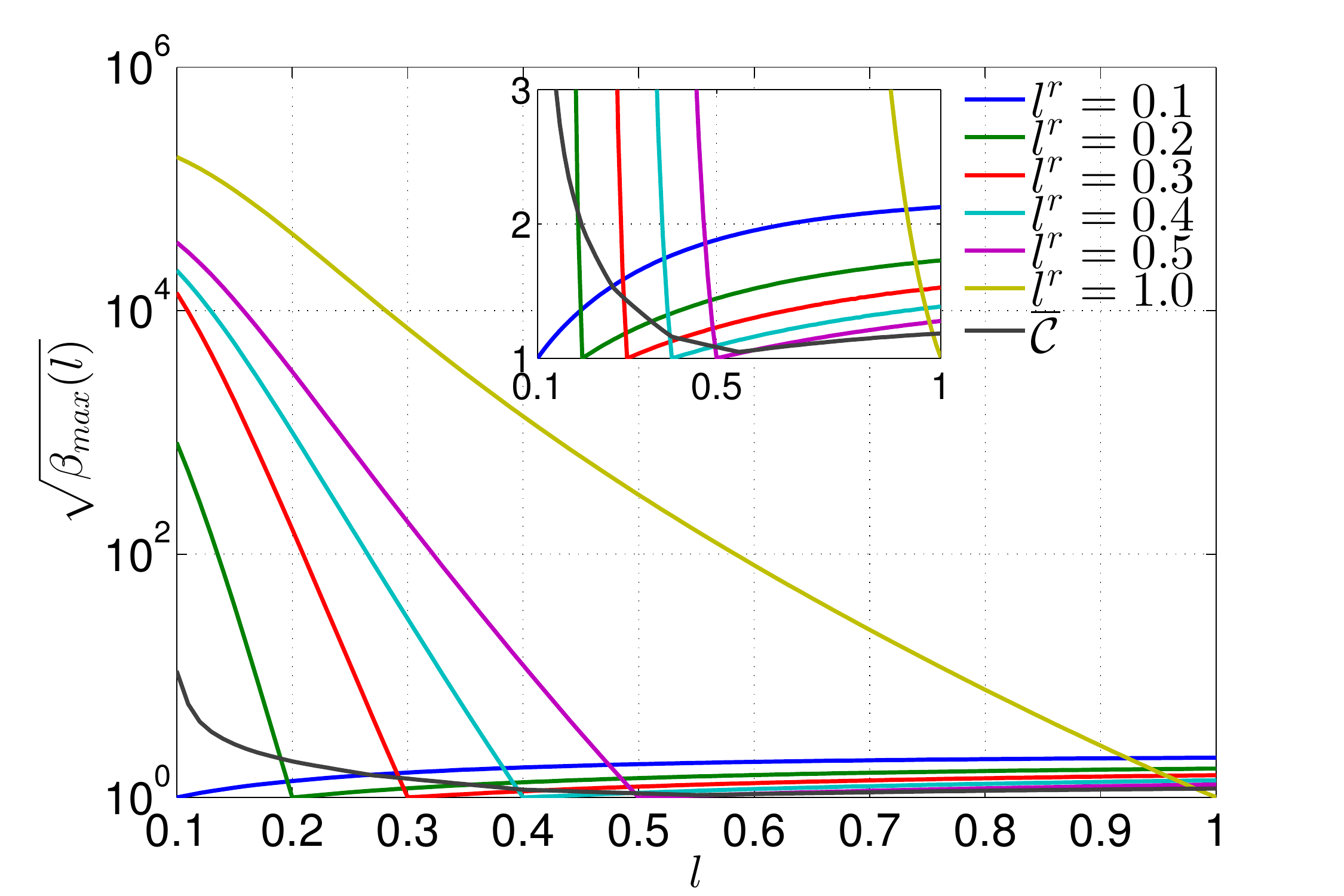}
\end{tabular}
\caption{(Left) Local approximation error $\epsilon_U(o,K,l)$, for $o=10$ and $K=15$. 
(Right) (log-scale) and inset (linear-scale): dependence on $l$ of the maximal stretching rate $\sqrt{\beta_{\rm max}(l)}$ induced by the coordinate 
transformation $\hat{\cal B}(\vec q)$. 
The curves correspond to different definitions of the reference covariance function ${\cal C}^r$ = ${\cal C}(l^r)$ with $l^r$ as indicated or 
the $\vec q$-averaged covariance function $\overline{\cal C}$.  
}
\label{fig:errorUks}
\end{figure}

Figure~\ref{fig:errorUks} presents in the left plot the normalized local error $\epsilon_U(o,K, l)$ for the case of $K=15$ and PC order $o=10$.
As mentioned previously, the local error combines the effects of approximating $M$ by $\hat M_K$, which has been reported in the right plot of 
Figure~\ref{fig:errorinM}, and the PC truncation error. Focusing first on the cases where ${\cal C}({\vec q}^r)$ is used as reference covariance,
we observe a more complex behavior of the local error with $l$, depending on the selected reference length-scale $l^r$. 
Specifically, the local error at some $l$ is always the lowest for the reference length-scale $l^r$ the closest to $l$. 
This is expected, as using $l^r$ is the optimal choice, given $K$ and $o$, to achieve the lower error at $l=l^r$. 
For ${\cal C}^r = \overline{\cal C}$, which ensures by construction the best compromise
over the $\vec q$-range, the local error remains below 2\%  over the whole range of hyper-parameters.
Further, using $l^r>0.2$, the local error first monotonically decreases with $l$ and then stagnates (except for $l^r=1$ where stagnation is not achieved).
 
Contrary to the local approximation error on $M$, the stagnation with $l\rightarrow 1$ occurs at an error level that strongly depends on $l^r$. 
This seems surprising as we have seen (right plot of Figure~\ref{fig:errorinM}) that for $l^r>0.2$ the approximation error on the process goes 
to zero as $l\rightarrow 1$, such that we could have expected an essentially constant local error $\epsilon_U$ for $0.2<l^r \lesssim l$, 
depending only on the PC expansion order $o$. 
But one has to take into account the mapping from $\vec \eta$ to $\vec \xi$ to understand the behavior of the local error. 
Specifically, the PC approximation is constructed to minimize the approximation error for the reference model problem based on 
$\hat M^{\rm PC}_K$ (or $K^{\rm PC}(\kappa)$) in Eq.~\eqref{def:RGPPC}, where the $\xi_k$'s are independent standard random variables. 
Therefore, the PC approximation aims at minimizing the error \emph{with respect to} the standard $K$-variates Gaussian measure. 
When querying the PC approximation for some specific hyper-parameters value $\vec q\ne {{\vec q}^r}$, 
$\xxi$ follows a conditional Gaussian distribution $p_\xi(\xxi|\vec q)$, induced by the transformation $\xxi = \hat{\cal B}(\vec q)\vec \eta$. 
In general, this conditional distribution differs from the standard Gaussian one, affecting the quality of the approximation depending on $\vec q$. 
To get better insight into this effect, we remark that the conditional density $p_\xi(\xxi|\vec q)$ is centered and Gaussian with covariance structure 
$\Sigma^2_\xxi(\vec q) = \hat{\cal B}^t(\vec q) \hat{\cal B}(\vec q)$. 

To measure the departure from the standard Gaussian multi-variates case, we present in the right plot of Figure~\ref{fig:errorUks} the largest 
eigen-value $\beta_{\rm max}(\vec q)$ of $\Sigma^2_\xxi(\vec q)$, as a function of $\vec q = \{ l \}$ and for the different reference covariance functions. 
$\sqrt{\beta_{\rm max}}$ measures of highest stretching rate induced by $\hat {\cal B}$.
For the results reported in Figure~\ref{fig:errorUks}, we used a thresholding parameter $\kappa = 10^{-12}$ in the definition of $\hat{\cal B}(\vec q)$. 
It is seen that when using ${\cal C}(l^r)$ for reference,
$\sqrt{\beta_{\rm max}}(l)$ increases exponentially fast with $l^r-l>0$, denoting a more and more stretched distribution for $\xxi(\vec \eta, l)$, along some direction, 
as $l$ decreases.
Interestingly, although the maximal stretching rate can reach values as high as $10^6$, its impact on the PC approximation error is clearly much less important.
The reason for the moderate sensitivity to coordinates stretching of the PC approximation error is that most of the stretching 
occurs along the directions associated with the lowest eigen-values $\lambda^r_k$, which have low to insignificant impacts on the model problem solution. 
In fact, our numerical experiments have demonstrated that the PC approximation error is essentially insensitive to $\kappa$, provided it is small enough. 
Indeed, a fast (exponential) decay of the successive KL modes' contributions to $U$ is expected for elliptic and parabolic model problems, as the effects of short-scale
fluctuations in the diffusivity field are filtered-out. 
However, coordinates stretching may yield robustness issues for other model types. 
Finally, it is seen that choosing $l^r$ equal to the minimal length-scale ($l=0.1$) yields a maximum stretching 
$\sqrt{\beta_{\rm max}}<3$ which is controlled over the whole range of $l$, while the case of $\overline{\cal C}$ yields a significant stretching 
(picking to $\approx 10$) around the minimal length-scale, but quickly decays with $l$ and remains close to $1$ (see left plot of Figure~\ref{fig:errorUks}). 
These findings confirm the appropriateness of the reference covariance function for the construction of the PC surrogate.

Concerning the PC approximation of the model-problem solution, we would like to stress the following points.
First the approach can be readily extended to alternative and more elaborated PC constructions methods, including non-intrusive ones; in particular 
considering adaptive techniques where the set of polynomials used in the PC expansion is determined as to minimize the approximation error, 
instead of proceeding from PC basis with uniform truncation order $o$, would clearly be beneficial, especially for problems involving high 
numbers of KL modes $K$ and requiring high polynomial order along certain $\xi_k$'s and not others. 
Second, the numerical tests have focused on length-scale uncertainty only, which is indeed the hardest source of uncertainty as it affects both the magnitude 
\emph{and shape} of the KL modes. In contrast, uncertainty in the process variance $\sigma_f^2$ in the Gaussian covariance family
only manifests itself in the magnitude of the eigen-values. Therefore, uncertainty in the pre-exponential factor $\sigma_f^2$ of the Gaussian covariance 
can be handled through either an additional dimension to the PC expansion, as performed in~\citep{MarzoukNajm2009,Tagade:2014}, or directly through our proposed 
change of coordinates approach based on the reference $\overline{\cal C}$, which amounts to take the averaged variance as the reference one. 
Similar to the problem for uncertain length-scale $l$, numerical tests (not shown)  have demonstrated that the $\vec q$-averaged definition of $\overline{\cal C}$ 
leads to globally lower errors in presence of variance uncertainty, compared to a definition of the reference based on some $q^r$. 
This has motivated the use of the $\vec q$-averaged definition of the reference covariance function in the remainder of the paper.

Finally, the PC expansion of the full model-problem solution has been considered here; there may be other situations were expansion of the full model-problem solution is not necessary. For instance, if the nature of the observations are known prior to constructing the PC expansion, the direct expansion of the model predictions 
$\vec u(\xxi)$ could be considered. 
If in addition the measurements have been performed, considering the direct PC expansion of the measurements to model-predictions discrepancy, 
$\Delta_d (\xxi) = \sum_{i=1}^\Nobs |d_i - u_i(\xxi)|^2$ (in the case of identically distributed additive noise) could be advantageous.

%%%%%%%%%%%%%%%%%%%%%%%%%%%%%%%%%%%%%%%%%%%%%%%%%%%
%%%%%%%%%%%%%%%%%%%%%%%%%%%%%%%%%%%%%%%%%%%%%%%%%%%%
%%%%%%%%%%%%%%%%%%%%%%%%%%%%%%%%%%%%%%%%%%%%%%%%%%%%%%%
\section{Application examples}
\label{sec:results}

In this section, we illustrate the benefit of considering prior Gaussian fields with parametrized covariance function in the inference of the diffusivity field in the transient diffusion problem introduced in Section~\ref{pc:example}.
We first present in Section~\ref{sec:setup} the inference problem, and introduce 3 cases that will serve to investigate the proposed method.
We also provide details on the exploitation the PC surrogate constructed in Section~\ref{sec:pc}, and on the PC accelerated formulation of the inference problem. 
For comparison purposes, we first solve in Section~\ref{sec:infnohyp} the Bayesian inference problem for a fixed covariance prior, that is without inferring the covariance hyper-parameters, using instead preassigned values. Then, in Section~\ref{sec:infhyp}, we considered the inference with hyper-parameters covariance and illustrate its advantage and behavior with respect to noise level, number of observations and surrogate polynomial order.

\subsection{Set-up of the inference problem}
\label{sec:setup}
The proposed method will be illustrated for the inference of a log-diffusivity field, using the transient diffusion model problem corresponding to Eq.~\eqref{eq:diffusion_eq}. 
To test the proposed method we consider three different log-diffusivity fields $m(x) = \log(\nu - \nu_0)$, to be inferred: 
\begin{itemize}
	\item Sinusoidal profile: $m^{\rm sin}(x) = \sin(2\pi x)$,
	\item Step function: $m^{\rm step}(x) = \begin{cases} -1/2, & x<0.5 \cr 1/2, &x\ge 0.5\end{cases}$,
	\item Random profile: $m^{\rm ran}(x)$ drawn at random from ${\cal GP} (0,{\cal C})$ where $\cal C$ is the Gaussian covariance with 
	length-scale $l=0.25$ and variance $\sigma_f^2=0.65$.
\end{itemize}

\begin{figure}[hbt]
\centering
\begin{tabular}{c}
\includegraphics[width=0.45\textwidth]{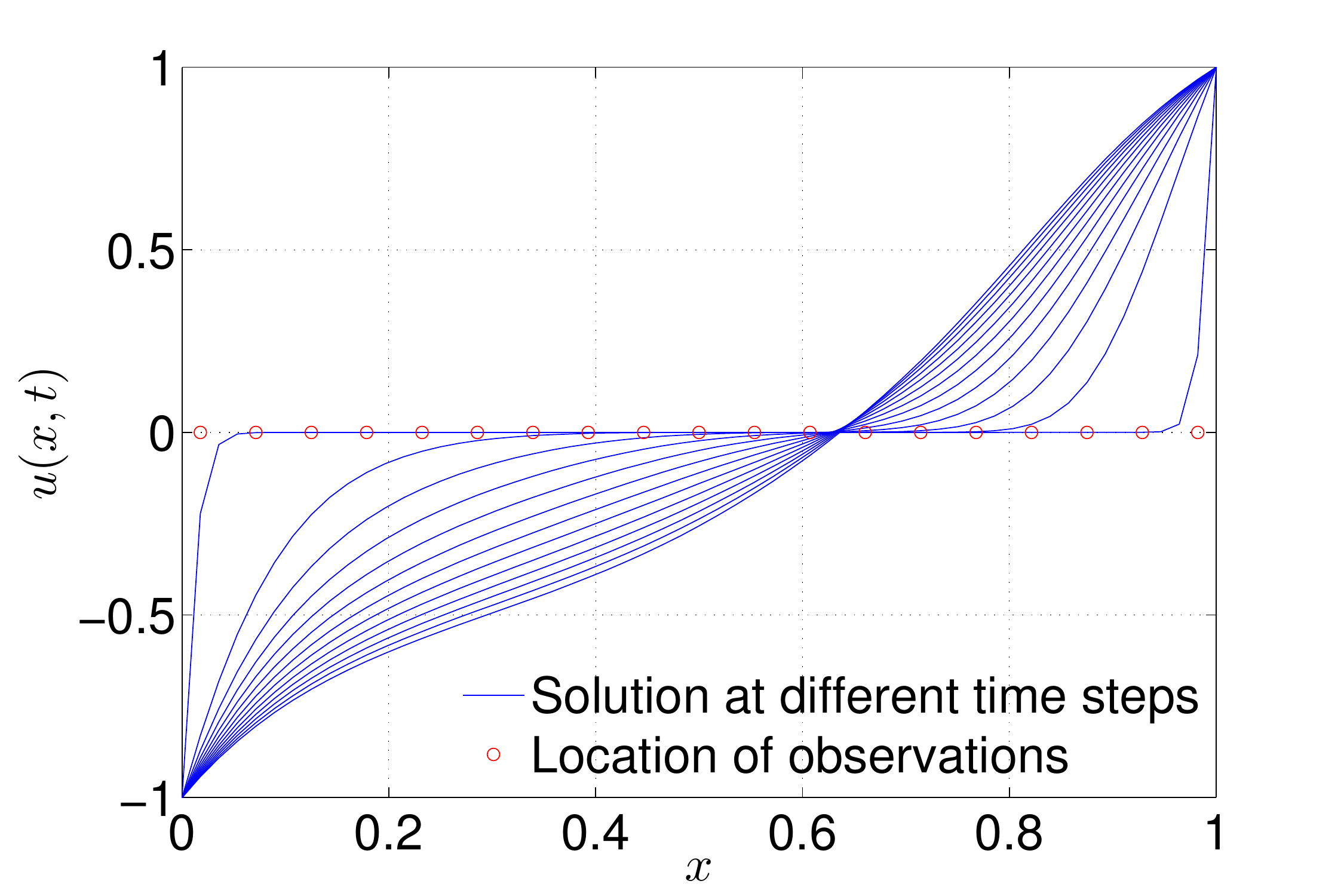} 
\end{tabular}
\caption{Illustration of inference problem for $m^{\rm sin}$. Plotted are the $n_x=19$ observation points and the solution of the diffusion equation with profile $m^{\rm sin}$ at the $n_t$ observation times.}
\label{fig:obs}
\end{figure}

The inferences are performed on sets of data, $\{d_i, i=1,\ldots,\Nobs\}$, consisting of noisy measurements of the solution to the diffusion equation for the three profiles. The measurements are taken at a set of $n_x$ spatial locations $x_i$ uniformly distributed inside $D=(0,1)$, and for $n_t$ times $t_i$ uniformly distributed in $(0,T)$. The total number of observations is then $\Nobs = n_x\times n_t$. The observations are synthetically generated by perturbing the respective model solutions for the 3 fields tested with a measurement noise $\epsilon_i$ randomly and independently drawn from the Gaussian distribution ${\mathcal N}(0,\sigma_\epsilon^2)$. 
To avoid the so-called inverse crime~\citep{Kaipio2007}, the solutions used to generate the observations are computed with a significantly finer spatial and temporal discretization than for the construction of the PC approximation. 
Unless otherwise specified, we use $n_x=19$, $n_t=13$ (so $\Nobs = 247$), with $T=0.05$ and a Gaussian noise with $\sigma_\epsilon^2=0.01$.
Figure~\ref{fig:obs} depicts the location of the observation points and the solution of the diffusion equation for $m^{\rm sin}$ at the different observation times.

For the inference, we consider in all cases the Gaussian prior $M \sim {\cal GP}(0,{\cal C(\vec q)})$, where the covariance function ${\cal C}(\vec q)$ has hyper-parameter ${\vec q}=\{ l,\sigma_f^2\}$. 
The prior is then fully characterized once we have selected the prior of the hyper-parameters. 
We choose a uniform prior for $l$ over the range $[l_{\min},l_{\max}]$, with as previously $l_{\min}=0.1$ and $l_{\max}=1$, and an 
inverse Gamma prior~\citep{Zabaras2005,GelmanEtAl:2004,gelman2006} for $\sigma_f^2$ with parameters $\alpha=3$ and $\beta=1$.
The prior of $\sigma_f^2$ thus has a long-tailed distribution with mean value $\beta/(\alpha-1)=0.5$ and variance $\beta^2/(\alpha-1)^2(\alpha-2)=0.25$. 
Note that the existence of the first moment of $\sigma_f^2$ is enough to ensure the existence of the average covariance function, and 
that $M(\vec \eta,\vec q) \in L_2(D,p_\eta,p_q)$ (because the modes in its KL decomposition scales with $\sqrt{\sigma_f^2}$). In contrast, expanding the diffusion equation solution with respect to both the KL coordinates $\vec \eta$ and hyper-parameter $\vec q$ (as proposed in~\citep{MarzoukNajm2009,Tagade:2014}) could be problematic since, to our knowledge, there is no standard orthogonal polynomial family for the inverse Gamma distribution function and the solution $U$ may not have second moment ($\exp(\sqrt(y))$ with $y\sim {\rm Inv}\Gamma(3,1)$ has unbounded second moment).
Using the notation of Section~\ref{sec:inference}, the prior of $\vec q$ is then
\begin{equation}
	p_q(\vec q) = p_q(l,\sigma^2_f) = \begin{cases}\displaystyle \frac{1}{|l_{\min}-l_{\max}|\Gamma(3)}(\sigma_f^2)^{-4} \exp \left(  -\frac{1}{\sigma_f^2}\right), & l\in[l_{\min},l_{\max}], \sigma_f^2>0\cr
	0, & \mbox{otherwise}. \end{cases}
\end{equation}
As for the noise hyper-parameter, we use the uninformative, improper, Jeffrey's prior 
\begin{equation}
	p_o(\sigma_o^2) \propto \frac 1 {\sigma_o^2}.
\end{equation}

\begin{figure}[htb]
\centering
\includegraphics[width=0.7\textwidth]{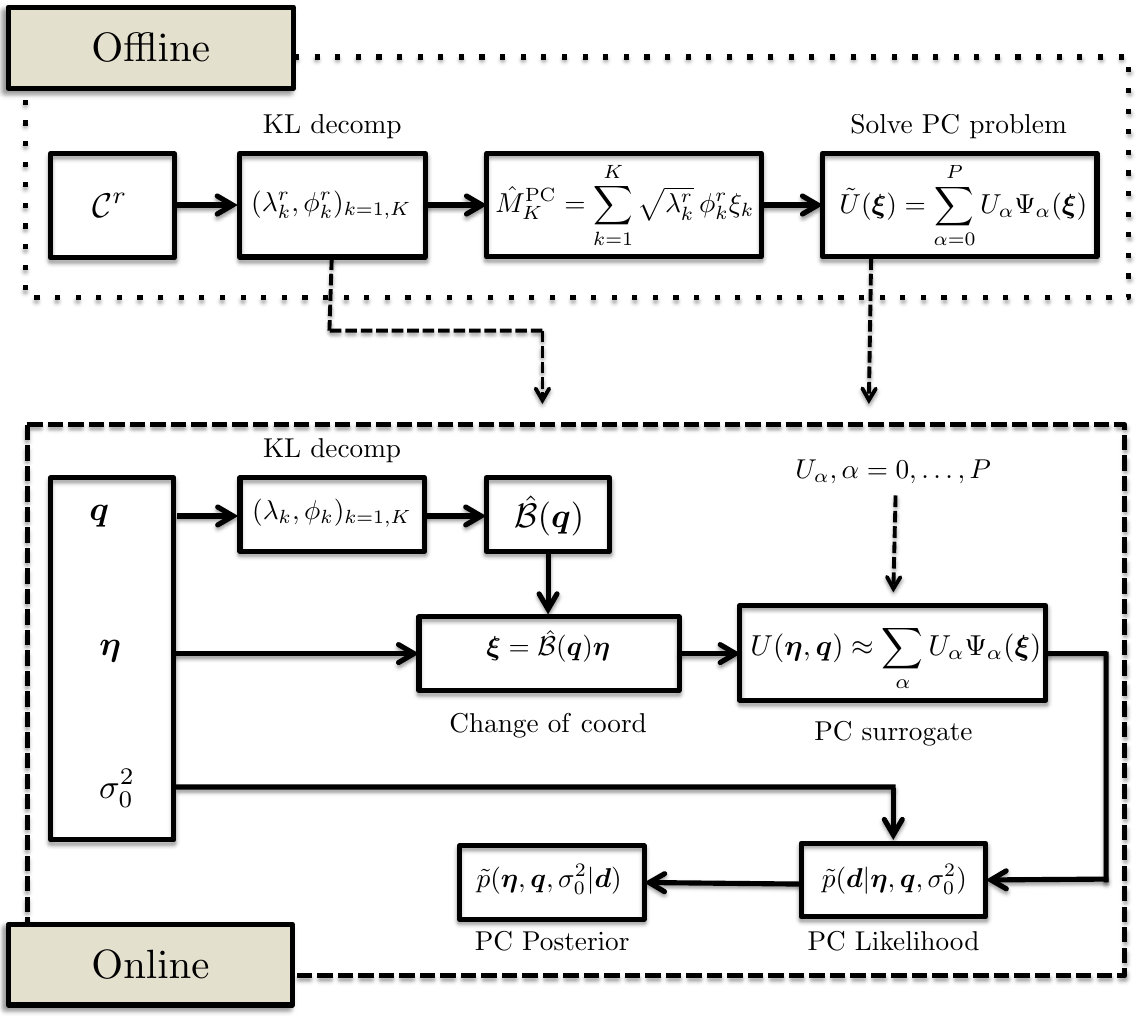}
\caption{Offline step (surrogate construction) of the accelerated MCMC sampler and Online step of the PC surrogate based evaluation of the posterior.}
\label{fig:flowchart2a}
\end{figure}

Having specified all priors, the determination of the Bayesian posterior $p(\vec \eta,\vec q, \sigma_o^2 | \vec d)$ requires the evaluation of the likelihood 
of the data $\vec d$ given $(\vec \eta, \vec q, \sigma_o^2)$. Instead of following the computational flow-chart presented in Figure~\ref{fig:flowchart1}, 
which would require the solution of a deterministic model problem for each new sample of $(\vec \eta, \vec q)$, we rely on coordinate transformation and PC approximation as introduced in the previous sections.
Following the findings of the previous section, the reference model problem is based on the stochastic process $\hat M^{\rm PC}_K$ corresponding to the $\vec q$-averaged covariance function $\overline{\cal C}$, whose KL decomposition is truncated to the first $K=15$ dominant modes. Solving this reference problem, we obtain the approximation 
$\tilde U(\xxi)=\sum_{\alpha=0}^PU_\alpha \Psi_\alpha(\xxi)$ of the reference model problem solution $U(\xxi)$. Unless stated otherwise we use in the following results a PC order $o=10$ with a spatial discretization involving 56 finite elements. From the approximate solution $\tilde U$, we can extract the PC approximations of the model predictions, $\tilde{\vec u}(\xxi)$, whose components are 
\begin{equation}
	\tilde u_i (\xxi) = \tilde U(x_i,t_i,\xxi) = \sum_{\alpha=0}^P U_\alpha(x_i,t_i)\Psi_\alpha(\xxi), \quad i=1,\dots, \Nobs.
\end{equation}
This constitutes the offline step of the proposed PC-accelerated sampler. 
Once the PC approximation has been determined, one can use $\tilde{\vec u}(\xxi(\vec \eta,\vec q))$  as a surrogate of the model predictions $\vec u(\vec \eta,\vec q)$ in the (online) computation of the likelihood:
\begin{equation}
	p(\vec d|\vec \eta, \vec q, \sigma_o^2) \approx \tilde p(\vec d|\vec \eta, \vec q, \sigma_o^2) \doteq \prod_{i=1}^{\Nobs} p_\epsilon (d_i - \tilde u_i(\xxi(\vec \eta, \vec q)),\sigma_o^2),
\end{equation}
where $\xxi(\vec \eta,\vec q)$ is given by Eq.~\eqref{chgBasis:PC1} and $p_\epsilon$ is defined in Eq.~\eqref{eq:likelihood1}. For the actual definition of the coordinate transformation $\hat{\cal B}(\vec q)$ in Eq.~\eqref{chgBasis:PC}, we set $\kappa=0$ since $\lambda^r_{k\le15}/\lambda^r_1$ remains large enough for the present settings.
Finally, multiplying by the hyper-parameter priors and prior of $\vec \eta$, one obtains (up to a constant normalization factor) the approximation $\tilde p(\vec \eta,\vec q,\sigma_o^2|\vec d)$  of the posterior distribution. 
The computational structure for the change of coordinates method and PC acceleration is schematically illustrated in Figure~\ref{fig:flowchart2a}, distinguishing 
between offline and online steps.
The online step is imbedded in an adaptive Metropolis-Hasting algorithm to generate samples of $(\vec \eta,\vec q,\sigma_o^2)$ following the posterior density.

\subsection{Inference with fixed covariance parameters}
\label{sec:infnohyp}
For comparison purposes, the Bayesian inference problems are first solved for the case of Gaussian prior covariance function having pre-assigned parameters $l = 0.5$ and $\sigma_f^2 = 0.5$.
The problem therefore consists in inferring only the $15$ coordinates $\vec \eta$ and the noise hyper-parameter $\sigma_o^2$.
Note that in this case $\hat{\cal B}$ is the identity, as the PC approximation is based on the prior process with pre-assigned covariance function.

\begin{figure}[hbt]
\centering
\begin{tabular}{ccc}
\includegraphics[width=0.3\textwidth]{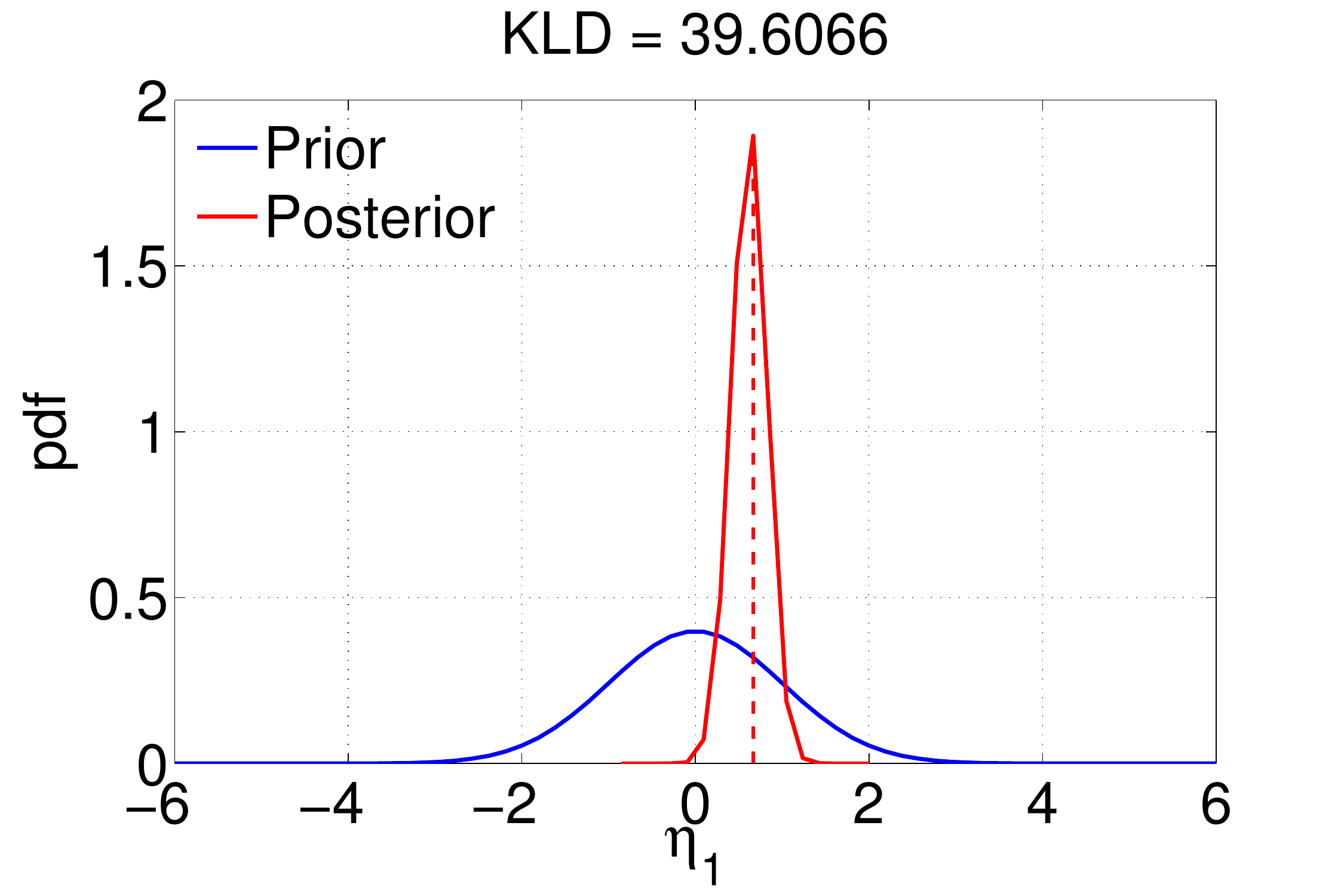}  &
\includegraphics[width=0.3\textwidth]{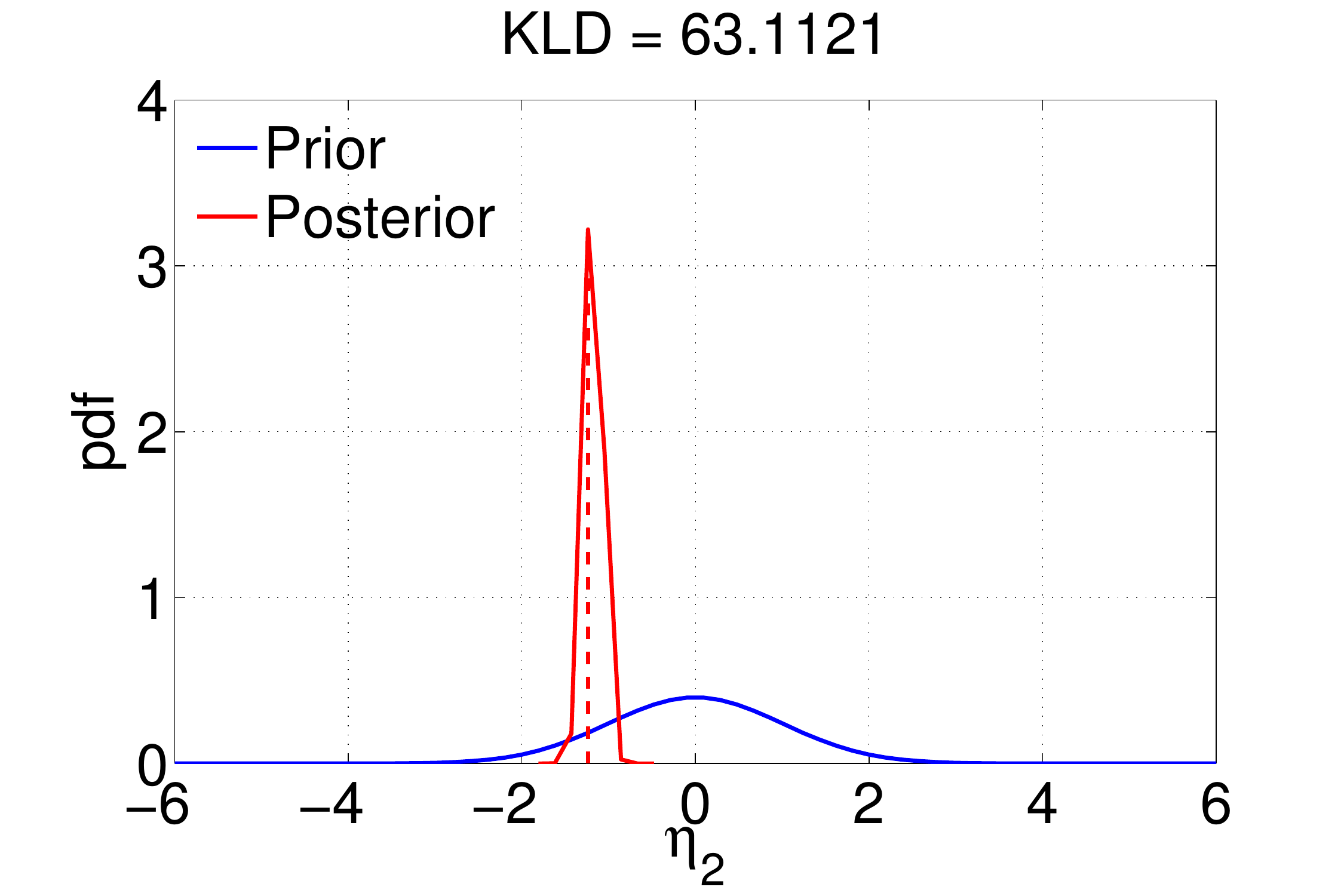}  &
\includegraphics[width=0.3\textwidth]{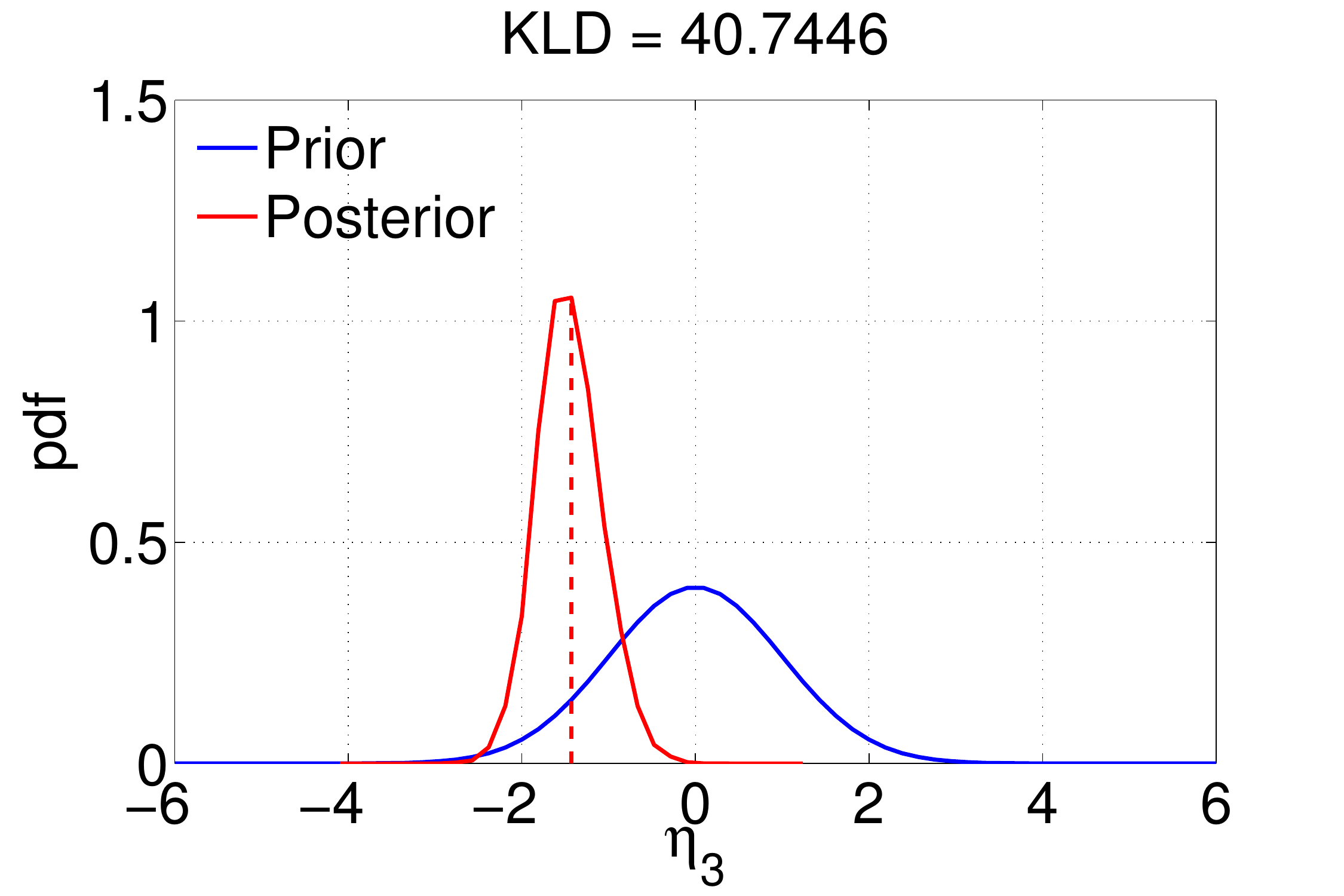}  \\
\includegraphics[width=0.3\textwidth]{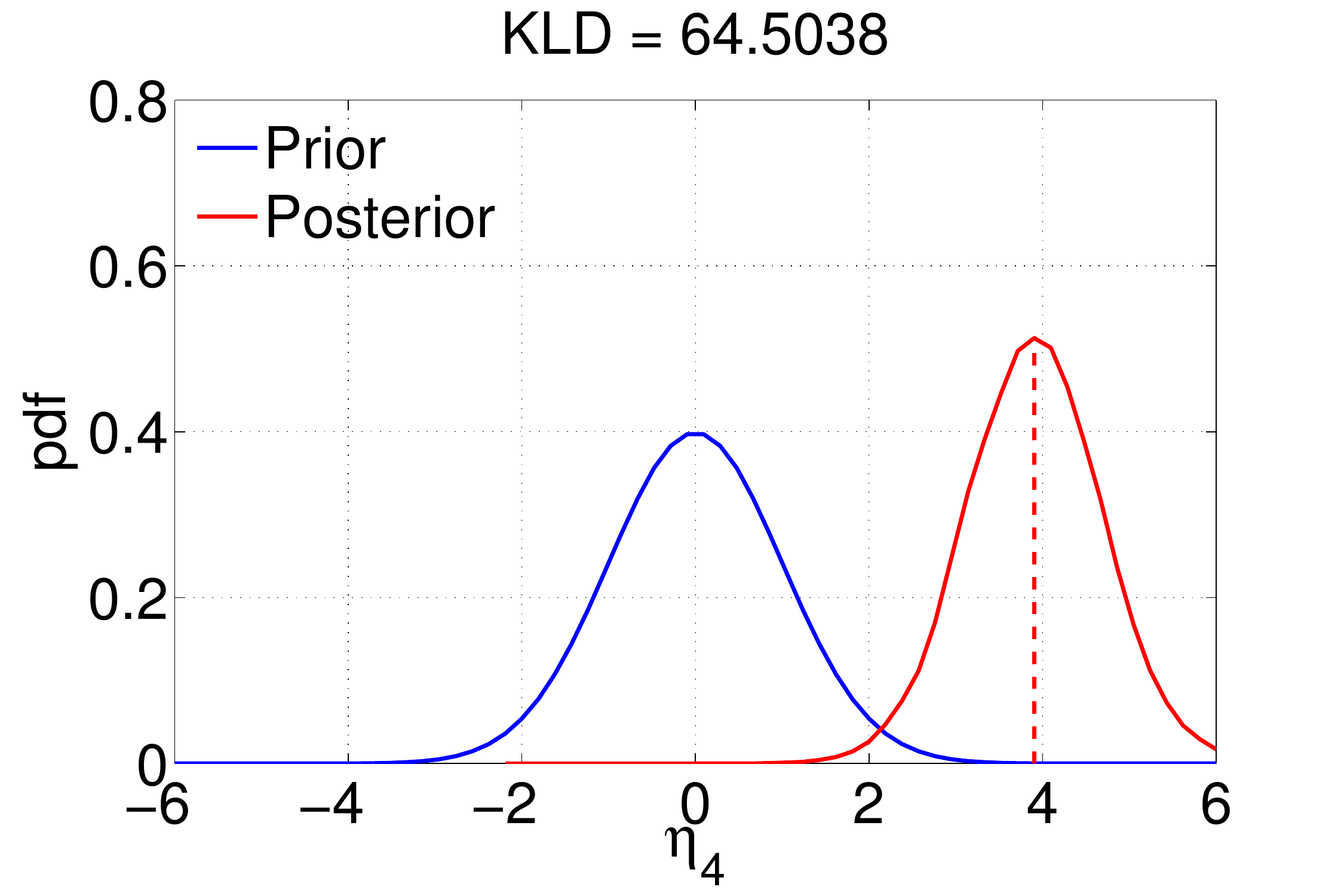}  &
\includegraphics[width=0.3\textwidth]{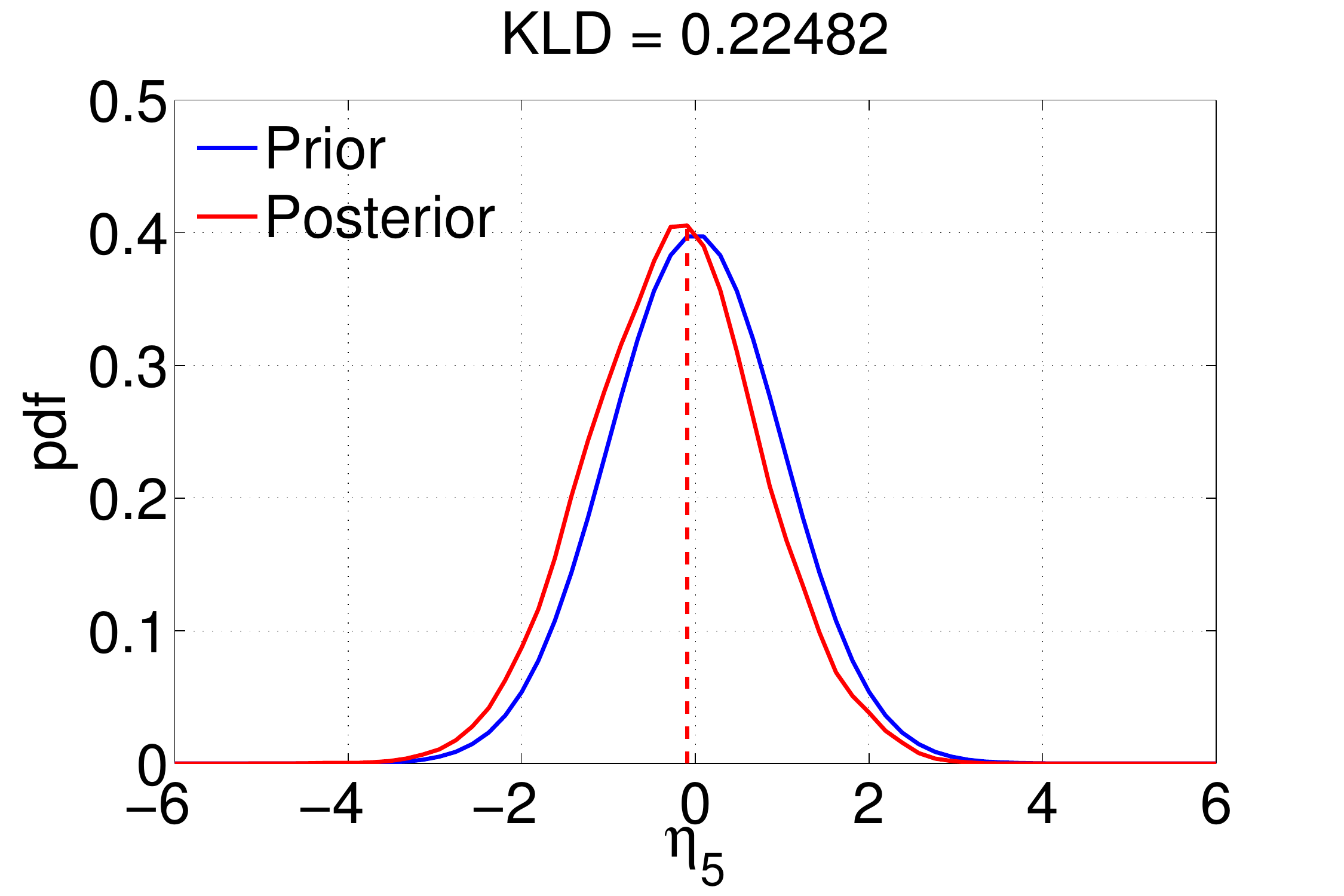}  &
\includegraphics[width=0.3\textwidth]{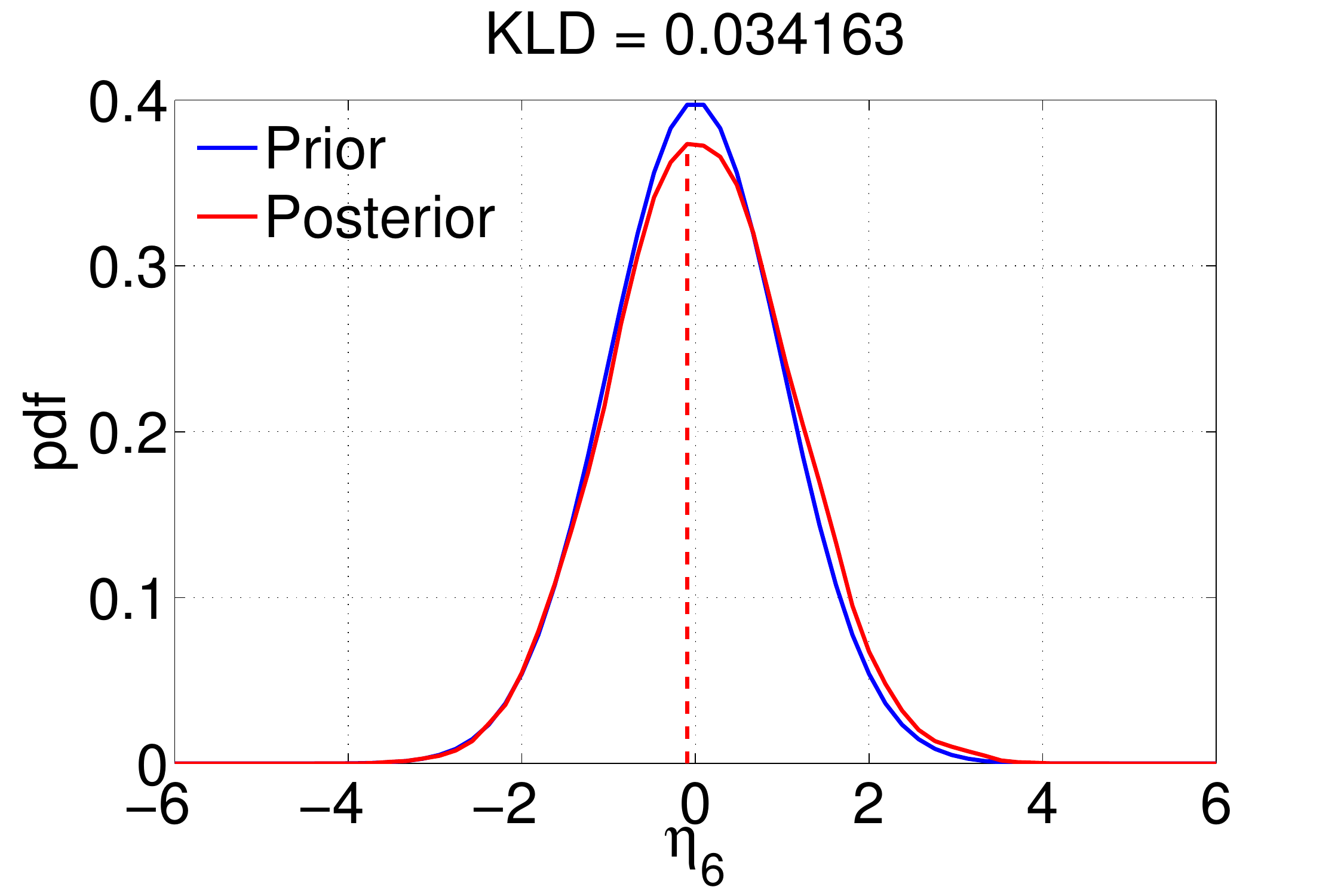}  \\
\includegraphics[width=0.3\textwidth]{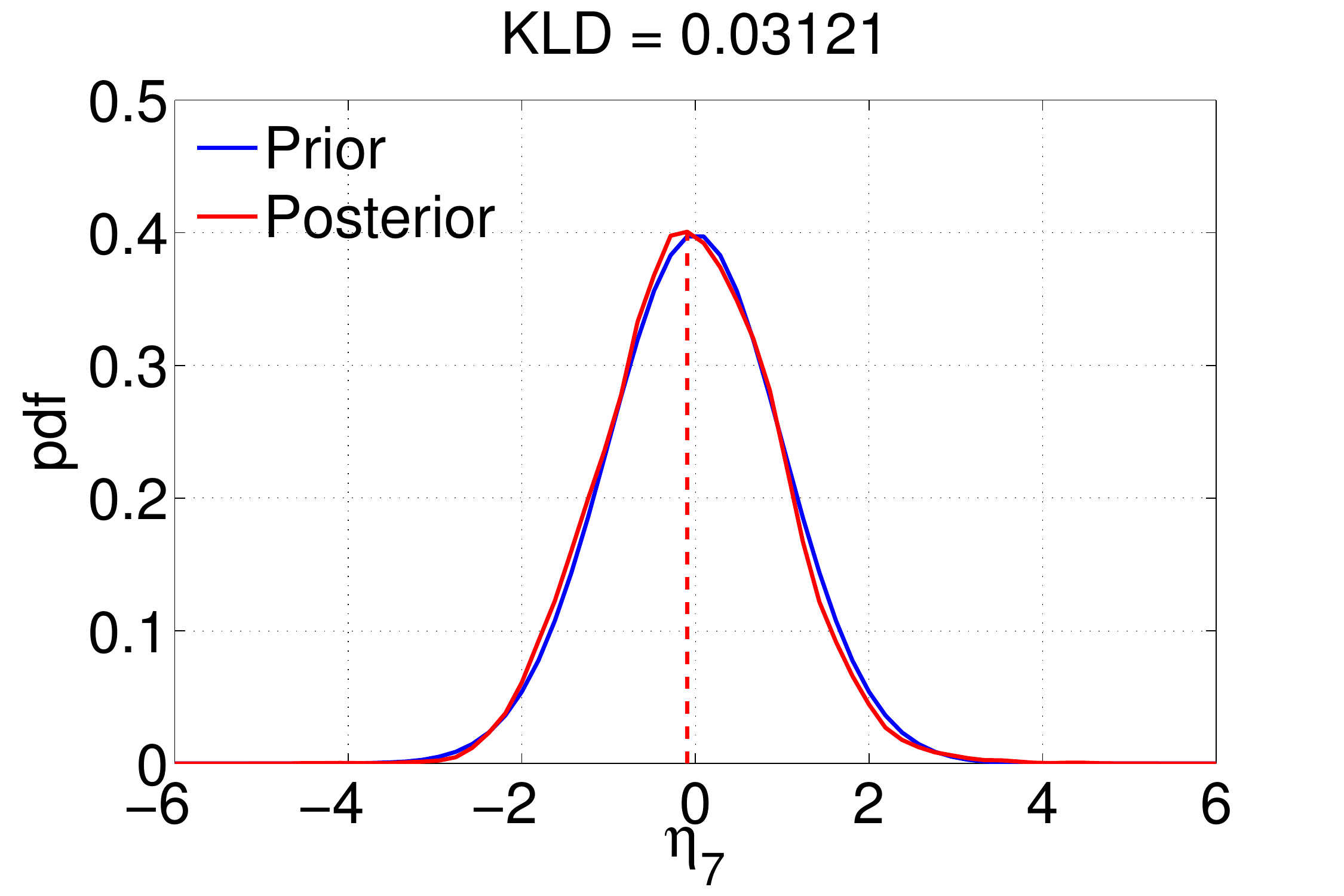}  &
\includegraphics[width=0.3\textwidth]{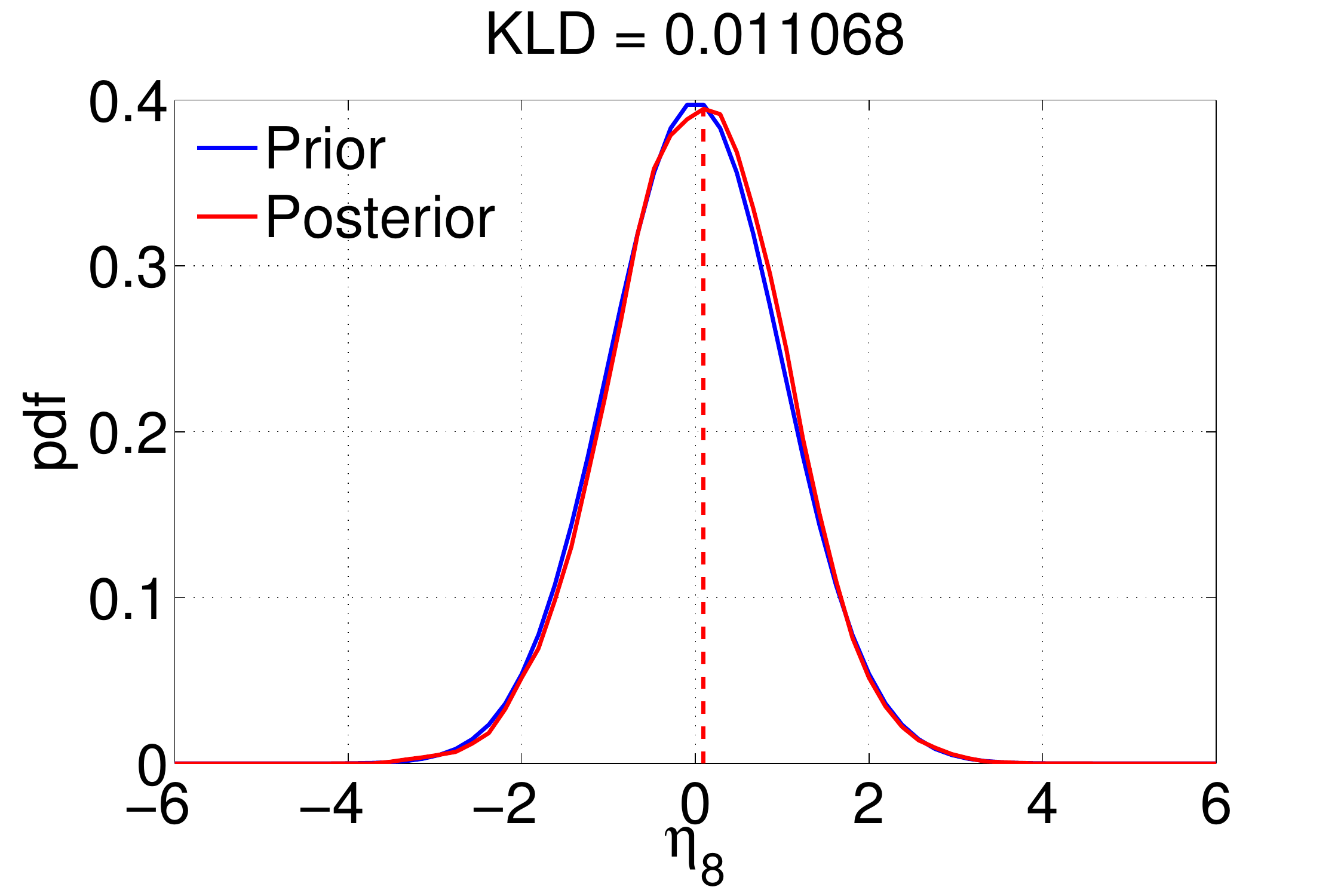}  &
\includegraphics[width=0.3\textwidth]{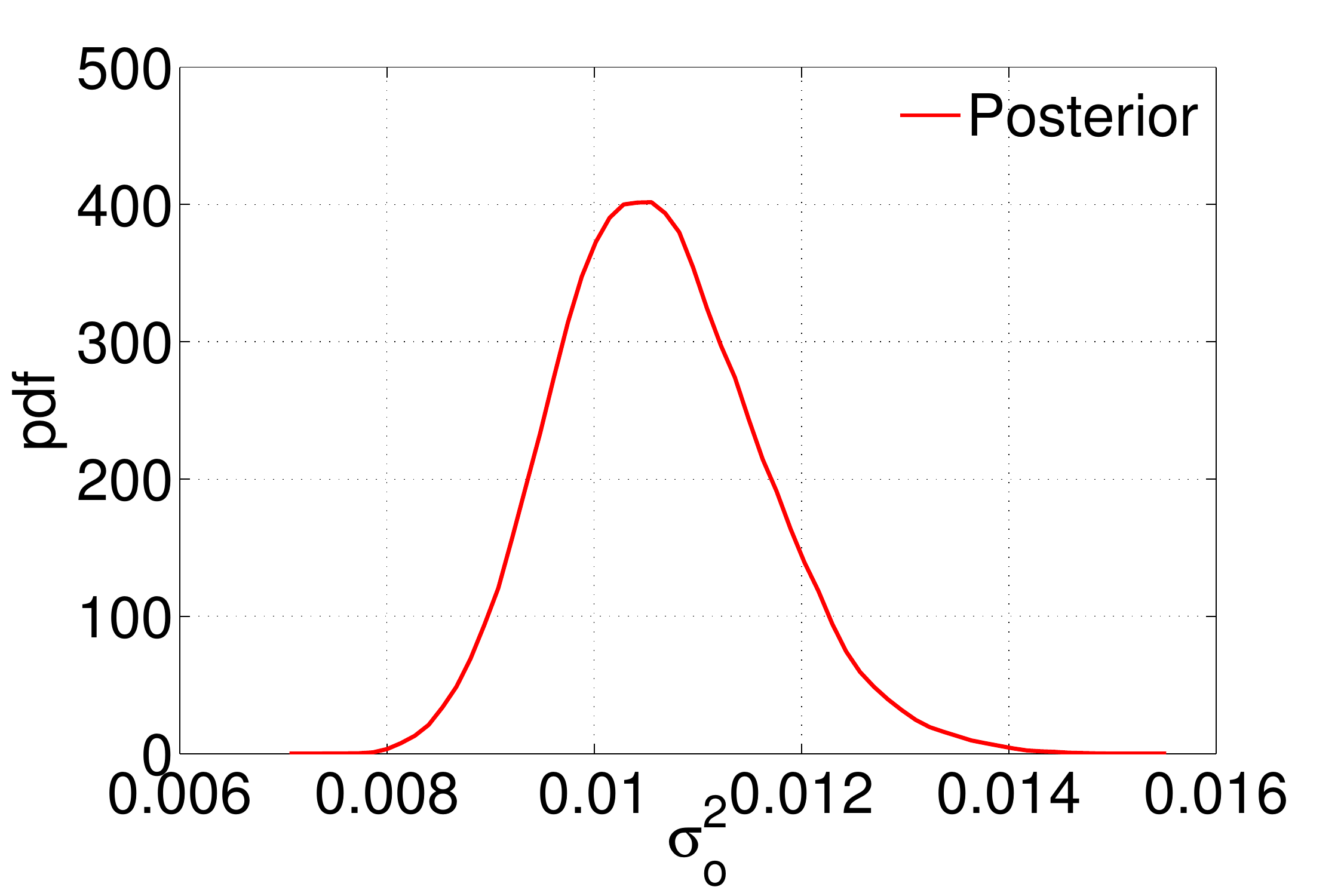}  
\end{tabular}
\caption{Comparison of the priors and (marginal) posteriors of the first 8 KL coordinates $\eta_k$ and noise hyper-parameter $\sigma^2_o$ (posterior only) for the inference of $m^{\rm sin}$ without using covariance hyper-parameters (a Gaussian covariance with $l=0.5$ and $\sigma_f^2=0.5$ is assumed). 
The corresponding Kullback-Leibler Divergences (KLD) for the KL coordinates are also indicated on top of each plot.}
\label{fig:pdfnhyp}
\end{figure}

A total number of $2.5 \times 10^5$ MCMC steps were deemed necessary for the pre-assigned hyper-parameters case to properly explore the posterior. 
The resulting chain of the KL coordinates were observed to be well-mixed (not shown).
The marginal posteriors, estimated using a standard Kernel Density Estimation (KDE) method~\citep{Silverman1986,Parzen1962},
of the first 8 KL coordinates for the inference of $m^{\rm sin}$ are shown in Figure~\ref{fig:pdfnhyp}. 
These posteriors are compared with their respective priors (standard Gaussian distributions). 
We notice that only the first 4 coordinates $\eta_k$ show significant improvement in their posterior distributions. 
This improvement can be quantified using the Kullback-Leibler Divergence (KLD) which is a statistical measure that quantifies 
the distance between two probability distributions $p$ and $q$~\citep{Kullback1959}, defined according to:
\begin{equation}
KLD(p,q) = \int_{-\infty}^\infty p(x) \ln{\frac{p(x)}{q(x)}dx}
\end{equation} 
Here we calculate the KLD between the prior and the (marginal) posterior of each KL coefficient $\eta_k$.
The KLD is indicated on top of each plot and quantifies the information gain from the observations, which is found significant only for the 
first 4 KL coordinates. 
Figure~\ref{fig:pdfnhyp} also shows the posterior of the noise variance hyper-parameter (bottom right plot), $\sigma^2_o$, which exhibits 
a Maximum A Posteriori (MAP) value close to the value used to generate the data, $\sigma_\epsilon^2=0.01$.
The similar findings are reported for the inferences of $m^{\rm step}$ and $m^{\rm ran}$ (results not shown for brevity).

\begin{figure}[hbt]
\centering
\begin{tabular}{ccc}
\includegraphics[width=0.35\textwidth]{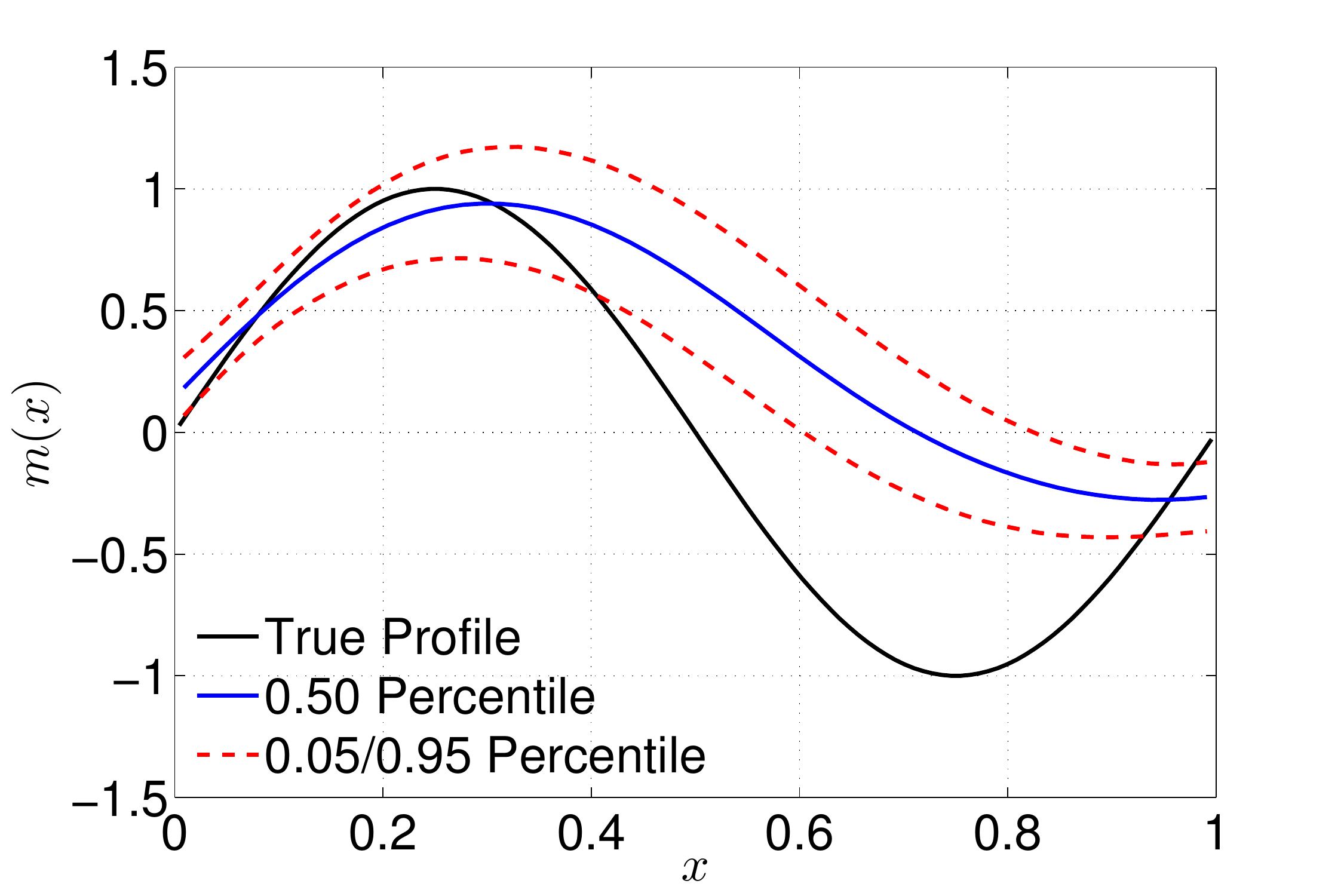}  &
\includegraphics[width=0.35\textwidth]{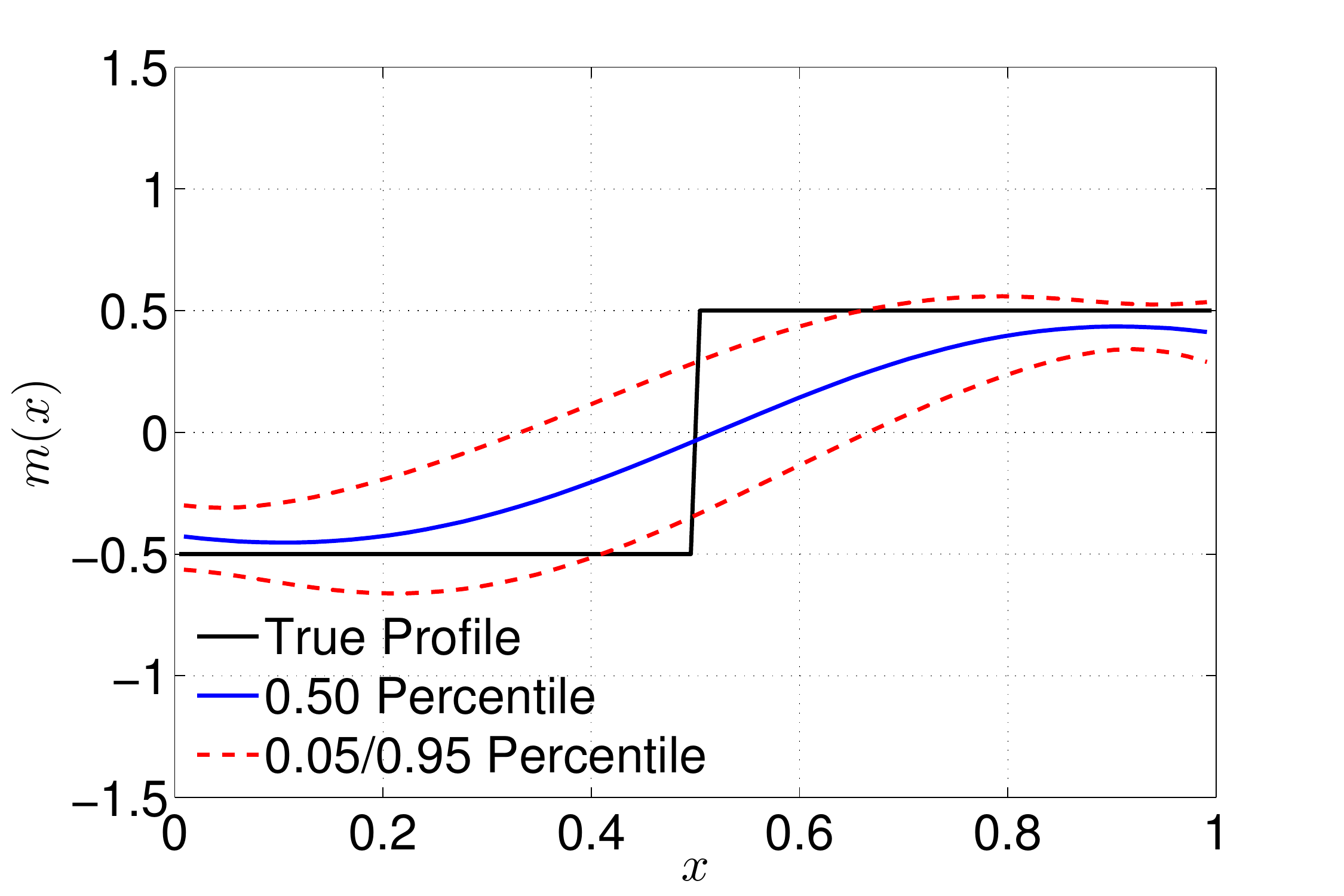} &
\includegraphics[width=0.35\textwidth]{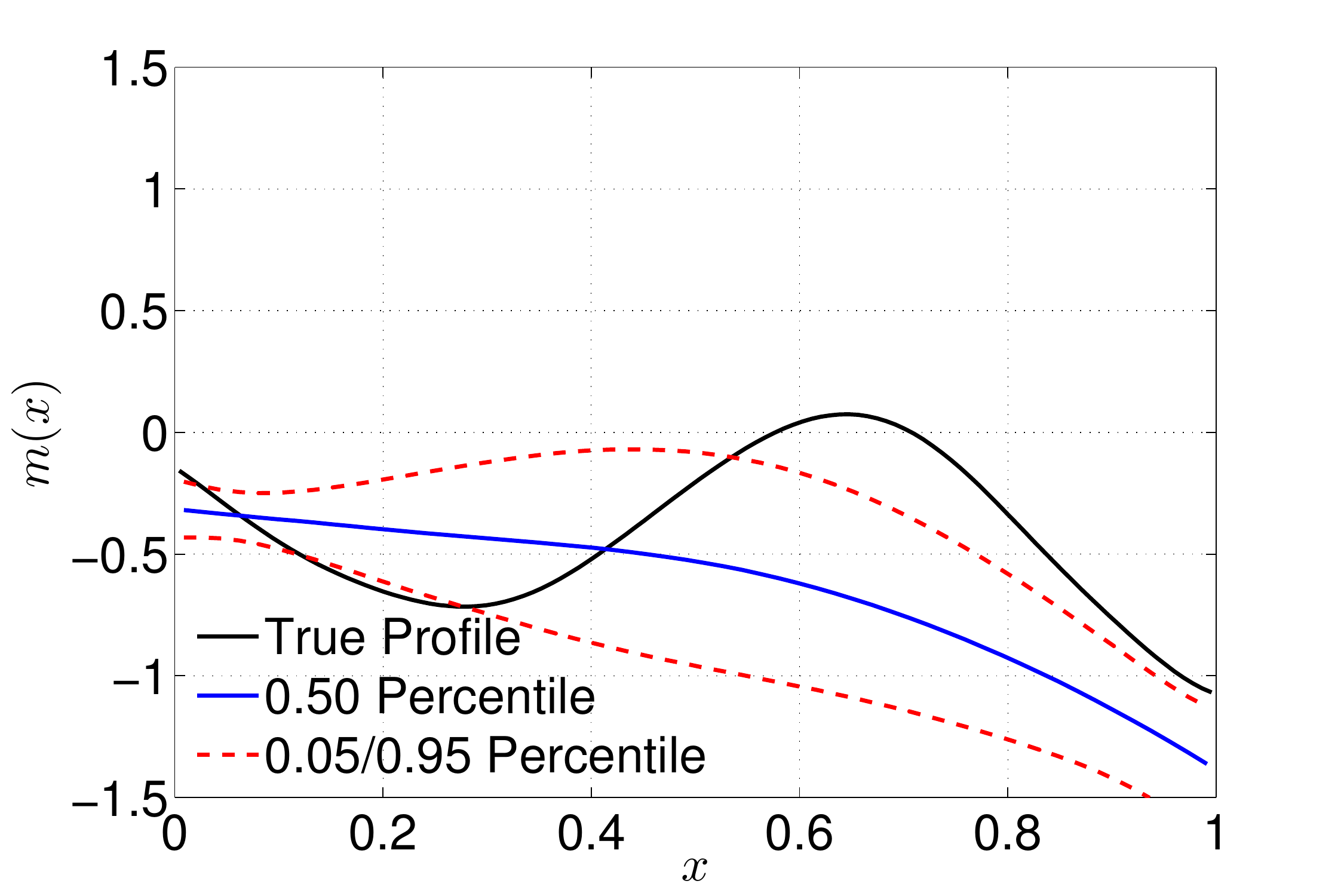} 
\end{tabular}
\caption{Comparison of the posterior of $m(x)$ with the true profile, for the cases of $m^{\rm sin}$, $m^{\rm step}$ and $m^{\rm ran}$ (from left to right).
The inferences use a fixed Gaussian covariance function with $l=0.5$ and $\sigma_f^2=0.5$. Shown are the median, $5\%$ and $95\%$ quantiles of the posterior and true profile.}
\label{fig:prfnhyp}
\end{figure}

To better analyze the quality of the inferred fields, we report in Figure~\ref{fig:prfnhyp}, for the 3 test cases, the median, $5\%$ and $95\%$ quantiles values of the posteriors of the inferred field $m(x)$. These statistical characterizations of $m$ are also compared with the true profiles.
For the inference of $m^{\rm sin}$ we notice that the $5\%$ to $95\%$ quantiles range does not contain the true profile for a large set of $x$. 
This mismatch can be attributed to the pre-assigned hyper-parameter values that are not suitable. 
The same observation can be made for the case of $m^{\rm ran}$.  
In contrast, $m^{\rm step}$ is nearly everywhere within the $5\%$-$95\%$ quantiles range of the inferred profile $m$.

\subsection{Inference with covariance hyper-parameters}
\label{sec:infhyp}
\begin{figure}[hbt]
\centering
\begin{tabular}{cc}
\includegraphics[width=0.4\textwidth]{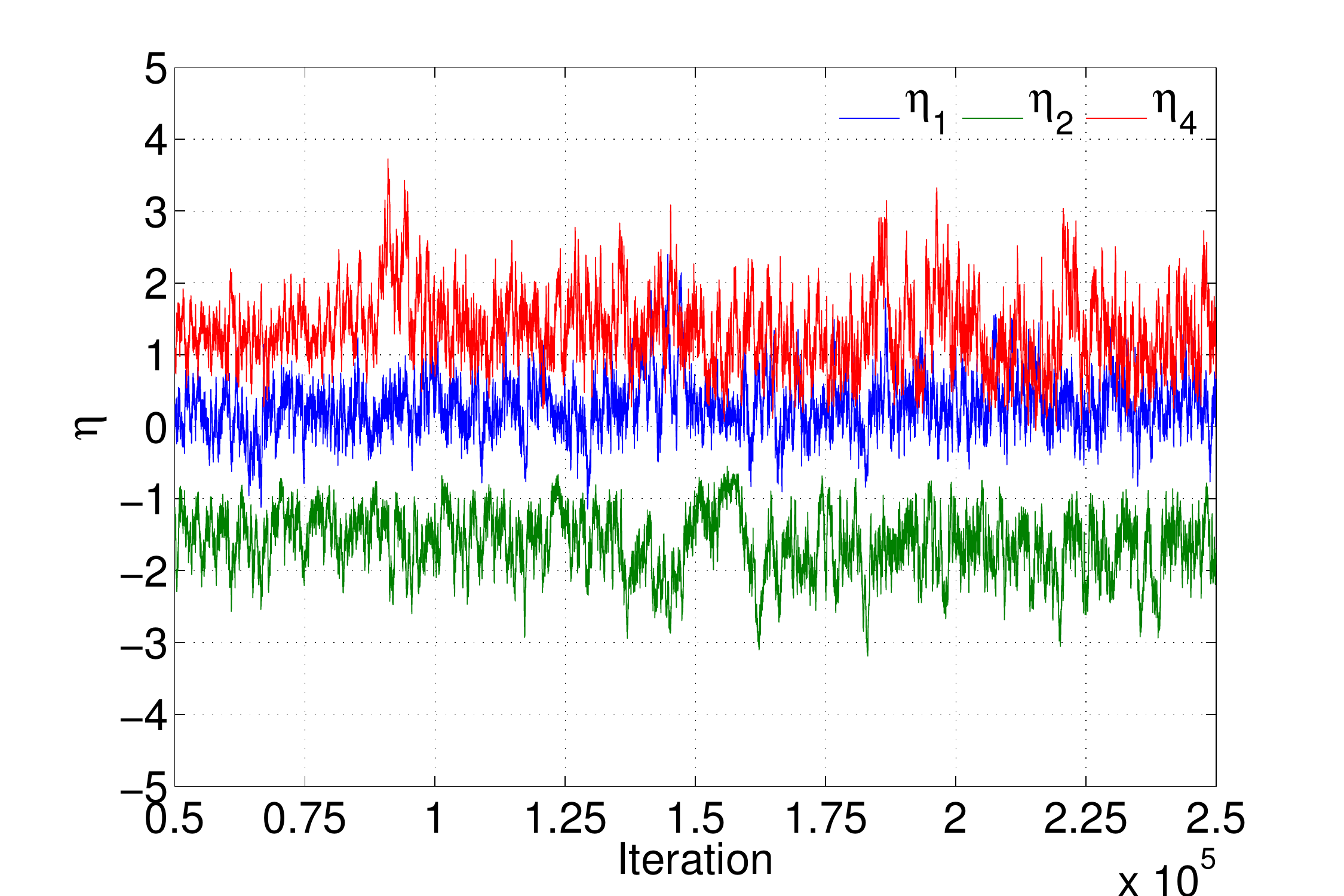}  &
\includegraphics[width=0.4\textwidth]{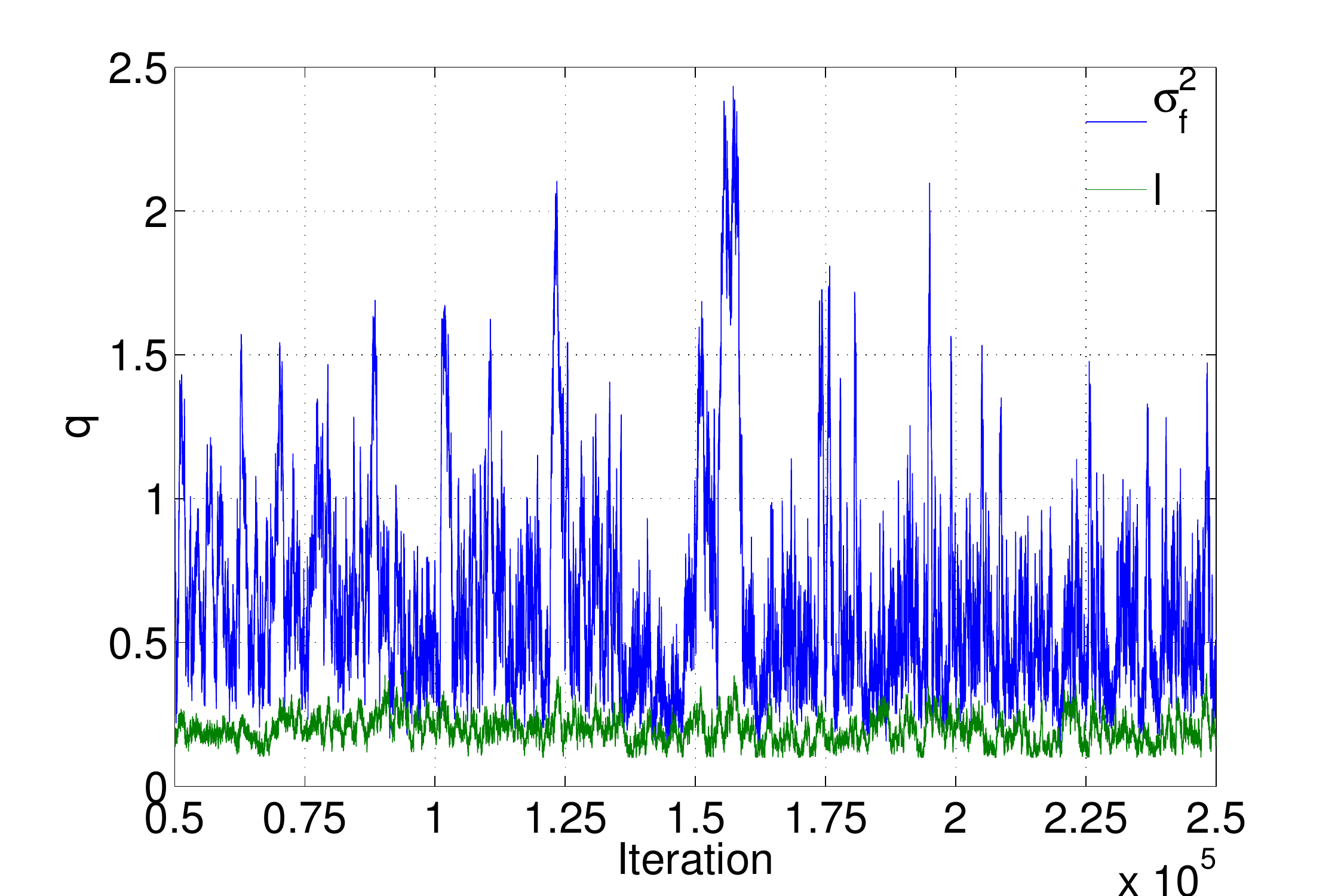}
\end{tabular}
\caption{Illustration of the chain generated by MCMC using PC surrogate and coordinate transformation: successive samples of 
(Left) few KL coordinates and (Right) hyper-parameters $\vec q$. Case of the inference of $m^{\rm sin}$.}
\label{fig:chain}
\end{figure}

Next, we repeat the previous inference problems but considering now the covariance hyper-parameters $l$ and $\sigma^2_f$ in addition to the $15$ KL modes and observation noise $\sigma_o^2$. For the sampling of the posterior, a total of $2.5\times 10^5$ MCMC steps was found also necessary to satisfactorily estimate the posterior statistics, same as for the case with pre-assigned parameters. 
The chains of all KL coordinates and hyper-parameters were observed to be well-mixed as illustrated in Figure~\ref{fig:chain}.

\begin{figure}[hbt]
\centering
\begin{tabular}{ccc}
\includegraphics[width=0.3\textwidth]{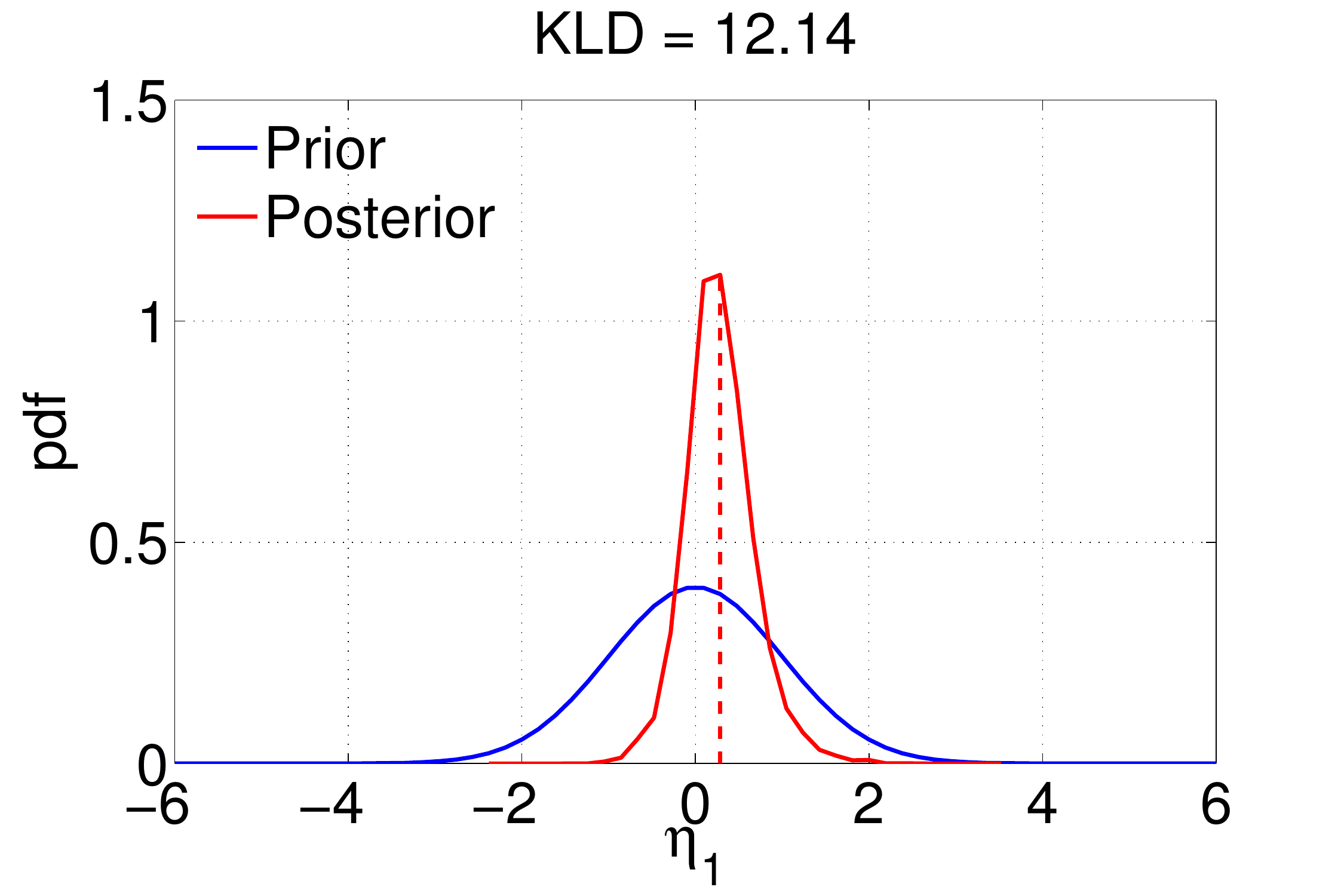}  &
\includegraphics[width=0.3\textwidth]{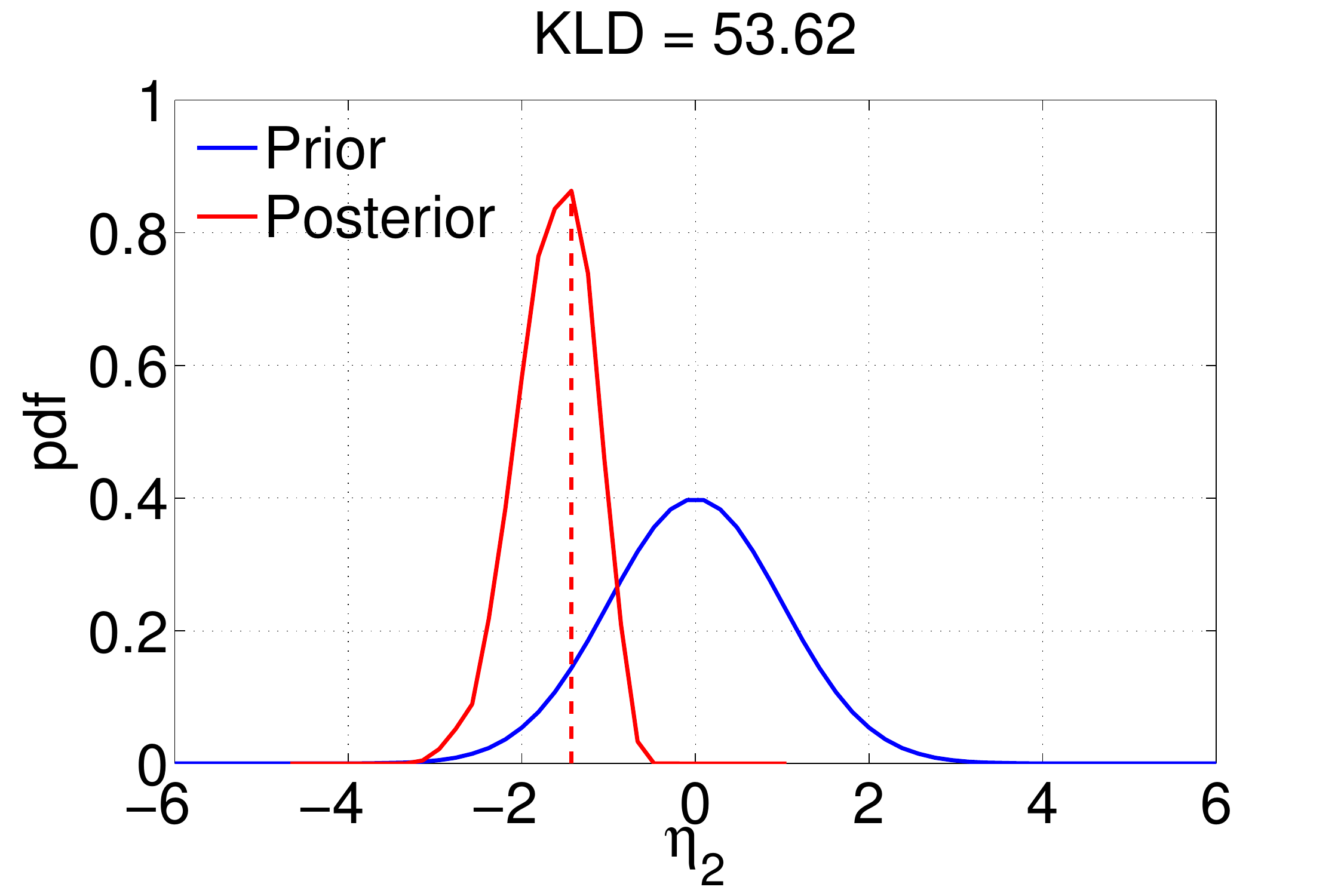}  &
\includegraphics[width=0.3\textwidth]{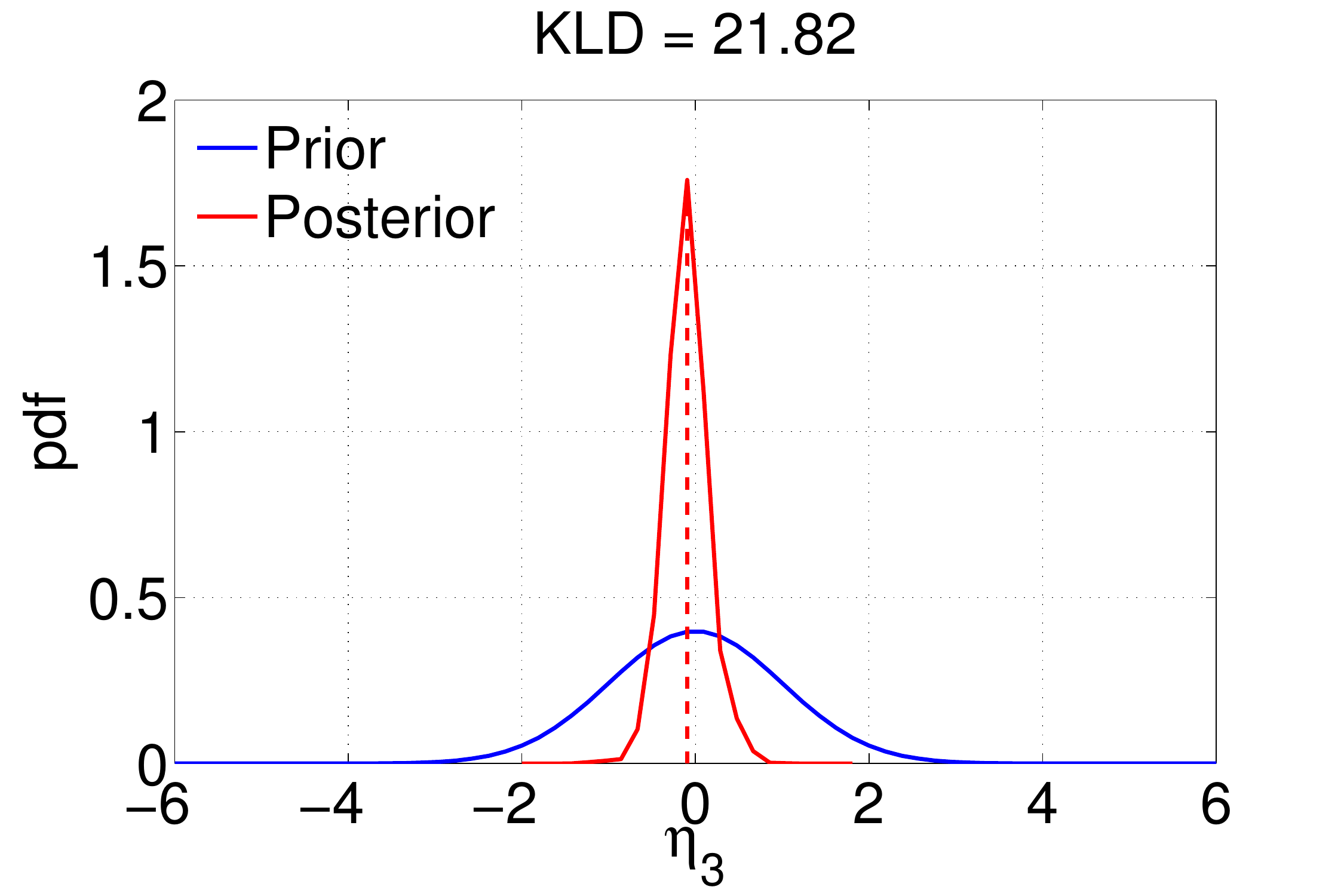}  \\
\includegraphics[width=0.3\textwidth]{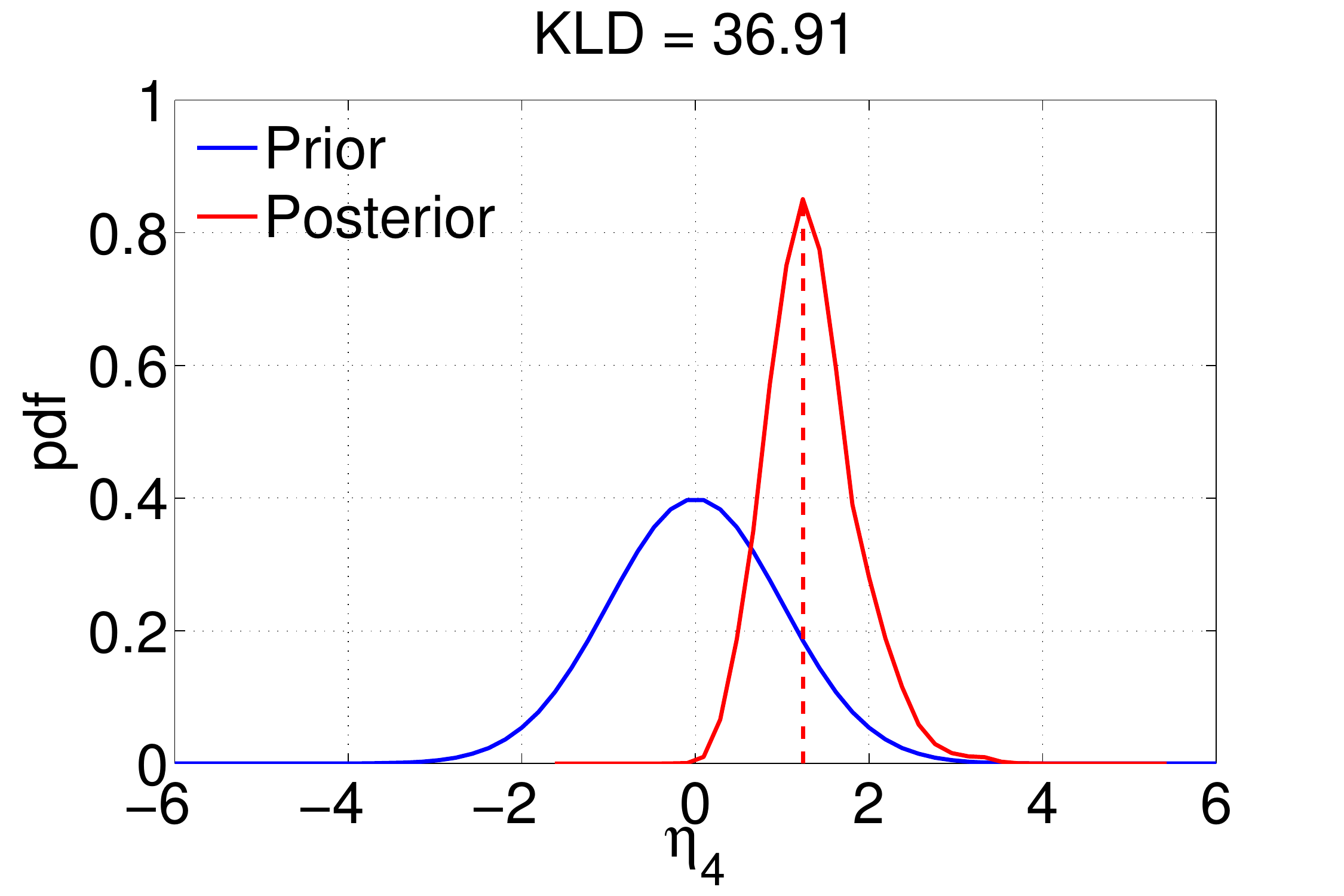}  &
\includegraphics[width=0.3\textwidth]{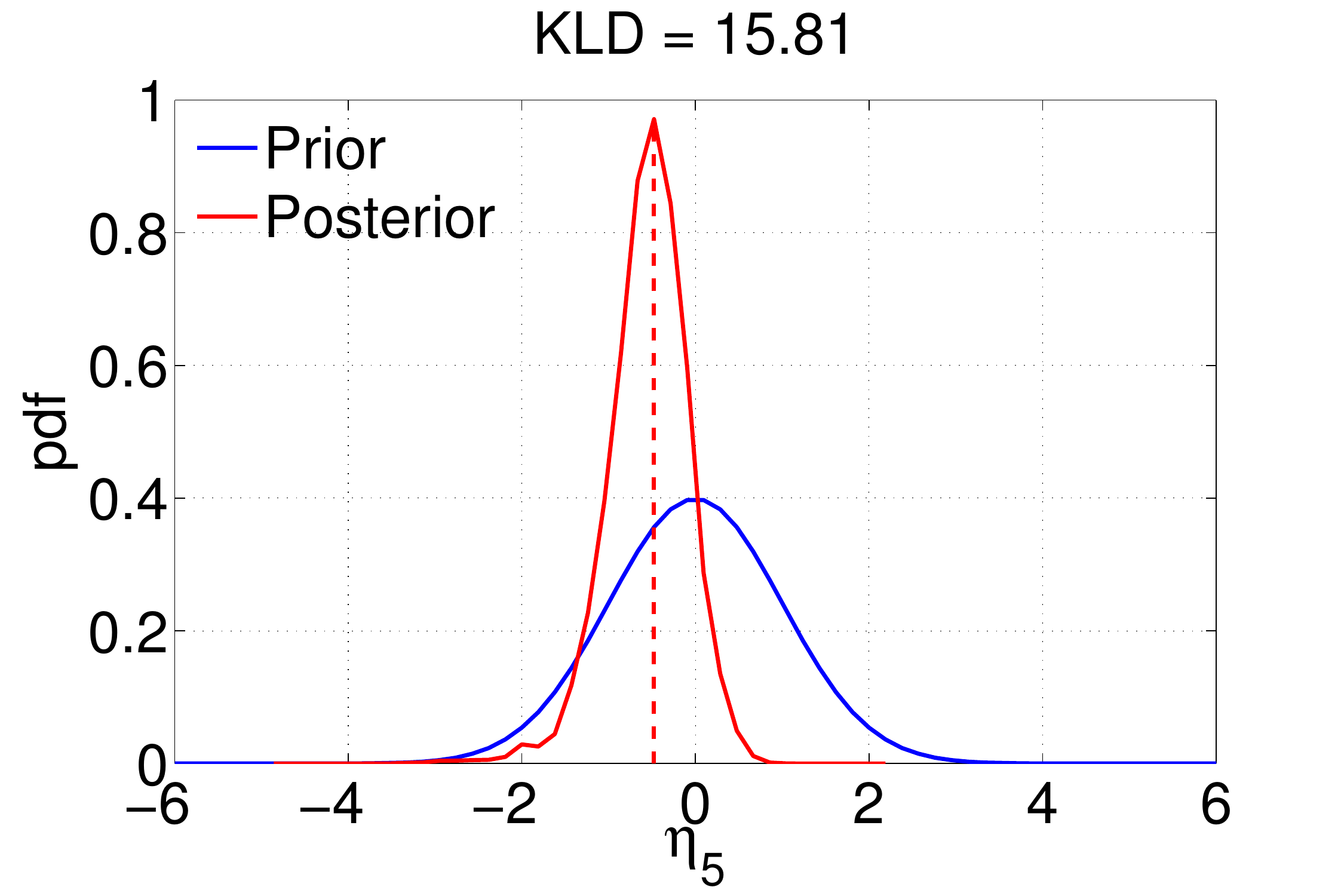}  &
\includegraphics[width=0.3\textwidth]{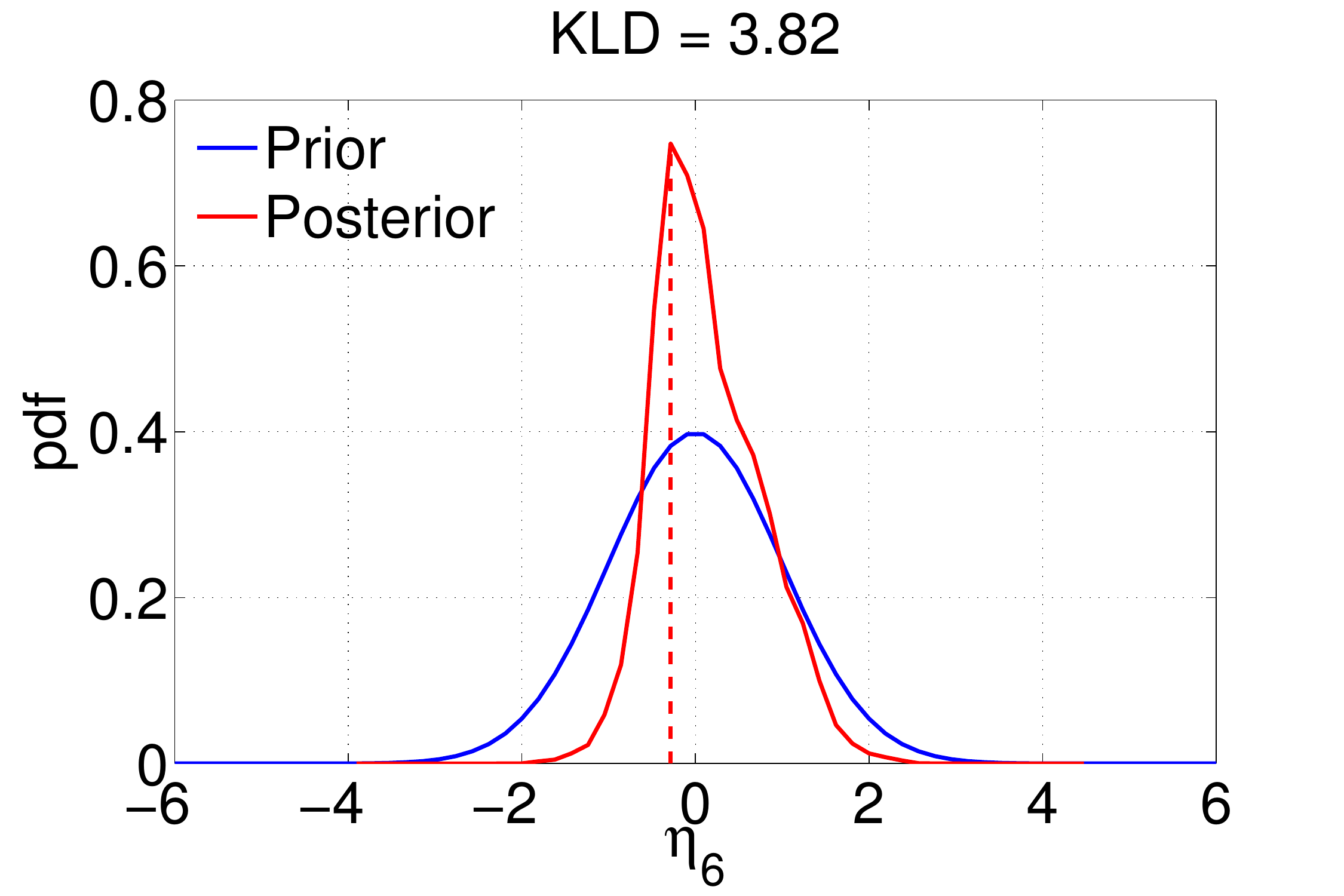}  \\
\includegraphics[width=0.3\textwidth]{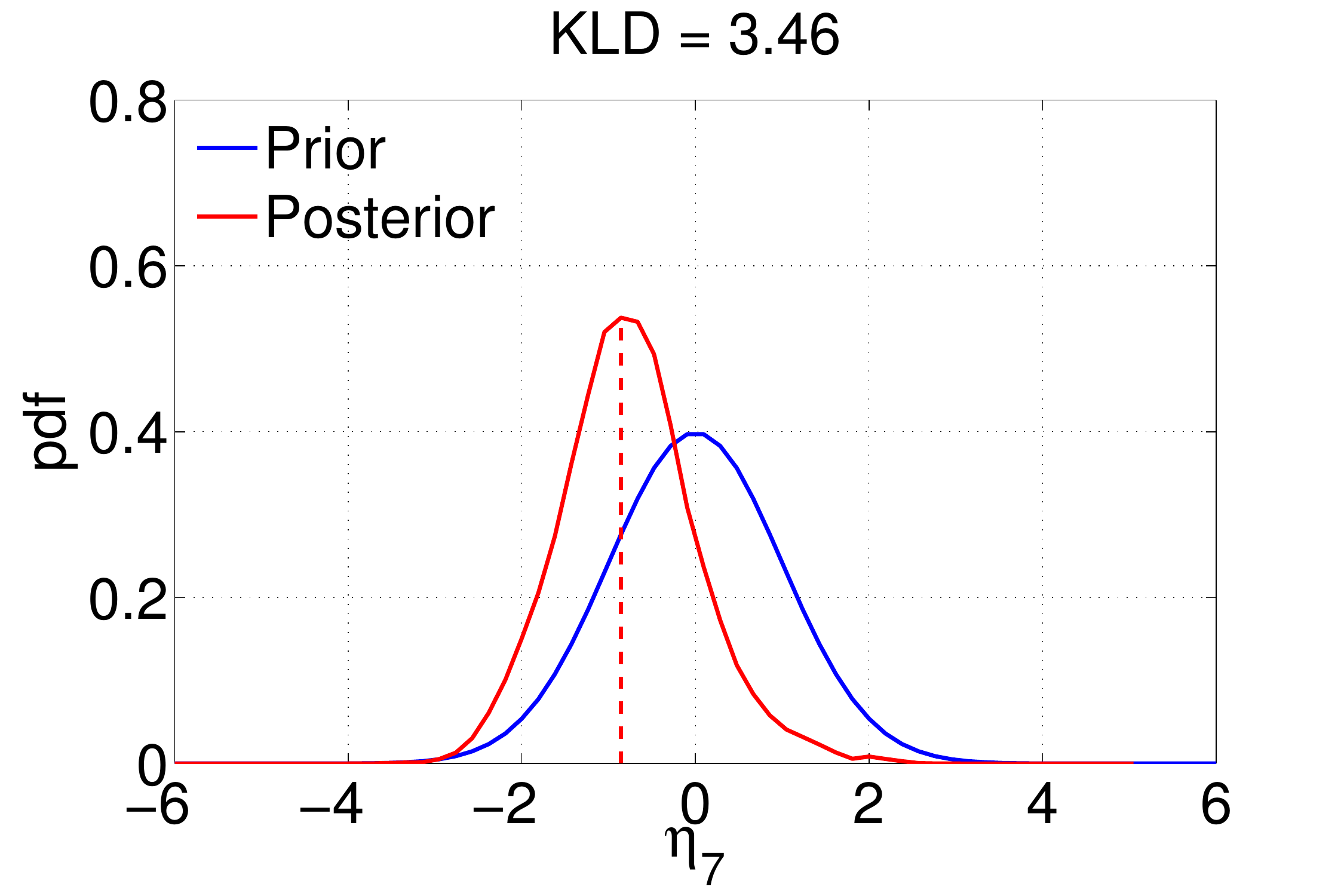}  &
\includegraphics[width=0.3\textwidth]{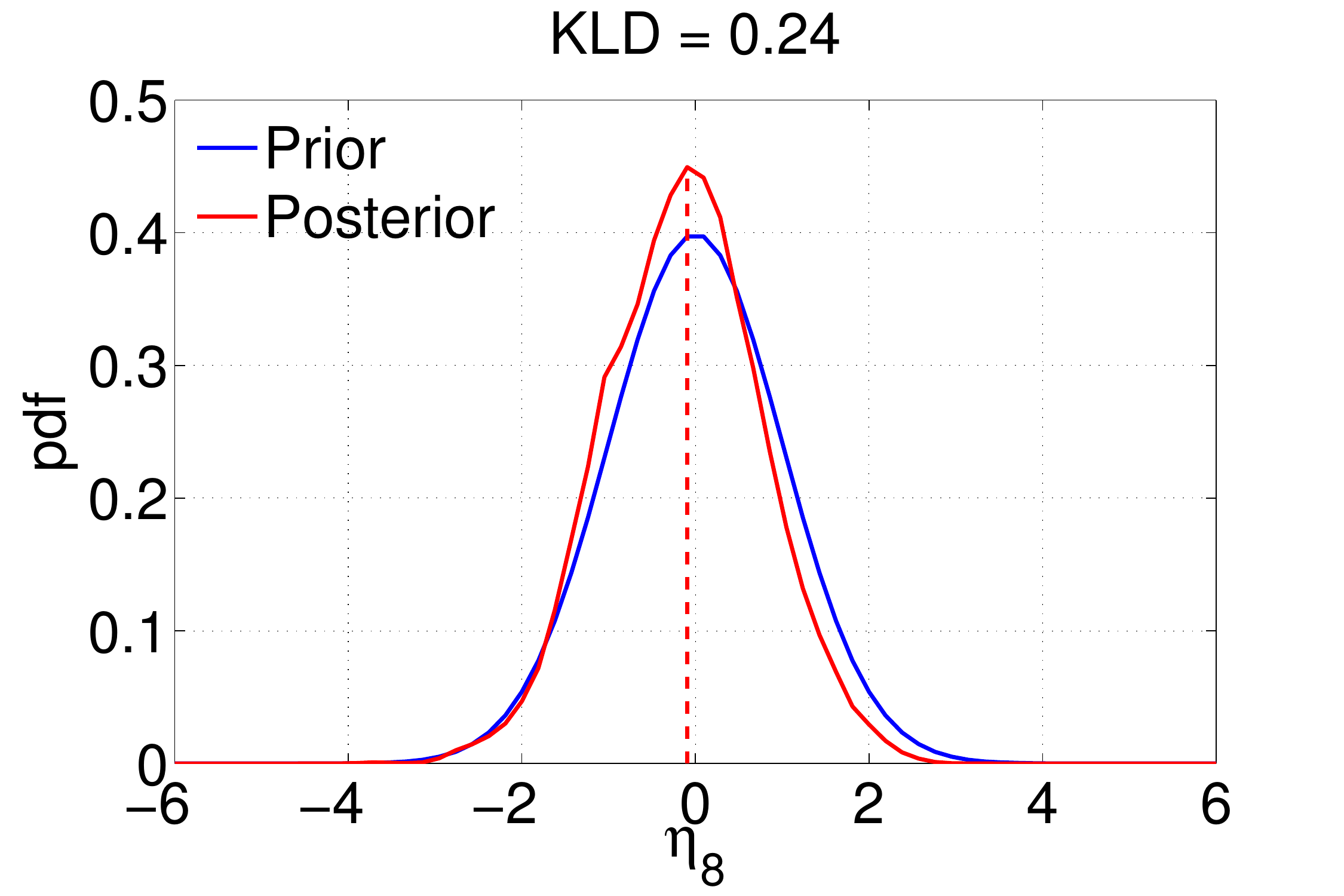}  &
\includegraphics[width=0.3\textwidth]{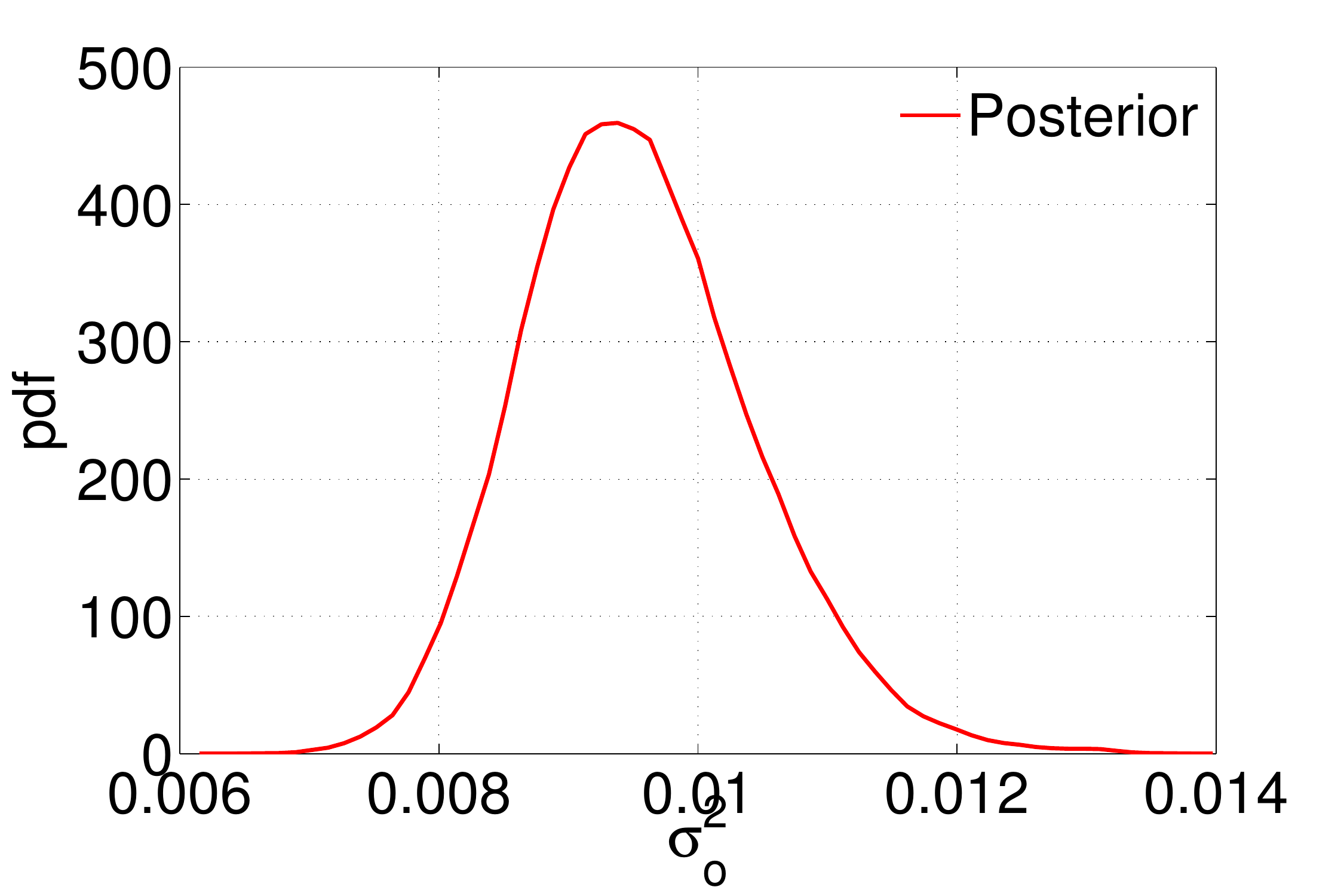}  
\end{tabular}
\caption{
Comparison of the priors and (marginal) posteriors of the first 8 KL coordinates $\eta_k$ and noise hyper-parameter $\sigma^2_o$ (posterior only) for the inference of $m^{\rm sin}$ with covariance hyper-parameters. 
The corresponding Kullback-Leibler Divergences (KLD) for the KL coordinates are also indicated on top of each plot.
}
\label{fig:pdfhyp}
\end{figure}

The marginal posteriors of the first 8 KL coordinates $\eta_k$ for $m^{\rm sin}$ are shown in Figure~\ref{fig:pdfhyp} together with their 
respective priors. 
The KLD values are also indicated on top of the plots. The results show a significant information gain for the first 7 KL coordinates, in contrast to only the first 4 KL coordinates when using pre-assigned parameters. In the same figure we show the marginal of the observation noise. The latter posterior has a MAP close to $\sigma^2_o = 0.01$, corresponding to the value used to generate the observations. Similar conclusions can be made for the cases of $m^{\rm step, ran}$ (results not shown for brevity).

\begin{figure}[hbt]
\centering
\begin{tabular}{cc}
	\includegraphics[width=0.4\textwidth]{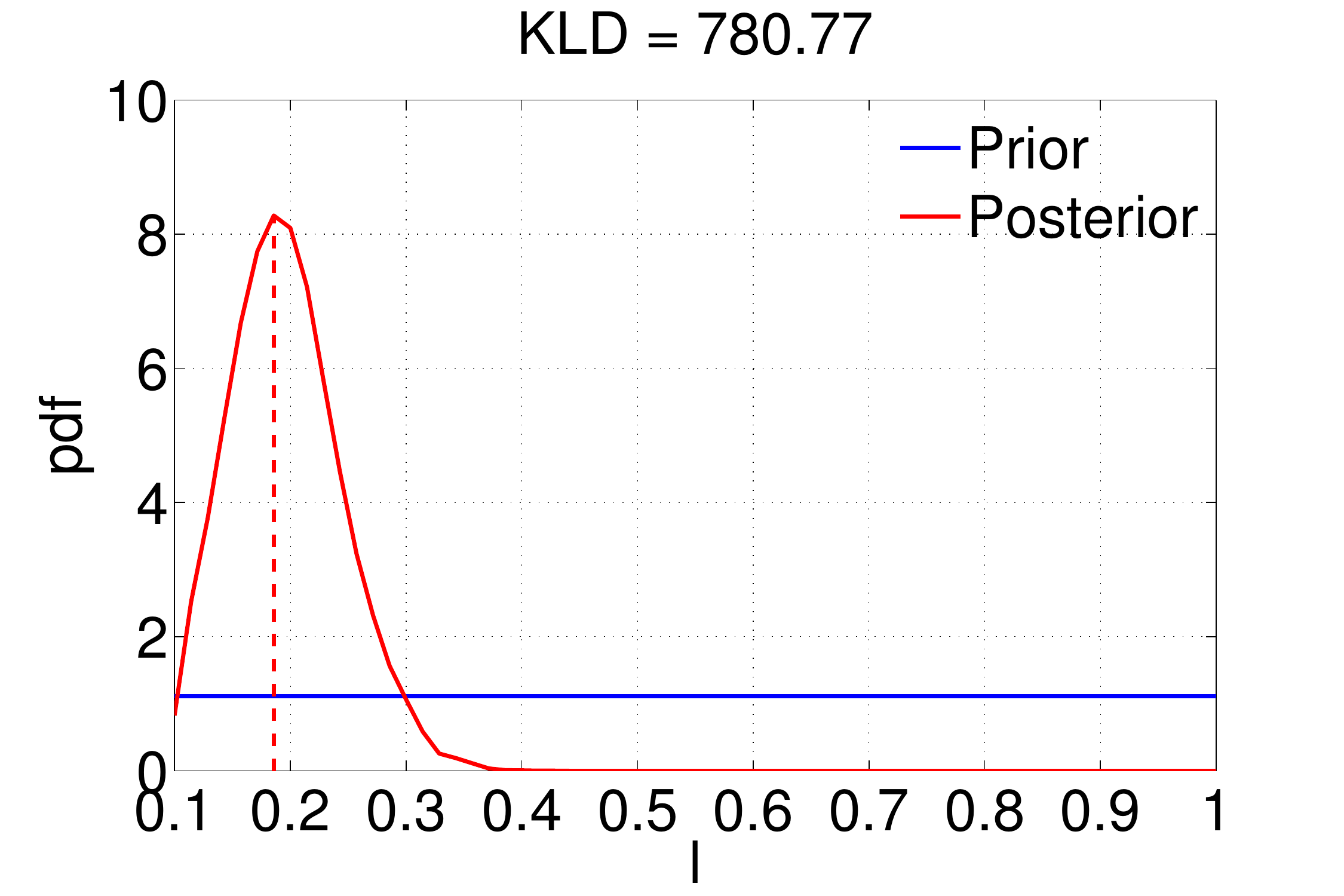}  &
\includegraphics[width=0.4\textwidth]{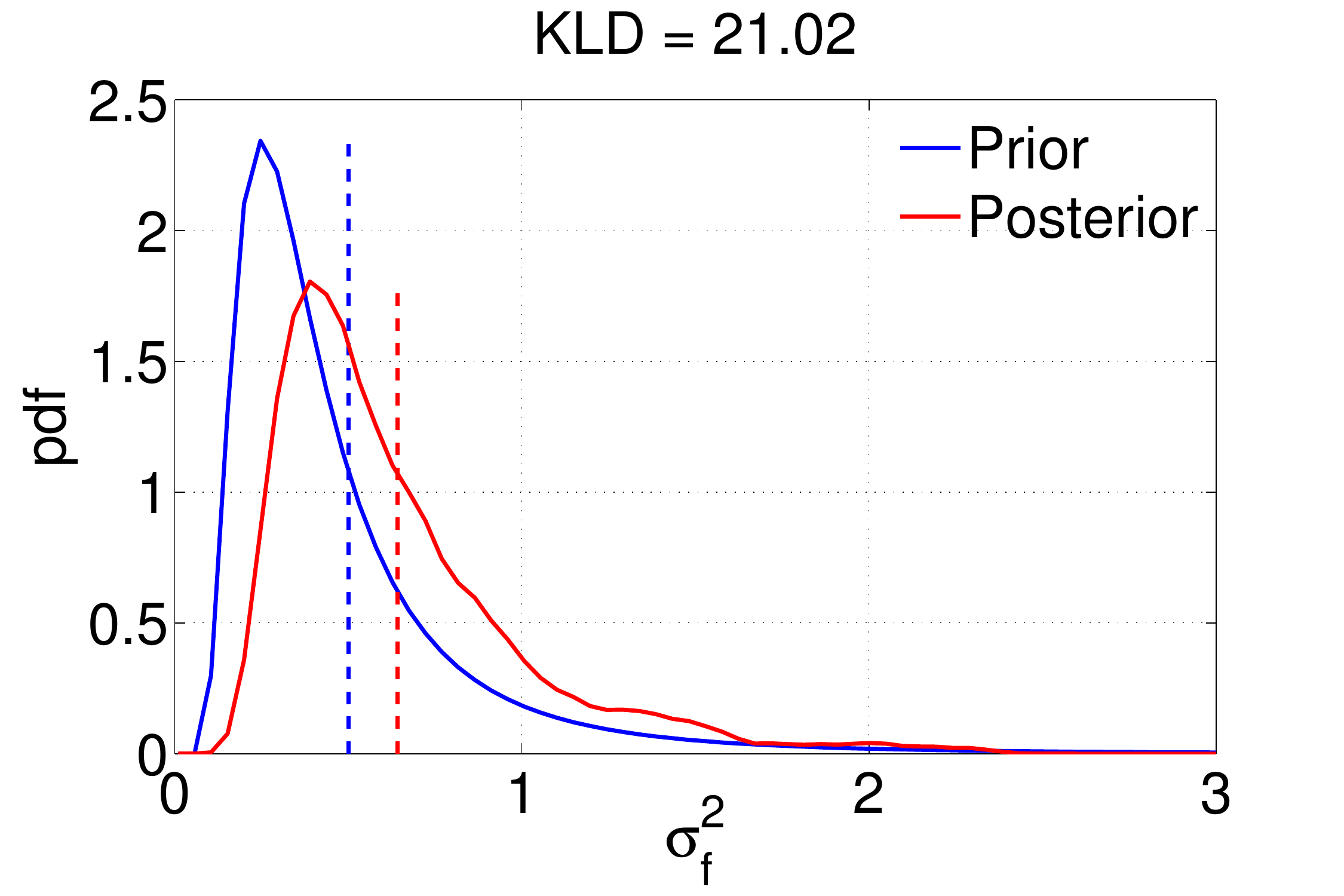}  \\
\end{tabular}
\caption{Priors and (marginal) posteriors of the covariance hyper-parameters $l$ (left plot) and $\sigma_f^2$ (right plot) for the inference of $m^{\rm sin}$.
Also indicated are the corresponding Kullback-Leibler Divergences (KLD).}
\label{fig:pdfhyper}
\end{figure}

The pdfs of the posterior of the hyper-parameters are shown in Figures~\ref{fig:pdfhyper} and compared with their priors for $m^{\rm sin}$.
The results show a significant difference between the prior and posterior of the covariance length scales $l$, with a MAP around $l=0.2$, while the posterior probability of $l>0.4$ is essentially zero. On the contrary, the posterior of the covariance variance $\sigma_f^2$ has a similar structure to that of its prior, 
with a shift of the expected (and MAP) value toward higher values. 

\subsubsection{Comparison with the inferences with and without hyper-parameters}

\begin{figure}[hbt]
\centering
\begin{tabular}{ccc}
\includegraphics[width=0.35\textwidth]{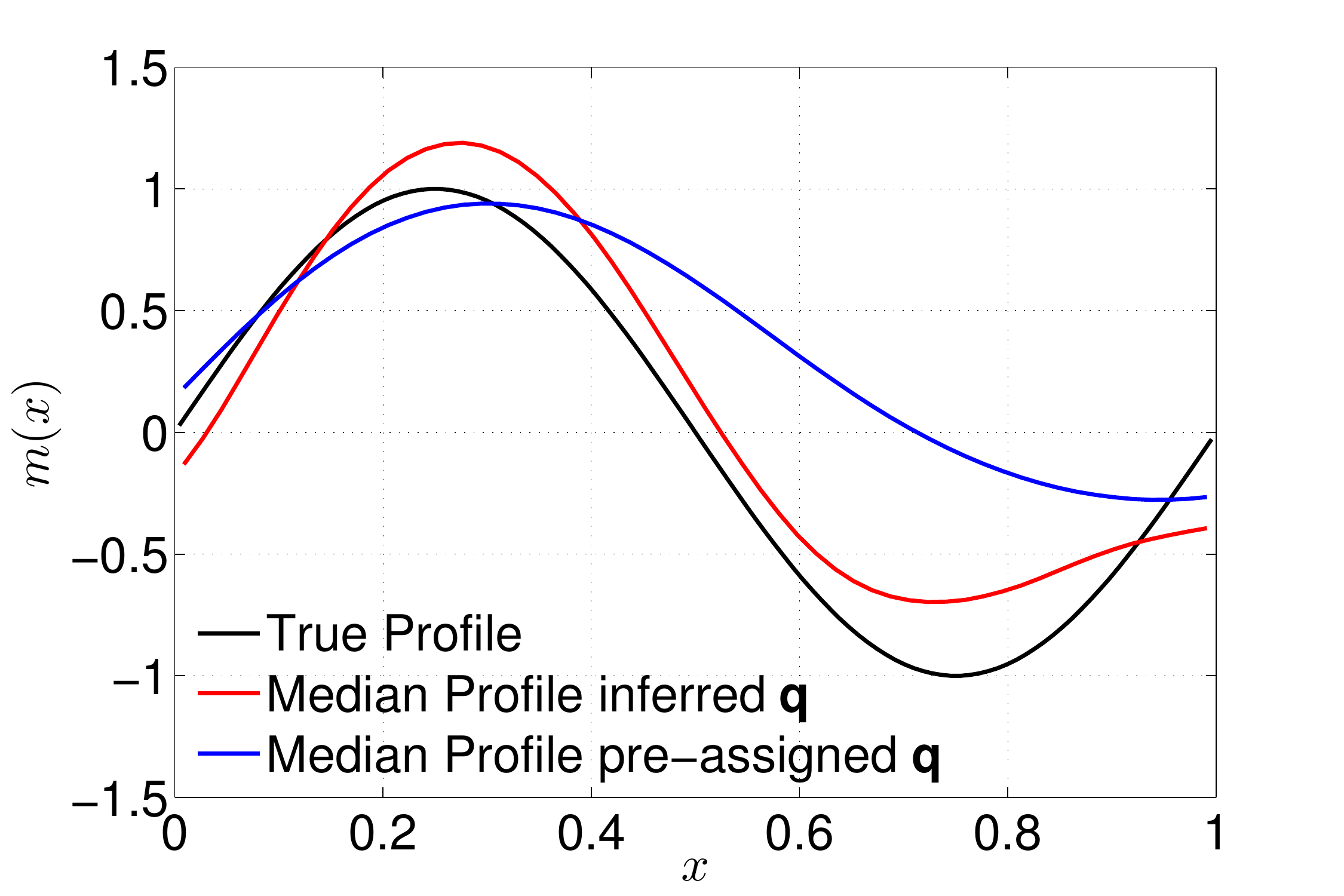} &
\includegraphics[width=0.35\textwidth]{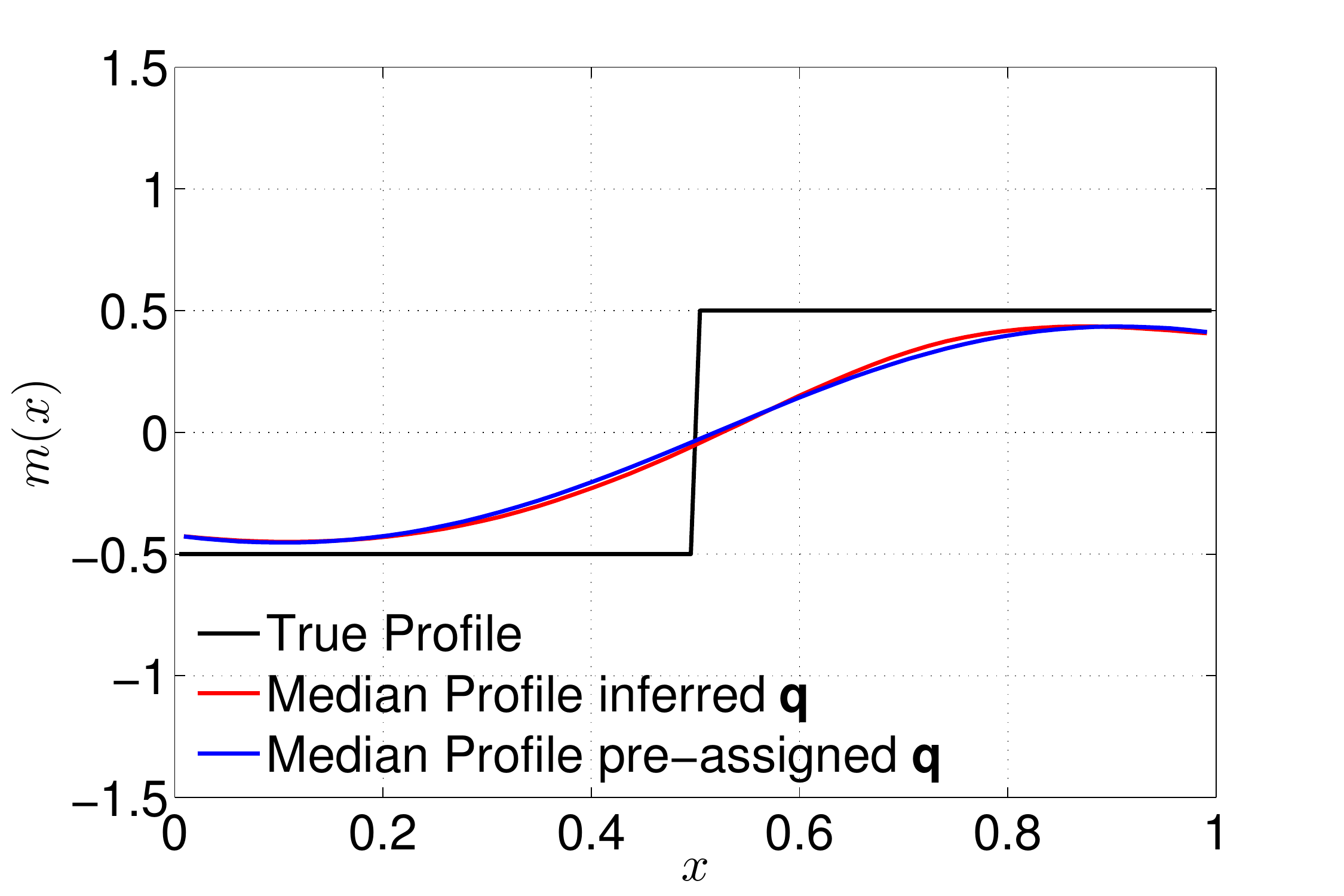} &
\includegraphics[width=0.35\textwidth]{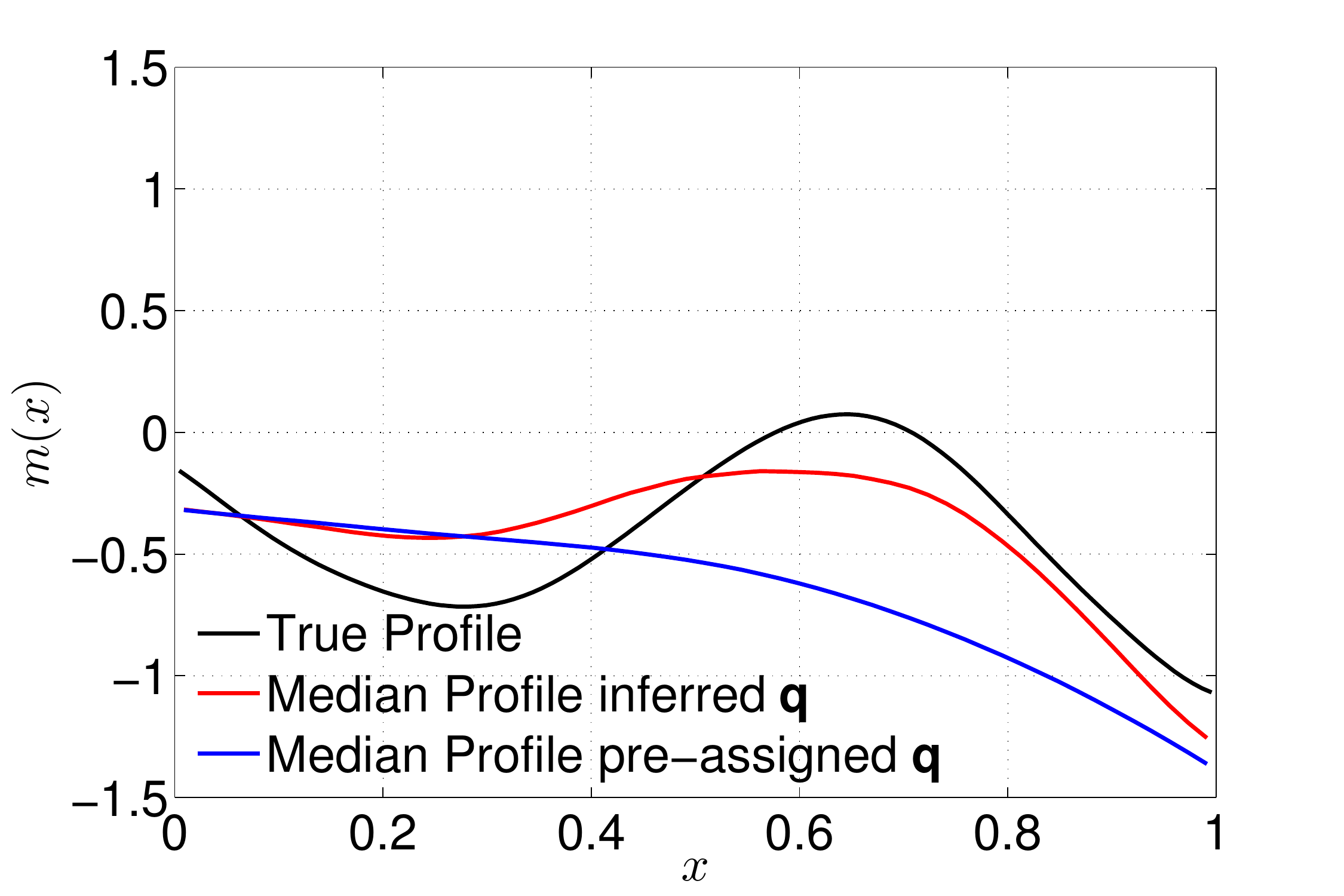} 
\end{tabular}
\caption{Comparison of the true log-diffusivity profiles with corresponding posterior medians for the inference with covariance hyper-parameters and preassigned covariance. Cases of $m^{\rm sin}$, $m^{\rm step}$ and $m^{\rm ran}$ from left to right.}
\label{fig:prfcom}
\end{figure}

To better appreciate the improvement resulting from the introduction of the covariance hyper-parameters, we first provide a comparison of the inferred median profiles, obtained by inferring covariance hyper-parameters or by using pre-assigned values.  The median profiles for all three cases are plotted in Figure~\ref{fig:prfcom}, which also depicts the true profiles. It is seen that introducing the covariance hyper-parameter significantly reduces the distance between the  median and true profiles in the smooth cases ($m^{\rm sin}$ and $m^{\rm ran}$), while having no significant impact on the inference of the piecewise constant profile $m^{\rm step}$. 
This behavior can be explained by the family of Gaussian processes considered, which is not well-suited for the inference of $m^{\rm step}$, and so the introduction of hyper-parameters does not help improving the inference.

\begin{figure}[hbt]
\centering
\begin{tabular}{ccc}
\includegraphics[width=0.35\textwidth]{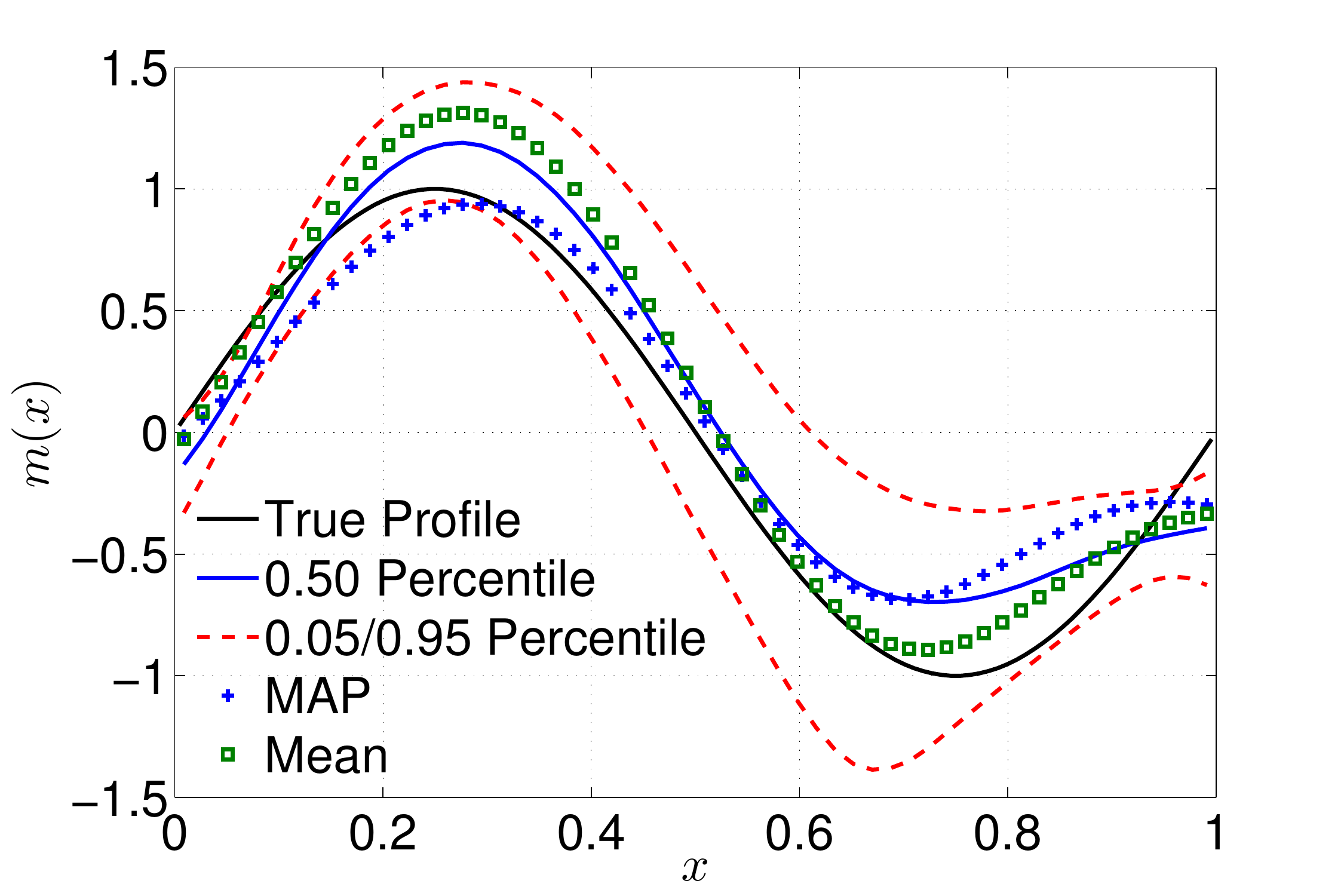}  &
\includegraphics[width=0.35\textwidth]{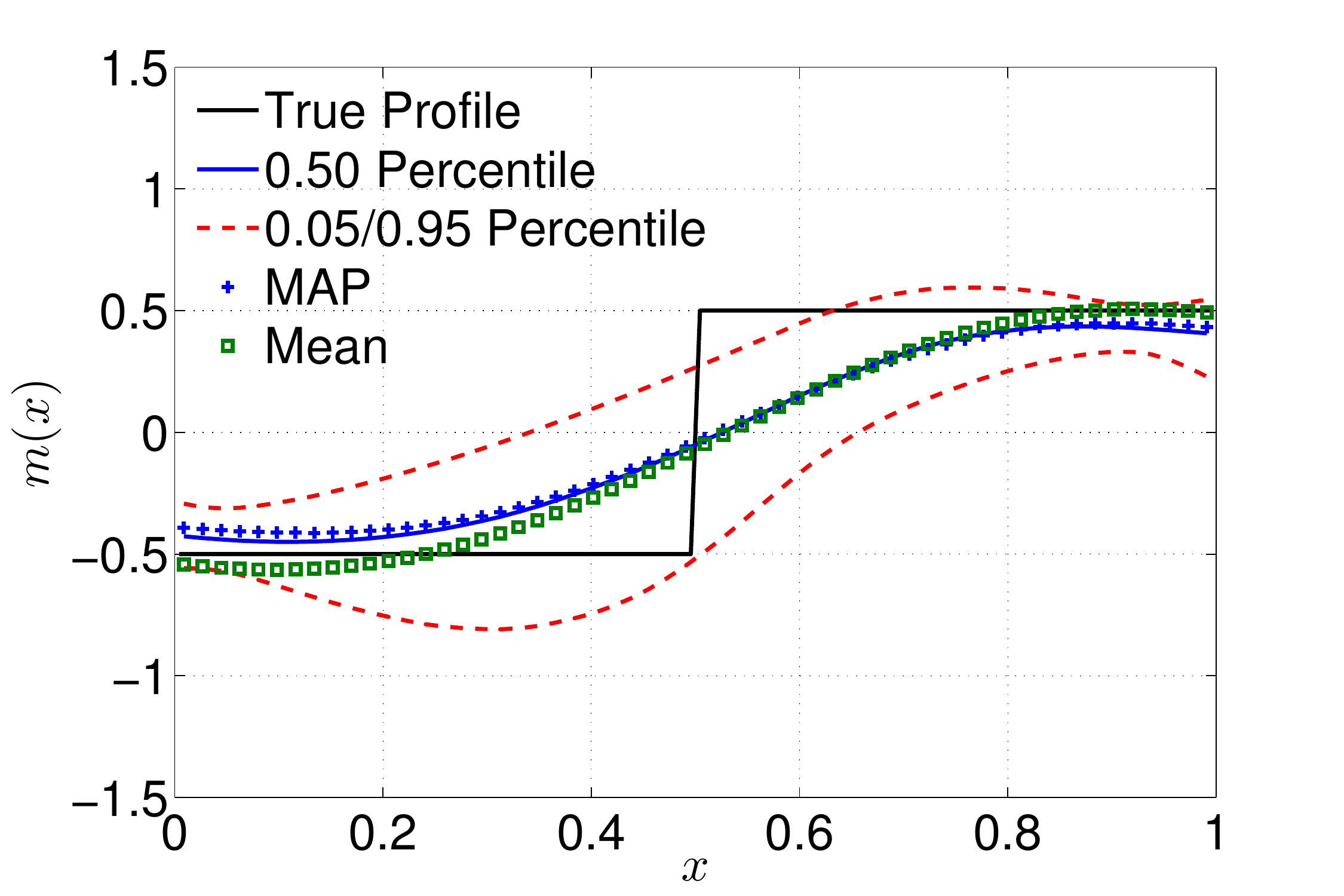} &
\includegraphics[width=0.35\textwidth]{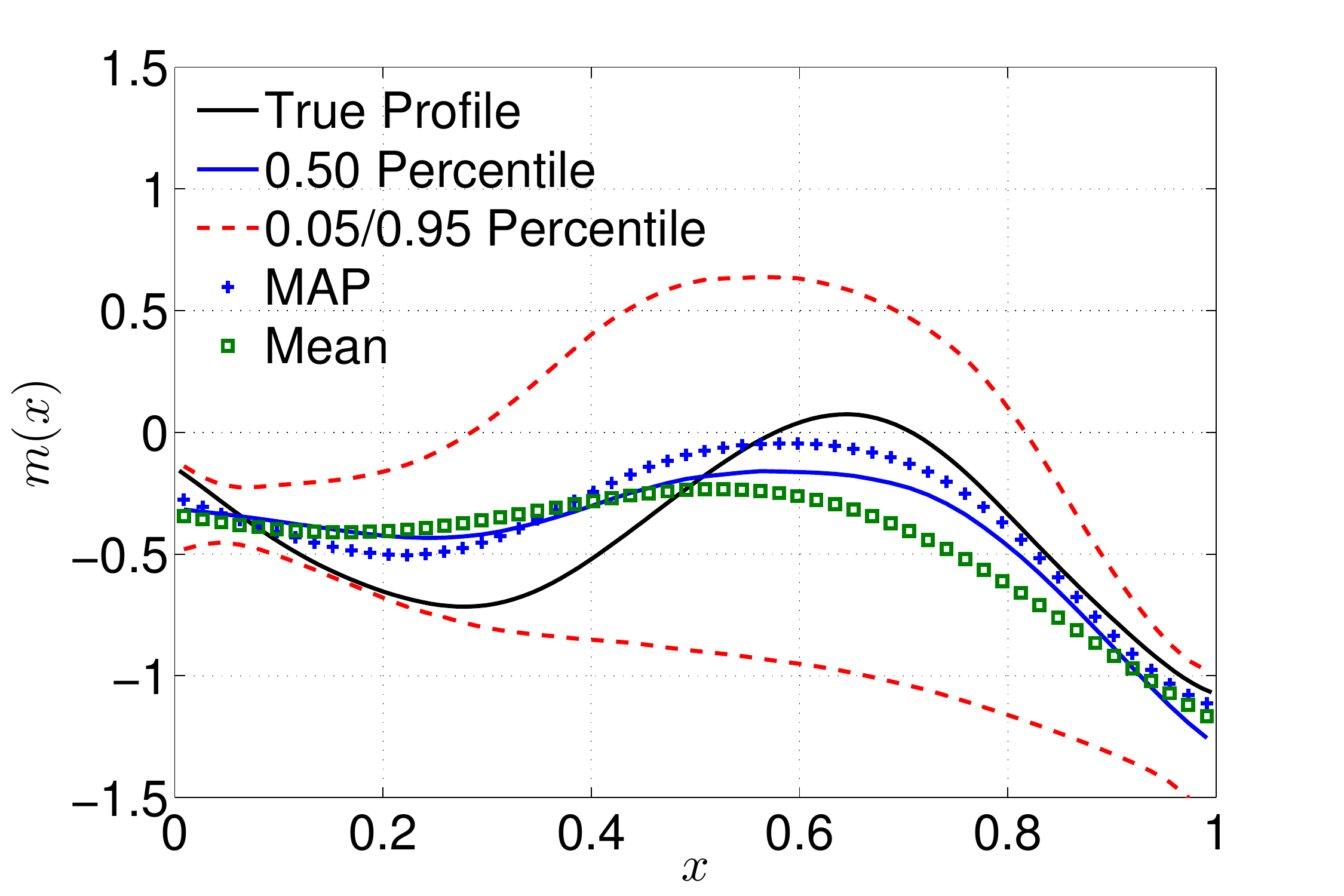} 
\end{tabular}
\caption{Comparison of the posteriors profiles with the true ones, for the cases of $m^{\rm sin}$, $m^{\rm step}$ and $m^{\rm ran}$ (from left to right).
The inferences use  covariance function with hyper-parameters. Shown in each plots are the median, mean, MAP, $5\%$ and $95\%$ quantiles of the posterior and true profiles.}
\label{fig:prfhyp}
\end{figure}

Second, the median, mean, MAP, $5\%$ and $95\%$ quantiles of the inferred log-diffusivity profiles are 
plotted in Figure~\ref{fig:prfhyp} and compared with the respective true profiles $m^{\rm sin, step, ran}$. 
These plots should be contrasted with the results shown in Figure~\ref{fig:prfnhyp}, obtained with pre-assigned covariance. 
Consistent with the previous observations on the medians, we observe that in the case of the discontinuous profile, $m^{\rm step}$, the 
inference of the hyper-parameters only affects slightly the 5\% and 95\% quantiles. 
On the contrary, for the smooth profiles $m^{\rm sin, ran}$ the 5\% and 95\% quantiles bounds now contain the true profiles for nearly every $x$. 
This significant improvement is due partly to the better agreement between the true and median profiles, but also to a generally higher variability in the posterior when considering the hyper-parameters in the inference process. 
In other words, the inference of the covariance hyper-parameters appears to yield a more flexible approach than when using a fixed covariance assumption.

\subsubsection{Effects of measurement noise and number of observations}
\begin{figure}[hbt]
\centering
\begin{tabular}{ccc}
\includegraphics[width=0.35\textwidth]{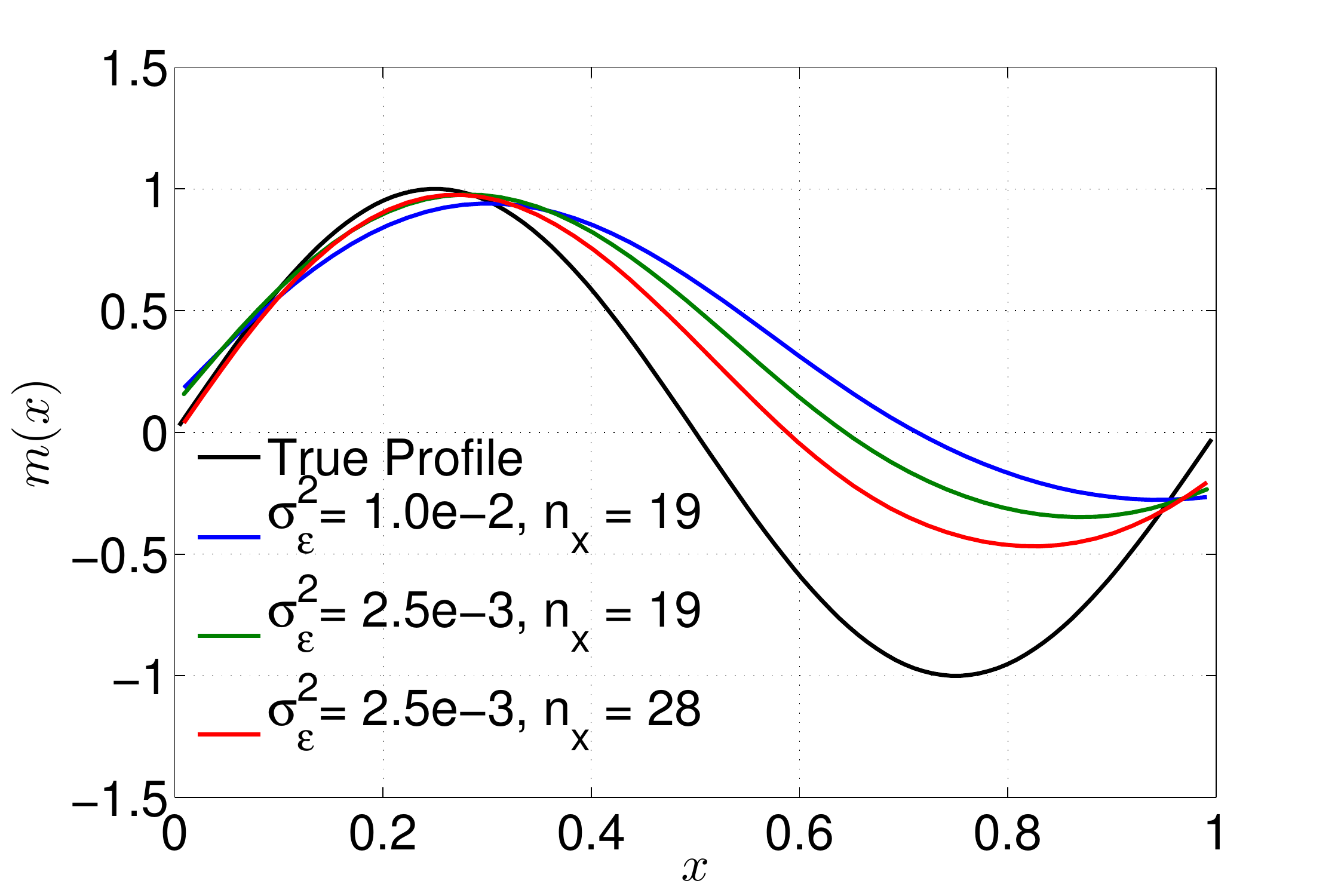}  &
\includegraphics[width=0.35\textwidth]{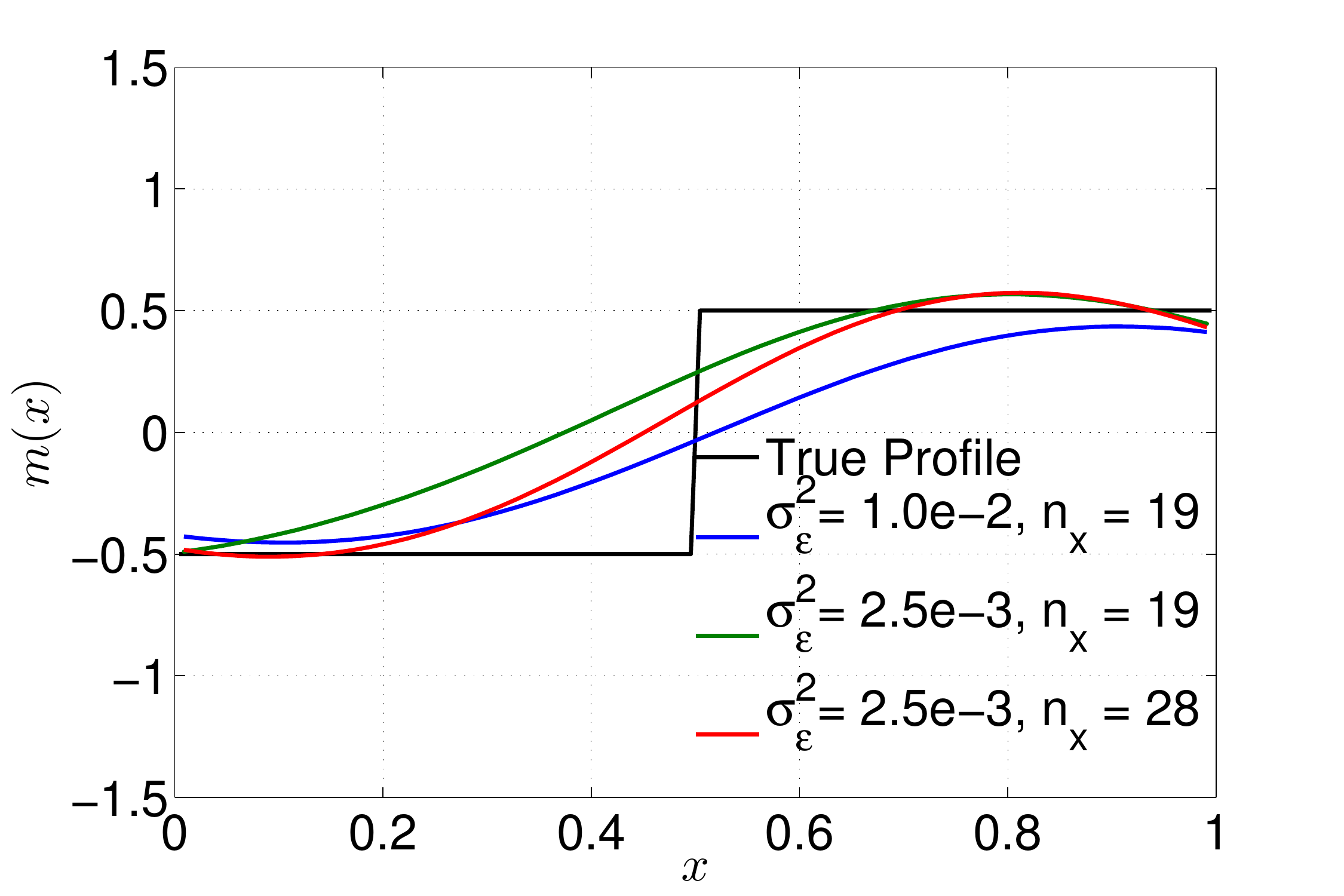}  &
\includegraphics[width=0.35\textwidth]{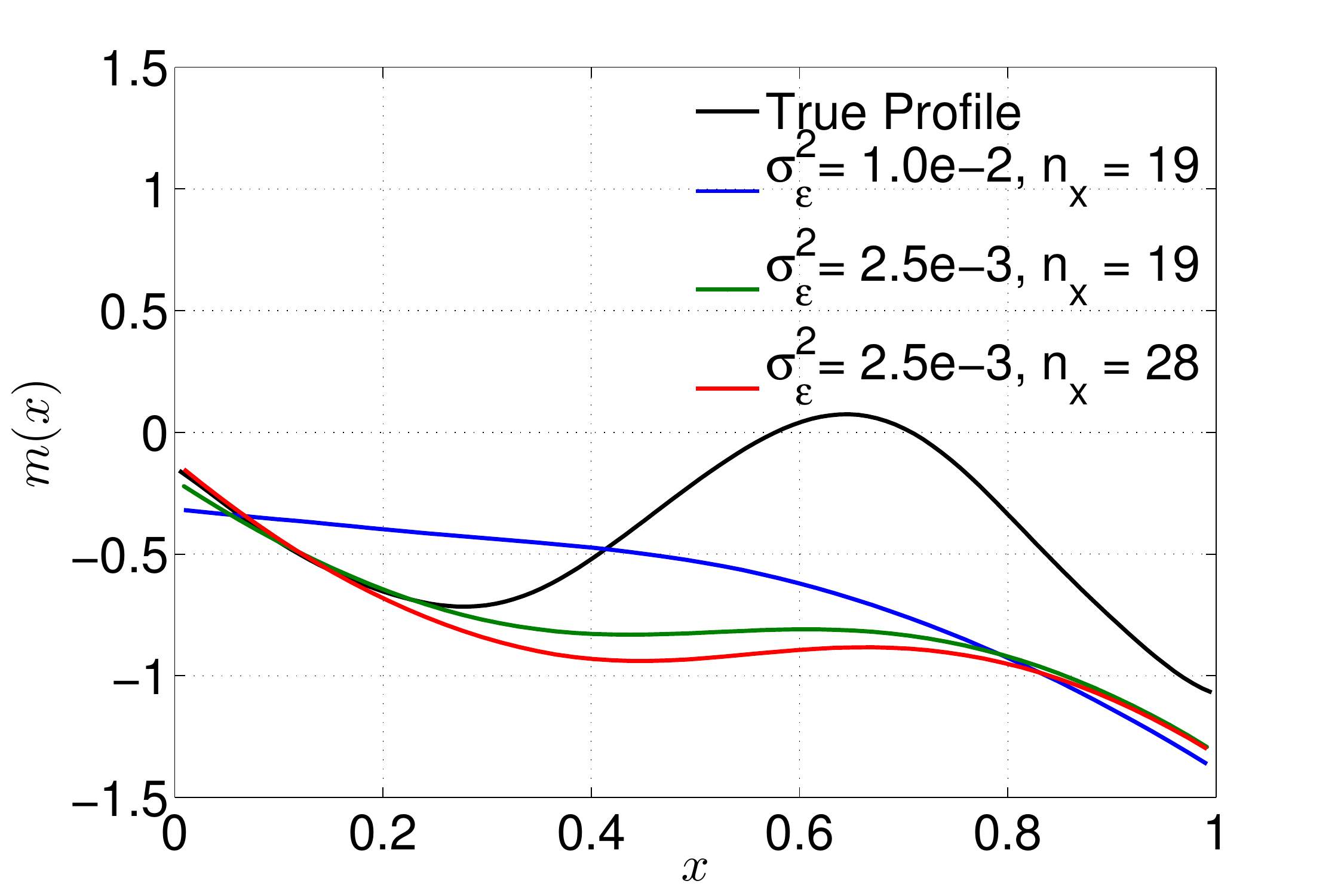}   \\
\includegraphics[width=0.35\textwidth]{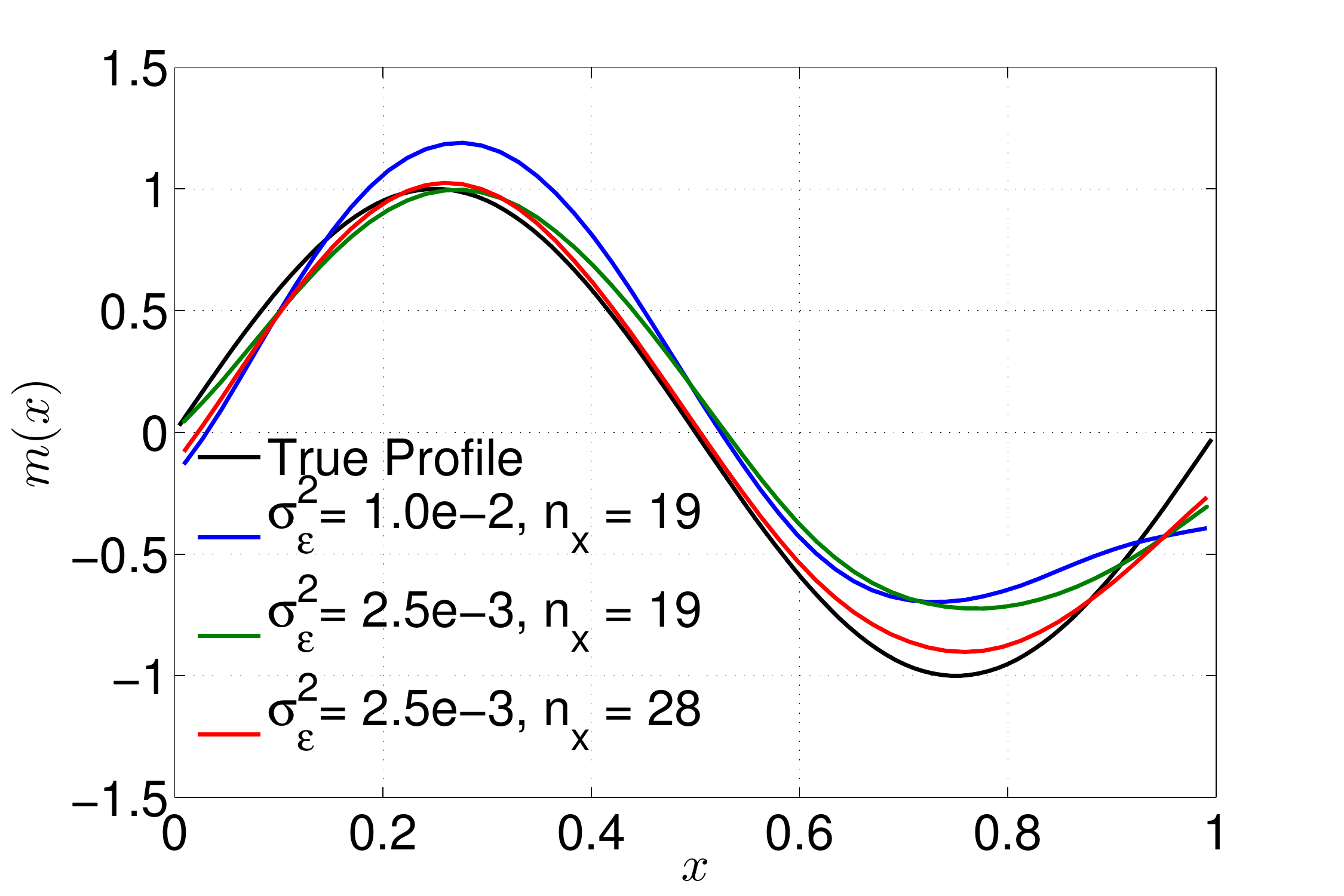}  &
\includegraphics[width=0.35\textwidth]{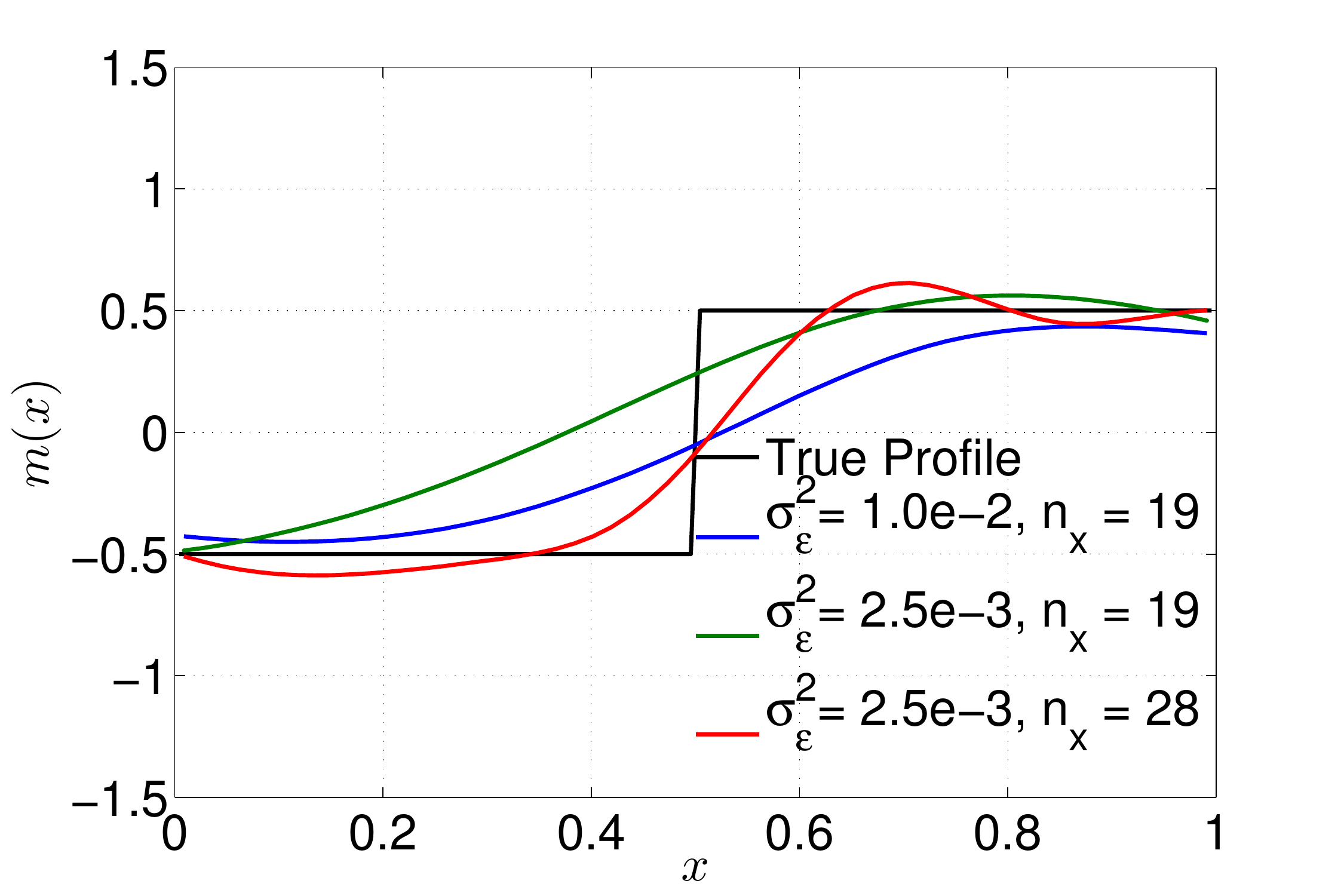}  &
\includegraphics[width=0.35\textwidth]{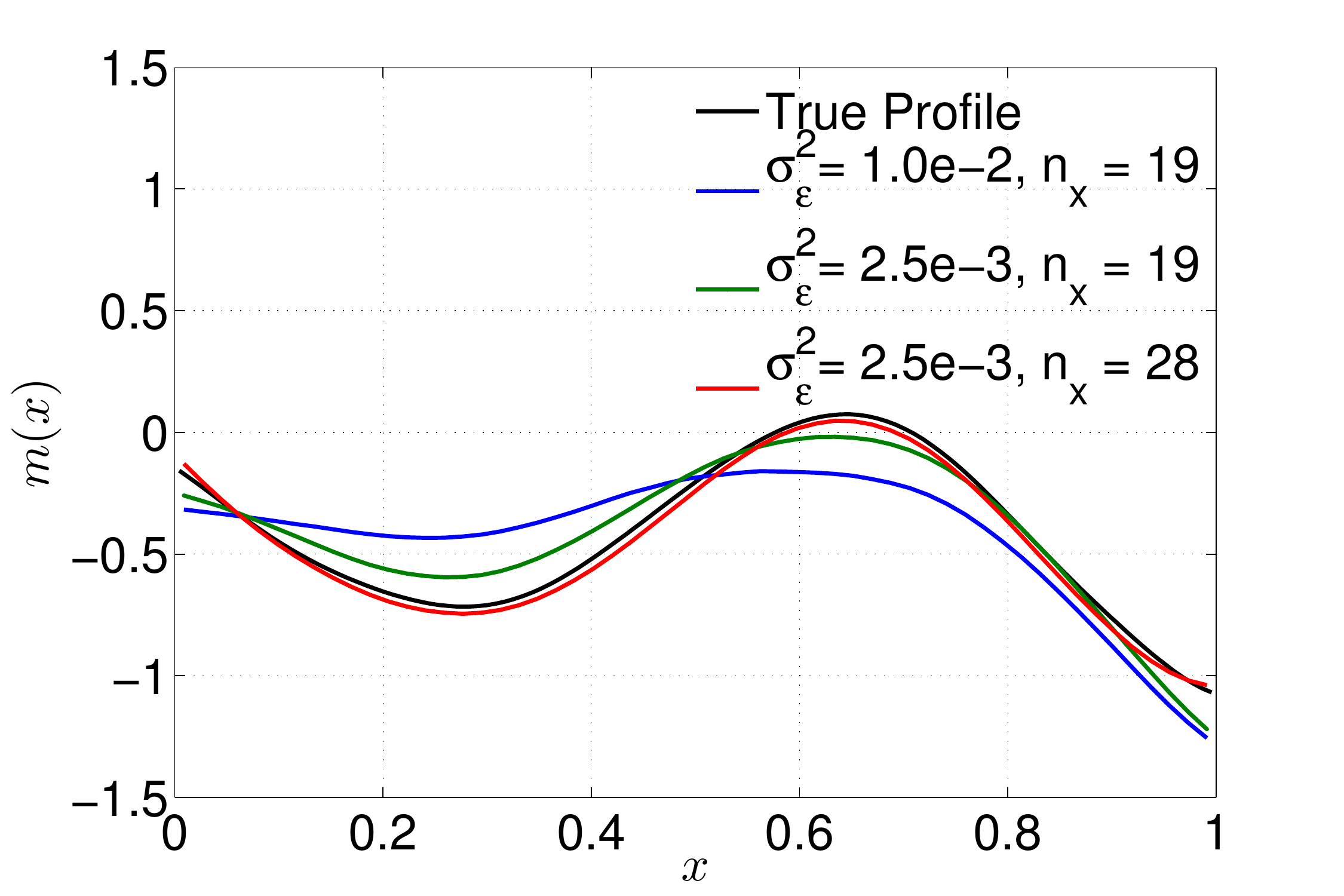}  
\end{tabular}
\caption{Effect of observations number and noise. Shown are the medians of the inferred profiles for the three test cases $m^{\rm sin, step, ran}$ (from left to right), and inferences for a pre-assigned covariance function (top row) or with hyper-parameters (bottom row).}
\label{fig:noise_effect}
\end{figure}

To investigate the impact of the observations on the inference processes, with or without covariance hyper-parameters, 
we repeat the previous inference problems for different noise level $\sigma_\epsilon^2$ in the observations and different number of spatial locations $n_x$.
The results are reported in Figure~\ref{fig:noise_effect} in terms of median profiles, for the three test profiles $m^{\rm sin, step, ran}$ (from left to right) and the inference without (top row) and with covariance hyper-parameters (bottom row). 
As expected, the plots indicate an improvement of the inferred (median) profiles when the noise level is lowered, and when the number of observation increases. The improvements are more significant in the cases of the smooth profiles ($m^{\rm sin, ran}$) than for the discontinuous one ($m^{\rm step}$), a result consistent with the previous observations.
In addition, for the smooth cases, the improvements carried by the introduction of the covariance hyper-parameters in the inference problem is seen to not only yield median profiles closer to the true ones, but also to significantly accelerate the convergence to the truth. The improvement of the convergence rate would require additional numerical experiments to be precisely measured, but it can already be safely asserted that more information is gained from the observations when considering the covariance hyper-parameters in the inference.

\subsubsection{Convergence with the PC surrogate order}
\begin{figure}[hbt]
\centering
\begin{tabular}{ccc}
\includegraphics[width=0.35\textwidth]{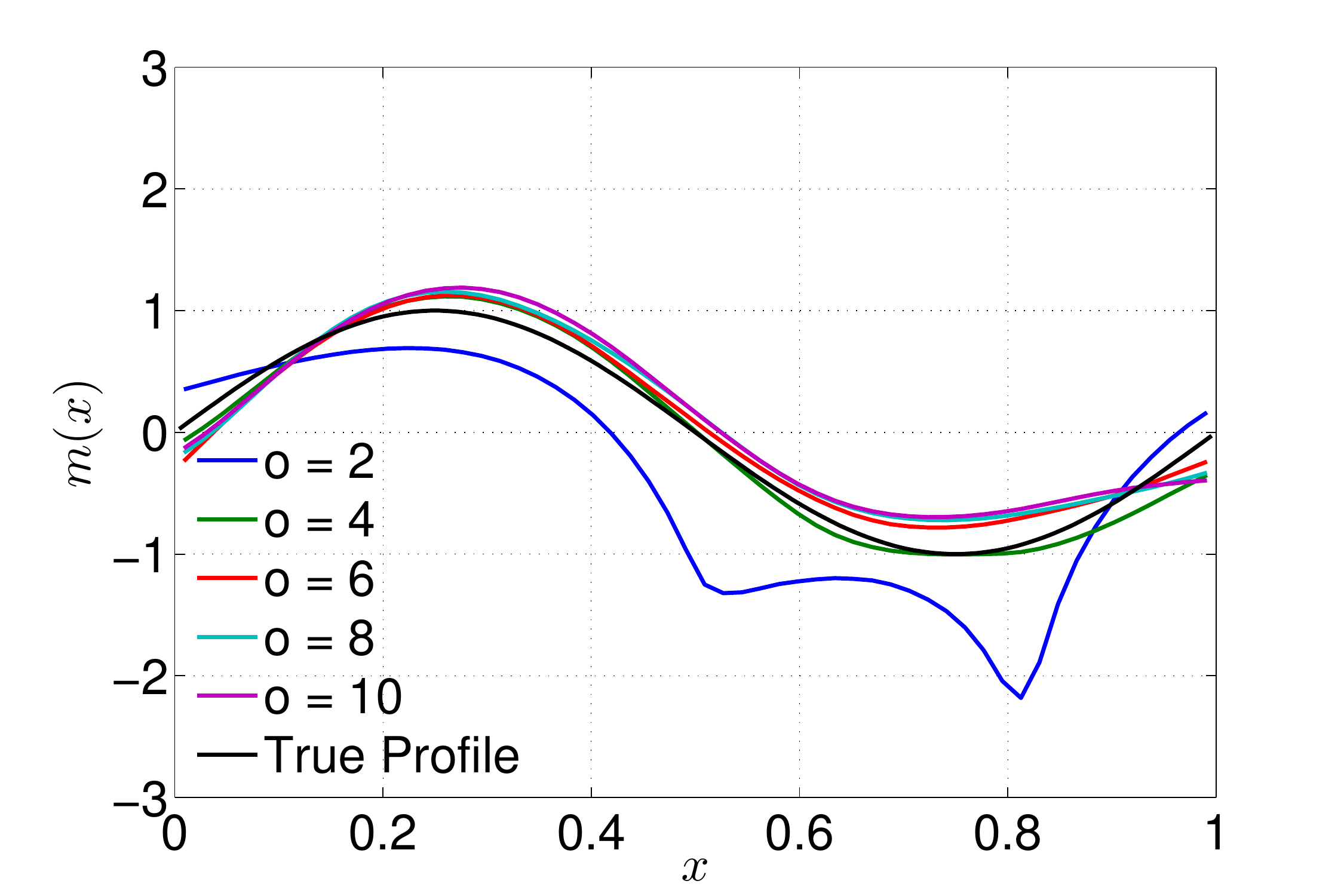}  &
\includegraphics[width=0.35\textwidth]{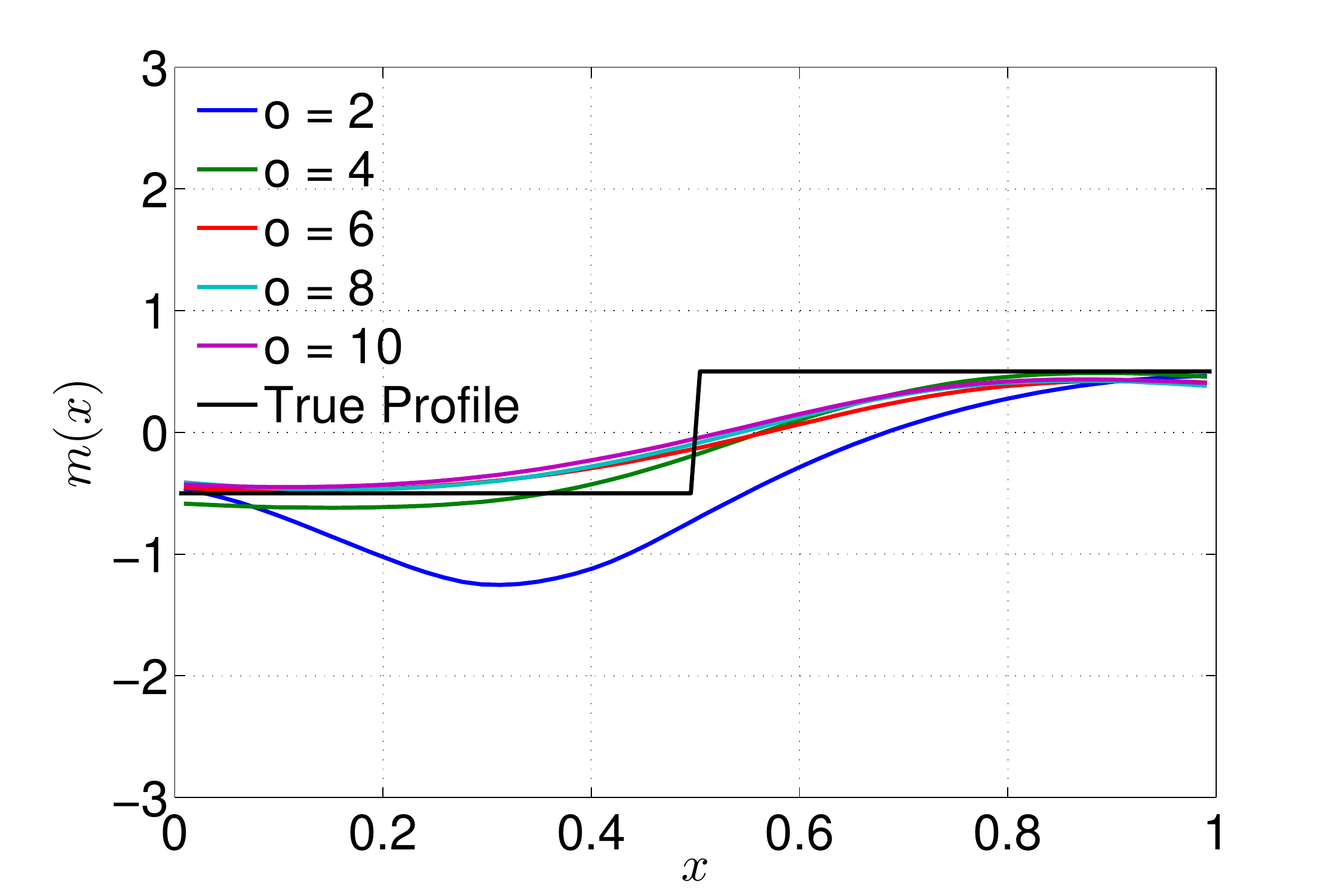}  &
\includegraphics[width=0.35\textwidth]{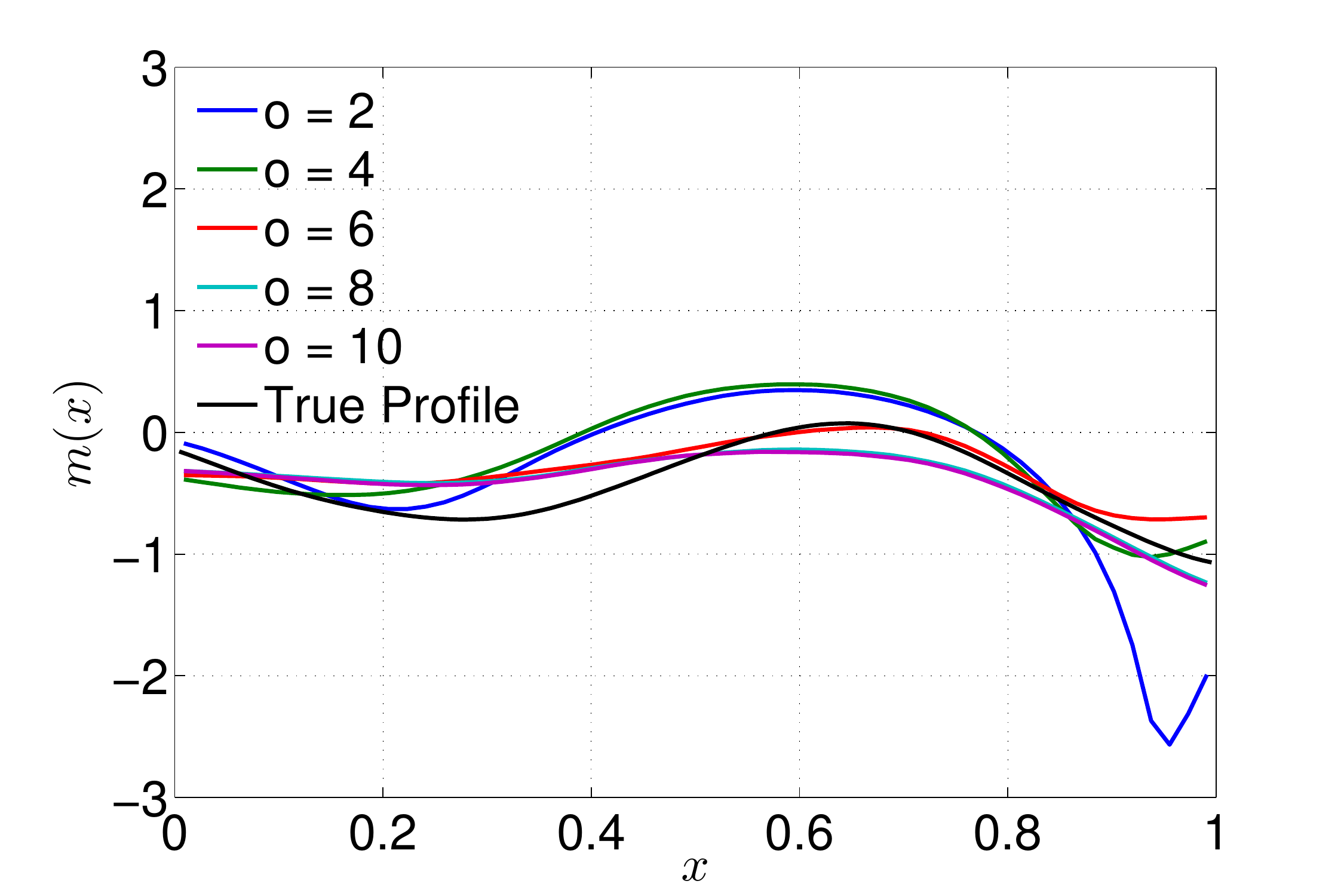} 
\end{tabular}
\caption{Effect of PC order $o$ on the inferred median of the posterior: cases of (Left) $m^{\rm sin}$, (Center) $m^{\rm step}$ and (Right) 
$m^{\rm ran}$.}
\label{fig:reffect}
\end{figure}

Finally, we illustrate in Figure~\ref{fig:reffect} the dependence of the inferred median profiles on the selected order $o$ for the PC surrogate model. The figure shows that, irrespective to the smoothness of the true profile, the inferred medians quickly converge as $o$ increases, demonstrating that the $L_2$ convergence of the PC surrogate with coordinate transformation reported in Section~\ref{pc:example} transfers to the inference problem. In fact, in view of the convergence curves shown in Figure~\ref{fig:errorU}, the differences in the inferred median profiles for $o=8$ and $o=10$ are more likely to come from sampling errors than from differences in the PC surrogates. 

%%%%%%%%%%%%%%%%%%%%%%%%%%%%%%%%%%%%%%%%%%%%%
%%%%%%%%%%%%%%%%%%%%%%%%%%%%%%%%%%%%%%%%%%%%%%%%%%
%%%%%%%%%%%%%%%%%%%%%%%%%%%%%%%%%%%%%%%%%%%%%%%%%
%!TEX root = paper.tex
\section{Discussion and Conclusion}
\label{sec:conc}

This paper presented a Bayesian approach to infer a parameter field from prior ${\cal GP}$ having a covariance function involving some hyper-parameters $\vec q$. 
The main contribution of the present work is the introduction of a coordinate transformation in order to represent the prior ${\cal GP}$ using a unique reference basis of spatial modes, while the effects of the covariance hyper-parameters is reflected by the (joint) prior probability density function of the random coordinates of the ${\cal GP}$ that becomes conditioned on $\vec q$. 
The coordinate transformation naturally leads to the construction of a unique polynomial surrogate for the forward model predictions; this surrogate model accounts for the dependence of the model predictions on the coordinates of the ${\cal GP}$ in the reference basis. For a Polynomial Chaos approximation, as considered in this paper, the construction of the surrogate amounts to solving a unique (stochastic) reference problem, assuming the independence of the ${\cal GP}$ coordinates. The stochastic dimensionality of the surrogate model is therefore equal to the dimensionality of the (truncated) ${\cal GP}$ representation, and is not augmented by the number of hyper-parameters intervening in the covariance function parametrization. This fact has to be contrasted with the alternative approaches proposed in~\citep{MarzoukNajm2009,Tagade:2014} where the PC expansion explicitly incorporates the dependencies on $\vec q$. 
Another advantage of selecting a reference problem for the construction of the PC surrogate, compared to the direct expansion with respect to the covariance hyper-parameters, is that it can overcome issues related to hyper-parameters with complex distributions, \textit{e.g.} improper, no second-order moments, \dots\ for which classical PC bases may not exist. 

The surrogate model can then be substituted for the true model predictions in the definition of the likelihood of the observations appearing in Bayes' formula for the posterior of the ${\cal GP}$ coordinates and covariance hyper-parameters. 
The resulting approximate likelihood can in turn be imbedded in a MCMC sampler to greatly accelerate the sampling of the posterior distribution, with significant computational savings. In its present form, the proposed method however introduces some overhead during the sampling stage, compared to other approaches relying on PC acceleration with explicit dependence on $\vec q$: for any new proposed values of the hyper-parameters the coordinate transformation must be determined. The determination of the transformation, given $\vec q$, requiring the computation of the dominant subspace of the covariance function (given $\vec q$) may constitute a severe limitation for large scale problems (for the simplified problems presented in Section~\ref{sec:results}, the CPU time of the inference with hyper-parameters was found roughly three time as large as for the case without hyper-parameters). To remedy this point in the future, we plan to approximate the dependence of the coordinate transformation, $\hat{\cal B}$, on $\vec q$ using, again, a PC expansion. As for the construction of the PC surrogate of the model predictions, the approximation of the coordinate transformation will be computed off-line and subsequently used in-line within the sampler.

The numerical experiments presented in the paper, although based on a simple model, have highlighted the following points:
\begin{itemize}
	\item Using for reference basis the truncated set of dominant modes of the $\vec q$-averaged covariance function is not only optimal (on average) for the representation of the processes with variable $\vec q$, but it also appears as the best choice in terms of averaged error for the PC surrogate of the model prediction in our example.
	\item The control of the stretching induced by the coordinate mapping $\hat{\cal B}(\vec q)$ is crucial for the error control;
	while using the marginalized conditional density $p_{\bar\eta}(\bar{\vec \eta}|\vec q)$ appears to be an appropriate choice, other alternatives may be conceived. In particular, augmenting the variability of the reference process $M^{\rm PC}_K(\xxi)$ could improve the robustness of the surrogate PC model.
	\item The introduction of covariance functions with hyper-parameters clearly improved the inference results in the problems considered, particularly when inferring smooth profiles. In particular, information gain was observed for a larger set of coordinates. In addition, when covariance hyper-parameters was accounted for, the convergence rate of the inferred field with increased number of observations and reduced observation noise also seemed to improve. 
	\item The convergence with the PC surrogate order seems quite fast for the presented problems, suggesting to possibility of using moderate PC orders, particularly to balance PC error and posterior sampling errors.
\end{itemize}

On the basis of the present findings, we plan for future work to develop the coordinate transformation approach to further exploit the posterior structure involving the conditional prior probability of the transformed coordinates $\hat{\vec \eta}$ and derive samplers adapted to this particular structure. 
Regarding the construction of the PC surrogate model, consideration of adaptive constructions would be beneficial to reduce the computational cost of the off-line step, to increase accuracy, and further accelerate the sampler. 
Further, the PC approximation of the coordinate transformation appears to be a key element to make the whole approach effective to handle large scale problems. 
In addition, the proposed method, in particular the construction of the PC approximation of the model prediction, would certainly benefit from fitting the procedure to the posterior distributions (of coordinates and hyper-parameters) rather than to the prior ones, especially when the observations are informative. Since these posterior distributions are not known \textit{a priori}, iterative constructions are needed.
Pursuit of these avenues is currently considered on a complex problem arising in subsurface geological models and earthquake model.

%%%%%%%%%%%%%%%%%%%%%%%%%%%%%%
%%%%%%%%%%%%%%%%%%%%%%%%%%%%%%%%%%
%%%%%%%%%%%%%%%%%%%%%%%%%%%%%%%%%%%%%%

\section*{Acknowledgments}
Research reported in this publication was supported by the King Abdullah University of Science and Technology (KAUST).  OLM and OK also acknowledge partial support provided the US Department of Energy (DOE), Office of Science, Office of Advanced Scientific Computing Research, under Award Number DE-SC0008789. 
%\nomenclature{$a$}{}%
%\appendix

\bibliographystyle{elsarticle-num}

%\bibliography{biblio}
\end{document}